\documentclass[reqno,11pt]{amsart}

\usepackage{amsmath,amssymb,bbm,dsfont,url}
\usepackage{mdframed} 
\usepackage[margin=1.03in]{geometry}
\usepackage[english]{babel}
\usepackage[latin1]{inputenc}   
\usepackage{hyperref}
\usepackage{amssymb}
\usepackage{relsize}
\usepackage{mathtools}
\usepackage[title]{appendix}
\usepackage{soul}
\usepackage[all,cmtip]{xy}
\usepackage[normalem]{ulem} 
\usepackage{enumerate}
\usepackage{tikz-cd}
\usetikzlibrary{matrix,arrows}
\usepackage{multirow}
\usepackage{xcolor}
\usepackage{etoolbox}
\usepackage{paracol}
\usetikzlibrary{calc}
\usetikzlibrary{decorations.pathmorphing}
\tikzset{curve/.style={settings={#1},to path={(\tikztostart)
    .. controls ($(\tikztostart)!\pv{pos}!(\tikztotarget)!\pv{height}!270:(\tikztotarget)$)
    and ($(\tikztostart)!1-\pv{pos}!(\tikztotarget)!\pv{height}!270:(\tikztotarget)$)
    .. (\tikztotarget)\tikztonodes}},
    settings/.code={\tikzset{quiver/.cd,#1}
        \def\pv##1{\pgfkeysvalueof{/tikz/quiver/##1}}},
    quiver/.cd,pos/.initial=0.35,height/.initial=0}
\tikzset{tail reversed/.code={\pgfsetarrowsstart{tikzcd to}}}
\tikzset{2tail/.code={\pgfsetarrowsstart{Implies[reversed]}}}
\tikzset{2tail reversed/.code={\pgfsetarrowsstart{Implies}}}
\tikzset{no body/.style={/tikz/dash pattern=on 0 off 1mm}}
\usepackage{comment}
\hypersetup{
  colorlinks   = true, 
  urlcolor     = blue, 
  linkcolor    = blue, 
  citecolor   = black 
}

\usepackage{amsthm}
\newtheorem{counter}{}[section]
\theoremstyle{definition}
\newtheorem{definition}         [counter]{Definition}

\theoremstyle{plain}
\newtheorem{lemma}              [counter]{Lemma}

\newtheorem{proposition}        [counter]{Proposition}
\newtheorem{theorem}            [counter]{Theorem}
\newtheorem{corollary}          [counter]{Corollary}

\newtheorem*{theorem*}          {Theorem}

\newtheoremstyle{named}{8pt}{8pt}{\itshape}{}{\bfseries}{.}{.5em}{#1 #3}
\theoremstyle{named}
\newtheorem{namedtheorem}{Theorem}

\makeatletter
\newcommand{\setcustomlabel}[1]{\def\@currentlabel{#1}}
\makeatother

\theoremstyle{remark}
\newtheorem{question}           [counter]{Question}
\newtheorem*{question*}          {Question}
\newtheorem{remark}             [counter]{Remark}
\newtheorem*{remark*}           {Remark}
\newtheorem{example}            [counter]{Example}

\numberwithin{equation}{section}

\usepackage{setspace}
\usepackage{graphicx}
\author{Kenichi Shimizu}
\address{Department of Mathematical Sciences, Shibaura Institute of Technology, 307 Fukasaku, Minuma-ku, Saitama-shi, Saitama 337-8570, Japan}
\email{kshimizu@shibaura-it.ac.jp}

\author{Harshit Yadav}
\address{Department of Mathematics, University of Alberta, Edmonton, AB, Canada T6G 2G1}
\email{hyadav3@ualberta.ca}

\makeatletter
\@namedef{subjclassname@2020}{%
  \textup{2020} Mathematics Subject Classification}
\makeatother

\allowdisplaybreaks

\newcommand{\YD}{\mathcal{YD}}

\DeclareMathOperator{\End}{End}

\newcommand{\cA}{\mathcal{A}}
\newcommand{\C}{\mathcal{C}}
\newcommand{\B}{\mathcal{B}}
\newcommand{\D}{\mathcal{D}}
\newcommand{\E}{\mathcal{E}}

\newcommand{\Z}{\mathcal{Z}}
\newcommand{\W}{\mathcal{W}}

\newcommand{\M}{\mathcal{M}}

\newcommand{\cO}{\mathcal{O}}

\newcommand{\bZ}{\mathbb{Z}}

\newcommand{\bA}{\mathbf{A}}

\newcommand{\bq}{\mathbf{q}}

\newcommand{\fp}{\mathfrak{p}}

\newcommand{\fu}{\mathfrak{u}}

\newcommand{\fg}{\mathfrak{g}}

\newcommand{\fsl}{\mathfrak{sl}}
\newcommand{\fB}{\mathfrak{B}}

\newcommand{\kk}{\Bbbk}

\newcommand{\ev}{\mathrm{ev}}
\newcommand{\coev}{\mathrm{coev}}

\newcommand{\id}{\mathrm{id}}
\newcommand{\Rep}{\mathrm{Rep}}

\newcommand{\Hom}{\mathrm{Hom}}
\newcommand{\ns}{\mathrm{ns}}
\newcommand{\Ind}{\mathrm{Ind}}

\newcommand{\unit}{\mathbbm{1}}

\newcommand{\Vect}{\mathsf{Vec}}
\newcommand{\uHom}{\underline{\sf Hom}}

\newcommand{\ra}{{\sf ra}}

\newcommand{\rra}{{\sf rra}}

\newcommand{\FPdim}{{\sf FPdim}}

\newcommand{\Inv}{{\sf Inv}}

\newcommand{\ucoev}{\underline{\mathrm{coev}}}

\newcommand{\loc}{\mathrm{loc}}

\newcommand{\op}{\mathrm{op}}

\newcommand{\oC}{\overline{\mathcal{\C}}}
\newcommand{\oD}{\overline{\mathcal{\D}}}
\newcommand{\oE}{\overline{\mathcal{\E}}}

\definecolor{darkgreen}{RGB}{55,138,0}
\definecolor{burntorange}{RGB}{180,85,0}
\definecolor{navyblue}{RGB}{18,40,180}
\definecolor{cyan(process)}{rgb}{0.0, 0.6, 1.0}

\begin{document}

\title{Commutative exact algebras and modular tensor categories}

\begin{abstract}
Inspired by the study of vertex operator algebra extensions, we answer the question of when the category of local modules over a commutative exact algebra in a braided finite tensor category is a (non-semisimple) modular tensor category. Along the way we provide sufficient conditions for the category of local modules to be rigid, pivotal and ribbon. We also discuss two ways to construct such commutative exact algebras. The first is the class of simple current algebras and the second is using right adjoints of central tensor functors. Furthermore, we discuss Witt equivalence and its relation with extensions of VOAs.
\end{abstract}
 
\maketitle


\section{Introduction}
Two-dimensional conformal field theory ($2$D-CFT) is the richest and most well-studied quantum field theory. Mathematically, a $2$D-CFT (its chiral half, to be precise) is understood using its symmetry algebra, which is a vertex operator algebra (VOA). VOAs $V$ are, in turn, understood using their category of representations $\Rep(V)$. 
The physical features of the CFT manifest into various properties of the category $\Rep(V)$ like it being braided monoidal, rigid, ribbon, nondegenerate, etc. This leads to the question of when a VOA admits a suitable category of modules which has the aforementioned features. 
As this is a difficult question in general, one looks at it through the lens of various constructions that relate VOAs. An important technique of understanding VOAs is by VOA extensions, which is just an inclusion $V\subset W$ of VOAs with the same conformal vector. In this case, if $\C=\Rep(V)$ is a braided monoidal category then $A:=W$ becomes a commutative algebra object in $\C$. Additionally, the category of $W$-modules, which we denote as $\D=\Rep(W)$, is a braided monoidal category that is braided tensor equivalent to the category $\C_A^{\loc}$ of so-called local $A$-modules in $\C$, namely $\D\simeq\C_A^{\loc}$ \cite{kirillov2002q,huang2015braided,creutzig2017tensor}.
Thus to understand $\Rep(W)$ using $\Rep(V)$ one needs to answer the following question.

\begin{question*}
 Let $\C$ be a braided monoidal category and $A$ a commutative algebra in $\C$. Under what conditions on $\C$ and $A$ is the category $\C_A^{\loc}$ braided, rigid, pivotal, ribbon, modular, etc.?
\end{question*}

Work on this question goes back to Pareigis \cite{pareigis1995braiding} who introduced the idea of local modules and proved that the category $\C_A^{\loc}$ is braided monoidal. Kirillov and Ostrik \cite{kirillov2002q} specialized to the case of semisimple categories and gave a complete answer to the question above using separable algebras satisfying additional conditions. Recent works that have generalized this include \cite{creutzig2017tensor} and \cite{laugwitz2023constructing}. While these works generalize the categories involved, the former by considering monoidal supercategories and the latter by considering non-semisimple tensor categories, the conditions on the algebra $A$ remained restrictive. This is a limitation because, in the non-semisimple setting, one finds many interesting algebra objects that are not separable. In this article, we work at the level of abelian rigid monoidal categories to answer the question. Upon specializing to finite tensor categories, we obtain new constructions of modular tensor categories (MTCs) via local modules.


\subsection*{Main results}
The background material for the following discussion can be found in Section~\ref{sec:preliminaries}. In the following, we proceed step by step to explain how we establish that $\C_A^{\loc}$ is rigid and ribbon, and $\C_A$ is pivotal, using various properties of $\C$ and $A$.

\underline{Rigidity}:
In view of future applications, we start in Section~\ref{sec:closed} with a closed monoidal category (one whose tensor product admits a right adjoint, called the internal Hom functor). This is weaker than assuming that $\mathcal{C}$ is rigid. If $\mathcal{C}$ is closed monoidal with equalizers and coequalizers, and $A$ is an algebra in $\mathcal{C}$, then the category of $A$-bimodules ${}_A\mathcal{C}_A$ is also a closed monoidal category.

If $A$ is a commutative algebra in $\mathcal{Z}(\mathcal{C})$ with half-braiding $\sigma$, one can form the category $\mathcal{C}_A^{\sigma}$ \cite{schauenburg2001monoidal}, the category $\mathcal{C}_A$ of right $A$-modules with a monoidal structure defined using $\sigma$. We show that $\mathcal{C}_A^{\sigma}$ is a closed monoidal subcategory of ${}_A\C_A$ (in the sense that it is closed under the internal Hom functor of ${}_A\C_A$). When $\C$ is braided and $A$ is a commutative algebra in $\mathcal{C}$, this implies that $\mathcal{C}_A$ and $\mathcal{C}_A^{\mathrm{loc}}$ are also closed monoidal subcategories of ${}_A\mathcal{C}_A$.

From this, we derive the first main result of Section~\ref{sec:closed} regarding rigidity: if ${}_A\mathcal{C}_A$ is rigid, the subcategories $\mathcal{C}_A$ and $\mathcal{C}_A^{\mathrm{loc}}$ are also rigid. Note that by the work of Pareigis \cite{pareigis1995braiding}, we already know that $\mathcal{C}_A^{\mathrm{loc}}$ is also a braided monoidal category.

\underline{Ribbon structure}: If $\C$ is a braided category with twist $\theta$ and $A$ is a commutative algebra with $\theta_A=\id_A$, it is easy to see that $\theta$ yields a twist for $\C_A^{\loc}$ as well. We show that if $\C$ is ribbon with $\theta_A=\id_A$, then $A^*\in\C_A^{\loc}$. To get that $\C_A^{\loc}$ is ribbon, one needs to show that the twist respects duality. We prove that this is equivalent to the requirement that $A^*$ belongs to $(\C_A^{\loc})'$, that is, the M\"uger center of $\C_A^{\loc}$.

\underline{Pivotal structure}:
Next, in Section~\ref{sec:double-dual-rigid}, we move to discussing when $\C_A^{\sigma}$ is pivotal. This then yields pivotality results for $\C_A$ and $\C_A^{\loc}$. We suppose that $\C$ is rigid and $(A,\sigma)$ is a commutative algebra in $\Z(\C)$ such that ${}_A\C_A$ is rigid. Let $(-)^{\dagger}$ denote the duality functor of ${}_A\C_A$. Unlike prior works \cite{kirillov2002q,laugwitz2023constructing}, it can happen that, for $M\in\C_A^{\sigma}$, $M^{\dagger}$ is not isomorphic to $M^*$ as an object of $\C$. This makes an analysis of duality difficult. However, to understand pivotality, we need to understand the double dual functor $(-)^{\dagger\dagger}$. We show that $(M)^{\dagger\dagger} \cong M^{**}$ as an object but the $A$-action is twisted using $u := \fu_{(A,\sigma)}$, where $\fu$ is the Drinfeld isomorphism (see \S\ref{subsec:prelim-braided-monoidal} for the definition).
We write $(M^{**})_{(u)}$ to denote this $A$-module. 
Then, we prove the isomorphism $(-)^{\dagger\dagger}\cong (-)^{**}_{(u)}$ of monoidal functors.
Using this, one can conclude that if $\C$ is pivotal with pivotal structure $\fp$, $(A,\sigma)$ is a commutative algebra in $\Z(\C)$ such that ${}_A\C_A$ is rigid and $u = \fp_A$, then the rigid monoidal category $\C_A^{\sigma}$ is pivotal. This immediately implies that if $A$ is a commutative algebra in $\C$ satisfying $\fu_A = \fp_A$, then $\C_A$ and $\C_A^{\loc}$ are pivotal.

The above discussion can be summarized in the following result:
\begin{namedtheorem}[A]\setcustomlabel{A}\label{thm:intro-A}
Let $\C$ be a rigid braided monoidal category with equalizers and coequalizers. Let $A$ be a commutative algebra in $\C$. Then,
\begin{center}
  {\small
  \begin{tabular}{ c c c }
    ${}_A\C_A$ rigid & $\implies$ \; & $\C_A$ and $\C_A^{\loc}$ are rigid, \\
    $\C$ pivotal, ${}_A\C_A$ rigid and $\fu_A=\fp_A$ & $\implies$ \; & $\C_A$ and $\C_A^{\loc}$ are pivotal, \\
    $\C$ ribbon, ${}_A\C_A$ rigid, $\theta_A=\id$ and $A^*\in(\C_A^{\loc})'$ & $\implies$ \; & $\C_A^{\loc}$ is ribbon.
  \end{tabular}
  }
  \end{center}
\end{namedtheorem}

We call a non-degenerate ribbon finite tensor category a modular tensor category (MTC). MTCs are the algebraic input needed to construct $3$D topological quantum field theories ($3$D-TQFTs) \cite{de20223}. $3$D-TQFTs in turn yield invariants of $3$-manifolds \cite{reshetikhin1991invariants,kerler2001non,de20223} and mapping class group representations \cite{de2023mapping}. This inspired us to study the question of when $\C_A^{\loc}$ is an MTC. For this, we restrict to the setting of $\C$ being a finite tensor category.

When $\C$ is a finite tensor category (FTC), ${}_A\C_A$ is rigid if and only if $A$ is an exact algebra \cite{etingof2004finite}. We collect some basic results about commutative exact algebras and their category of local modules in Section~\ref{sec:finite-cats}. In the finite setting, if $\C$ is non-degenerate, we can also establish the non-degeneracy of $\C_A^{\loc}$, generalizing \cite{laugwitz2023constructing}. Furthermore, if $\C$ is an MTC and $A$ is a haploid, commutative and exact algebra, then the conditions $\theta_A=\id$ and $A^*\in(\C_A^{\loc})'$ are equivalent to $A$ being a symmetric Frobenius algebra.
Thus, we obtain the following result (see Section~\ref{subsec:MTC-main-result} for details). 

\begin{namedtheorem}[B]\setcustomlabel{B}\label{thm:intro-B}
Let $\C$ be a braided FTC and $A$ a haploid and commutative algebra in $\C$. Then,
\begin{center}
{\small
\begin{tabular}{ c c c }
$A$ exact & $\implies$ \; & $\C_A$ and $\C_A^{\loc}$ are FTCs, \\
$\C$ pivotal, $A$ exact, $\fu_A=\fp_A$ & $\implies$ \; & $\C_A$ and $\C_A^{\loc}$ are pivotal, \\
$\C$ ribbon, $A$ exact, $\theta_A=\id$, $A^*\in(\C_A^{\loc})'$ & $\implies$ \; & $\C_A^{\loc}$ is ribbon, \\
$\C$ MTC, $A$ exact symmetric Frobenius & $\implies$ \; & $\C_A^{\loc}$ is an MTC.
\end{tabular}
}
\end{center}
If $\C$ is an MTC and $A$ is a commutative, simple, symmetric Frobenius algebra, then so is $\C_A^{\loc}$.
\end{namedtheorem}

After our article appeared, in \cite{coulembier2025simple}, it was proven that simple algebras are exact. Using this, the above Theorem can be rephrased in terms of simple algebras (see Remark~\ref{rem:CSZ-simple} for details).

Previous works \cite{kirillov2002q,davydov2013witt} established that $\mathcal{C}_A^{\text{loc}}$ is a semisimple MTC when $\mathcal{C}$ is semisimple and $A$ is separable with trivial twist. This was extended to non-semisimple MTCs in \cite[Theorem~4.16]{laugwitz2023constructing} using haploid, commutative, special Frobenius algebras. While these algebras satisfy the conditions of Theorem~\ref{thm:intro-B}, we demonstrate that our result generalizes \cite{laugwitz2023constructing} by providing a family of algebras that satisfy our requirements but not those of the latter. Our results in Theorem~\ref{thm:intro-B} can also be specialized to monoidal supercategories and this is done in \cite{mcrae2025rigidity}.
See Section~\ref{subsec:relation-to-recent-work} for further details.

With the above result in hand, one needs examples of commutative exact algebras that satisfy the conditions of the theorem. In Section~\ref{sec:commutative-algebras}, we discuss two ways to construct such algebras.
The first is the class of simple current algebras, these are commutative algebras which are direct sums of simple invertible objects \cite{fuchs2004tft,creutzig2020simple}. We discuss this in Section~\ref{subsec:simple-currents} where we provide examples of such algebras in the category of representations of the (even order) small quantum group and the Drinfeld double of the Taft algebra. The second construction is using right adjoints $F^{\ra}$ of central tensor functors $F:\C\rightarrow\D$ between FTCs. This known construction \cite{davydov2013witt} produces a haploid, commutative and exact algebra, namely $A=F^{\ra}(\unit)$ in $\C$.
We discuss this in Section~\ref{subsec:comm-alg-via-adjoint} where we use this to construct an example of a commutative exact algebra in the category of representations of the (odd order) small quantum group.

Lastly in Sections~\ref{sec:Witt-equivalence} we use commutative exact algebras to discuss Witt equivalence in the non-semisimple setting. The use of such algebras, instead of separable commutative algebras as in \cite{laugwitz2023constructing}, allows us to generalize many results from the semisimple case \cite{davydov2013witt}. 
We also elaborate on the connection between extensions of logarithmic VOAs and Witt equivalence. We call a braided FTC \textit{completely anisotropic} if the only haploid, commutative and exact algebra in it is $\unit$, extending the notions defined in \cite{davydov2013structure,laugwitz2023constructing}.
For two examples of VOAs, the triplet VOA and even part of the symplectic fermion VOAs, we show that their representation categories are not completely anisotropic. Since this paper was posted, a non-semisimple completely anisotropic braided FTC over characteristic $0$ has been constructed in \cite{ostrik2026non}. This and other questions related to Witt groups are discussed in Section~\ref{subsec:Witt-questions}.


\subsection*{Outlook}
Our results here have interesting implications for extensions $V \subset W$ of VOAs, which we will briefly discuss now (refer to Section~\ref{subsec:VOA-extensions} for additional background). Among the most well-understood classes of VOAs are strongly rational VOAs. Given a strongly rational VOA $V$, $\Rep(V)$ forms a semisimple MTC \cite{huang2008rigidity}. Thus, establishing the strong rationality of VOAs is desirable. Using our rigidity results, we demonstrate in \cite{mcrae2025rigidity} that a simple extension $W$ of a strongly rational VOA $V$ is strongly rational. 

A natural generalization of strongly rational VOAs are strongly finite VOAs. These VOAs $V$ have a category of modules conjectured to be an MTC \cite{huang2009representations}. In this context, McRae \cite{mcrae2021rationality} proved that if $\Rep(V)$ is rigid, then it is an MTC. Therefore, establishing the rigidity of strongly finite VOAs is crucial. In \cite{mcrae2025rigidity}, we prove that if $V \subset W$ is an extension of strongly finite VOAs where $\Rep(V)$ is an MTC, then $\Rep(W)$ is rigid and thus an MTC. For more details, we refer the reader to \cite{mcrae2025rigidity}.

Many currently known examples of VOAs admit a rigid braided tensor category of representations that is not necessarily finite. Given an extension $V \subset W$ of such a VOA, the algebra $A := W$ is usually an object in the ind-completion of $\Rep(V)$. With this in mind, in an upcoming work, we will address the questions of being rigid, ribbon, non-degenerate, etc., for a certain subcategory of finitely generated objects of $\Ind(\mathcal{C})_A^{\mathrm{loc}}$.


\section{Preliminaries}\label{sec:preliminaries} 
Given a category $\C$, we write $\C^{\op}$ to denote the opposite category of $\C$. We will work over a (not necessarily algebraically closed) field $\kk$ and denote by $\Vect$ the category of finite-dimensional $\kk$-vector spaces.


\subsection{Monoidal categories}\label{subsec:monoidal}
As in \cite[\S2.1]{etingof2016tensor}, we will denote a monoidal category by a triple $(\C,\otimes,\unit)$ where $\C$ is a category, $\otimes$ is a bifunctor $\otimes:\C\times\C\rightarrow\C$ and $\unit$ is an object of $\C$.


\subsubsection*{Rigid monoidal categories}
Let $(\C,\otimes,\unit)$ be a monoidal category \cite[\S2.1]{etingof2016tensor}. 
Given an object $X\in\C$, a \textit{left dual} of $X$ is an object $X^*\in\C$ such that there exist maps $\ev_X:X^*\otimes X\rightarrow\unit$ and $\coev_X:\unit\rightarrow X\otimes X^*$ satisfying 
\[(\id\otimes \ev_X)(\coev_X\otimes\id) = \id_X, \quad
(\ev_X\otimes \id)(\id\otimes\coev_X) = \id_{X^*}.\] 
Similarly, a \textit{right dual} of $X$, denoted ${}^*X$, is defined \cite[\S2.10]{etingof2016tensor}. 
An object $X$ is called \textit{rigid} if it admits a left and a right dual.
A monoidal category is called \textit{rigid} if every object is rigid. 

A rigid monoidal category is called \textit{pivotal}, if it comes equipped with a monoidal natural isomorphism $\fp: \id_{\C}\Rightarrow (-)^{**}$ \cite[\S4.7]{etingof2016tensor}.
For an endomorphism $f: X \to X$ in a pivotal monoidal category with pivotal structure $\fp$, we define the pivotal trace of $f$ by
\begin{equation*}
  \mathrm{tr}(f) = \ev_{X^*} \circ (\fp_X f \otimes \id_{X^*}) \circ \coev_X : \unit \to \unit.
\end{equation*}
The pivotal dimension of $X \in \C$ is defined to be the pivotal trace of $\id_X$. When $\C$ is $\kk$-linear and $\End_{\C}(\unit) \cong \kk$, we often regard $\mathrm{tr}(f)$ as an element of $\kk$. 


\subsubsection*{Braided monoidal categories}
\label{subsec:prelim-braided-monoidal}
We refer the reader to \cite[\S8]{etingof2016tensor} for basic definitions related to braided monoidal categories.
Given a braided monoidal category $\B$ with braiding $c$, we denote by $\overline{\B}$ the monoidal category $\B$ with the braiding $\overline{c}$ given by $\overline{c}_{X,Y} = c_{Y,X}^{-1}$ ($X, Y \in \C$).

For a braided monoidal category $\C$ with braiding $c$ and its full subcategory $\D$, the \textit{centralizer} of $\D$ in $\C$ is the full subcategory of $\C$ defined as:
\[ \Z_{(2)}(\D\subset \C) = \{X\in\C \mid \text{$c_{Y,X}\circ c_{X,Y} =\id_{X\otimes Y}$ for all $Y\in\D$} \}. \]
The category $\Z_{(2)}(\C\subset\C)$ is written as $\Z_{(2)}(\C)$ or, more shortly, as $\C'$.

We assume that the braided monoidal category $\B$ is rigid. Then the \textit{Drinfeld isomorphism} of $\B$ is the natural isomorphism $\mathfrak{u}:\id_{\B}\rightarrow (-)^{**}$ given by
\begin{equation}
  \label{eq:def-Drinfeld-iso}
    \mathfrak{u}_X=(X \xrightarrow{\id\otimes\coev_{X^*}} X\otimes X^*\otimes X^{**} \xrightarrow{c\otimes \id} X^*\otimes X\otimes X^{**} \xrightarrow{\ev_X\otimes \id} X^{**} ) \quad (X\in\B).
\end{equation}
A \textit{twist} on $\B$ is a natural isomorphism $\theta:\id_{\B}\rightarrow \id_{\B}$ satisfying 
\begin{equation}
  \label{eq:def-twist}
  \theta_{X\otimes Y} = (\theta_X\otimes \theta_Y)\circ c_{Y,X}\circ c_{X,Y}
\end{equation} for all $X, Y \in \B$.
If further, $(\theta_X)^*=\theta_{X^*}$ holds for all $X\in\B$, then $\theta$ is called a \textit{ribbon structure}. A braided category equipped with a ribbon structure is called a \textit{ribbon category} \cite[\S 8.10]{etingof2016tensor}. 

Let $\C$ be a monoidal category.
Given an object $X\in\C$, a \textit{lax half-braiding} on $X$ is a natural transformation 
$\sigma = \{ \sigma_Y: X\otimes Y\rightarrow Y \otimes X\}_{Y\in\C}$
satisfying $\sigma_{\unit} = \id_{X}$ and
\begin{equation}
  \label{eq:lax-half-braiding}
  \sigma_{Y \otimes Z} = (\id_Y \otimes \sigma_Z) \circ (\sigma_Y \otimes \id_Z)
\end{equation}
for all $Y, Z \in \C$.
A \textit{half-braiding} is an invertible lax half-braiding.
The (lax) \textit{Drinfeld center} \cite[\S7.13]{etingof2016tensor} is the category whose objects are pairs $(X, \sigma)$ consisting of an object $X \in \C$ and a (lax) half-braiding $\sigma$ and whose morphisms are morphisms in $\C$ that are compatible with the (lax) half-braidings. The Drinfeld center of $\C$, denoted by $\Z(\C)$, is a braided monoidal category.

Assume that $\C$ is rigid. Then the lax Drinfeld center of $\C$ coincides with $\Z(\C)$. Indeed, any lax half-braiding $\sigma$ on $X \in \C$ is invertible with the inverse
\begin{equation}
  \label{eq:inverse-of-half-br}
  \sigma_Y^{-1} = (\ev_{{}^*Y} \otimes \id_X \otimes \id_Y)
  (\id_Y \otimes \sigma_{{}^*Y} \otimes \id_Y)
  (\id_Y \otimes \id_X \otimes \coev_{{}^*Y})
  \quad (Y \in \C).
\end{equation}
Moreover, $\Z(\C)$ is also a rigid monoidal category \cite[Section 7.13]{etingof2016tensor}.
For $\mathbf{X} = (X, \sigma) \in \Z(\C)$, the evaluation and the coevaluation are given by $\ev_{\mathbf{X}} = \ev_X$ and $\coev_{\mathbf{X}} = \coev_X$, respectively. 

\subsubsection*{Graphical calculus}
We sometimes use the techniques of graphical calculus to avoid extremely long expressions. Our convention is that the evaluation, the coevaluation and the braiding and its inverse are depicted as follows:
\begin{equation*}
    \tikzset{overstrand/.style={preaction={-,draw=white,line width=8pt}}}
    \ev_X = \begin{tikzpicture}[x = 10pt, y = 10pt, baseline = 0]
      \coordinate (S1) at (0,1); \node at (S1) [above] {$X^*$};
      \coordinate (S2) at (2,1); \node at (S2) [above] {$X$};
      \draw (S1) to [out = -90, in = -90, looseness=4] (S2);
    \end{tikzpicture} \qquad
    \coev_X = \begin{tikzpicture}[x = 10pt, y = 10pt, baseline = 0]
      \coordinate (S1) at (0,-1); \node at (S1) [below] {$X$};
      \coordinate (S2) at (2,-1); \node at (S2) [below] {$X^*$};
      \draw (S1) to [out = 90, in = 90, looseness=4] (S2);
    \end{tikzpicture} \qquad
    c_{X,Y} = \begin{tikzpicture}[x = 8pt, y = 8pt, baseline = 0]
      \coordinate (S1) at (0,2); \node at (S1) [above] {$X$};
      \coordinate (S2) at (3,2); \node at (S2) [above] {$Y$};
      \coordinate (T1) at (0,-2); \node at (T1) [below] {$Y$};
      \coordinate (T2) at (3,-2); \node at (T2) [below] {$X$};
      \draw (S1) to [out = -90, in = 90] (T2);
      \draw [overstrand] (S2) to [out = -90, in = 90] (T1);
    \end{tikzpicture} \qquad
    c_{X,Y}^{-1} = \begin{tikzpicture}[x = 8pt, y = 8pt, baseline = 0]
      \coordinate (S1) at (0,2); \node at (S1) [above] {$Y$};
      \coordinate (S2) at (3,2); \node at (S2) [above] {$X$};
      \coordinate (T1) at (0,-2); \node at (T1) [below] {$X$};
      \coordinate (T2) at (3,-2); \node at (T2) [below] {$Y$};
      \draw (S2) to [out = -90, in = 90] (T1);
      \draw [overstrand] (S1) to [out = -90, in = 90] (T2);
    \end{tikzpicture} .
\end{equation*}


\subsection{Tensor categories}
By a \textit{multi-tensor category}, we will mean a locally finite, $\kk$-linear, abelian, rigid monoidal category $\C$ with $\otimes:\C\times\C\rightarrow\C$ bilinear on morphisms. Additionally, if $\End_{\C}(\unit)\simeq \kk$, we call $\C$ a \textit{tensor category}. If in addition, $\C$ is a finite abelian category (in the sense of \cite[Definition 1.8.6]{etingof2016tensor}), we call it a \textit{finite tensor category} (FTC).

A \textit{tensor functor} is a $\kk$-linear, exact, strong monoidal functor between tensor categories. By \cite[Subsection 2.2]{bruguieres2011exact}, a tensor functor is faithful.
The \textit{image} of a functor $T: \mathcal{A} \to \mathcal{B}$ between abelian categories, denoted by $\mathrm{Im}(T)$, is the full subcategory of $\mathcal{B}$ consisting of all subquotients of objects of the form $F(X)$ for some $X \in \mathcal{A}$.
The image of a tensor functor between tensor categories is a \textit{tensor full subcategory} in the sense that it is closed under subquotients, duals and tensor product.
We say that a tensor functor $F:\C\rightarrow\D$ is \textit{surjective} if $\mathrm{Im}(F) = \D$ \cite[Definition 3.2]{bruguieres2011exact}.

Given a tensor functor $F: \C \to \D$ between tensor categories $\C$ and $\D$, we denote by $\mathrm{Im}_s(F)$ the full subcategory of $\D$ consisting of all subobjects of objects of the form $F(X)$ for some $X \in \C$.
A tensor functor $F : \C \to \D$ is said to be \textit{dominant} if $\mathrm{Im}_s(F) = \D$. We note:

\begin{lemma}
  \label{lem:dominant-and-surjective}
  For a tensor functor $F: \C \to \D$ between finite tensor categories $\C$ and $\D$, we have $\mathrm{Im}_s(F) = \mathrm{Im}(F)$.
\end{lemma}
\begin{proof}
This follows from \cite[Lemma~2.3]{etingof2017exact} applied to the tensor functor $\C \to \mathrm{Im}(F)$ induced by $F$.
\end{proof}

\begin{remark}\label{rem:tensor-reflect-exact}
  An exact functor $F$ between abelian categories is faithful if and only if $F(X) \ne 0$ whenever $X \ne 0$. By this fact and \cite[Exercise 4.3.11]{etingof2016tensor}, we see that the tensor product in a finite tensor category is faithful. As $\C$ is rigid, the tensor product is also exact.
  Additionally, an exact faithful functor $F: \mathcal{A} \to \mathcal{B}$ between abelian categories reflects exact sequences. Thus, tensoring reflects exact sequences.
\end{remark}


\subsubsection*{Braided tensor categories}

A tensor category with a braiding is called a braided tensor category.
We say that a braided tensor category $\C$ is \textit{non-degenerate} if $\C' \simeq \Vect$.
A \textit{modular tensor category} (MTC) is a non-degenerate ribbon finite tensor category. 

Let $\C$ be a finite tensor category. Then its Drinfeld center $\Z(\C)$ is a braided finite tensor category \cite[\S8.5]{etingof2016tensor} that is non-degenerate.
We have a forgetful functor from $U_{\C}:\Z(\C)\rightarrow\C$ given by $U_{\C}((X,\sigma)) = X$. It admits a right adjoint, which we will denote as $R_{\C}$.


\subsection{Algebras in monoidal categories}
Let $\C$ be a monoidal category. An \textit{algebra} is a triple $(A,\mu:A\otimes A \rightarrow A,\eta:\unit\rightarrow A)$ where $\mu$ is associative and unital. A left $A$-module is a pair $(M,a^l_M:A\otimes M\rightarrow M)$ consisting of an object $M\in\C$ and a map $a^l_M$ satisfying $a_M^l(\id_A\otimes a_M^l)=a^l_M(\mu\otimes \id)$ and $a^l_M(\eta\otimes\id)=\id_M$. Similarly, a right $A$-module $(M,a^r_M)$ is defined. Given two algebras $A$ and $B$, an $(A,B)$-bimodule is a triple $(M,a^l_M,b^r_M)$ such that $(M,a^l_M)$ is a left $A$-module and $(M,b^r_M)$ is a right $B$-module. For these morphisms, we use the following graphical expressions:
\[
\mu = \begin{tikzpicture}[x = 12pt, y = 12pt, baseline = 0]
\draw (0,1) to [out = -90, in = -90, looseness = 2] coordinate[midway] (M1) ++(2,0);
\node at (M1) {$\bullet$};
\draw (M1) -- (1,-1);
\end{tikzpicture} \qquad
\eta = \begin{tikzpicture}[x = 12pt, y = 12pt, baseline = 0]
\node at (0,1) {$\bullet$};
\draw (0,1) -- (0,-1);
\end{tikzpicture} \qquad
a_M^{l} = \begin{tikzpicture}[x = 12pt, y = 12pt, baseline = 0]
\draw (0,1) -- (0,-1);
\draw (-1,1) to [out = -90, in = 180] (0,0) node {$\bullet$};
\end{tikzpicture} \qquad
b_M^{r} = \begin{tikzpicture}[x = 12pt, y = 12pt, baseline = 0]
\draw (0,1) -- (0,-1);
\draw (1,1) to [out = -90, in = 0] (0,0) node {$\bullet$};
\end{tikzpicture} .
\]

For algebras $A$ and $B$ in a monoidal category $\C$, we denote by ${}_A\C$, $\C_B$ and ${}_A\C_B$ the categories of left $A$-modules, right $B$-modules and $A,B$-bimodules in $\C$, respectively.

\subsubsection*{Central algebras}
An algebra $A\in\C$ is called (lax) \textit{central} if it comes equipped with a (lax) half-braiding $\sigma$ such that $(A,\sigma)$ is an algebra in the (lax) Drinfeld center of $\C$.
By definition, for a (lax) central algebra $(A,\sigma)$, the multiplication $\mu$ and the unit $\eta$ of the algebra $A$ are compatible with $\sigma$ in the sense that the equations
\begin{equation*}
  \sigma_X \circ (\mu \otimes \id_X) 
  = (\id_X\otimes \mu) \circ (\sigma_X \otimes \id_A) \circ (\id_A \otimes \sigma_X),
  \quad \sigma_{X} \circ (\eta \otimes \id_X)
  = \id_X \otimes \eta
\end{equation*}
hold for all $X \in \C$.

\subsubsection*{Commutative algebras}
Let $\B$ be a (lax) braided monoidal category with lax braiding $c$. An algebra $(A,\mu,\eta)$ in $\B$ is called (lax) \textit{commutative} if $\mu = \mu \circ c_{A,A}$. A (lax) \textit{central commutative algebra} is just a commutative algebra in the (lax) Drinfeld center of $\C$.

\subsubsection*{Monoidal category of bimodules}

Let $\C$ be a monoidal category such that coequalizers exist in $\C$ and the tensor product of $\C$ preserves them. Then the tensor product of a right $A$-module $(M,a^r_M)$ and a left $A$-module $(N,a^l_N)$ over $A$ is defined by
\begin{equation*}
M\otimes_A N = \mathrm{coequalizer} \big(
\begin{tikzcd}[column sep = 70pt]
  M\otimes A\otimes N
  \arrow[r, yshift = .3em, "{a_M^r\otimes N}"]
  \arrow[r, yshift = -.3em, "M\otimes a^l_N"']
  & M\otimes N
\end{tikzcd} \big)
\end{equation*}
\cite[Definition~7.8.21]{etingof2016tensor}.
The category ${}_A\C_A$ of $A$-bimodules in $\C$ is a monoidal category with respect to the tensor product over $A$ and with unit $A$.

Let $A$ be an algebra equipped with a lax half-braiding $\sigma$. We define $\C^{\sigma}_A$ to be the full subcategory of ${}_A \C_A$ consisting of all $A$-bimodules $X$ such that
\begin{equation*}
  a_X^{\ell} = a_X^{r} \circ \sigma_X:A\otimes X\rightarrow X.
\end{equation*}
If $A$ is lax central commutative, then the category $\C_A^{\sigma}$ is in fact a monoidal full subcategory of ${}_A\C_A$ \cite{schauenburg2001monoidal}.


\subsubsection*{Braided case}
Suppose that $\C$ is braided and $A$ is a commutative algebra in $\C$.
A right $A$-module is called \textit{local} ($=$ dyslectic \cite{pareigis1995braiding}) if the equation $a^r_M\circ c_{A,M}\circ c_{M,A} = a^r_M$ holds. The full subcategory of $\C_A$ consisting of local $A$-modules is denoted as $\C_A^{\loc}$. 
It is a braided monoidal category with braiding $c^A$ induced by the braiding of $\C$.

\begin{remark}\label{rem:local-int}
There are two functors $F_{\pm}: \C_A\rightarrow {}_A\C_A$ given by
\begin{equation*}\label{eq:F-+}
  F_+ ((M,a^r_M)) = (M,a^r_M , a^r_M \circ c_{A,M}) , \quad F_- ((M,a^r_M)) = (M, a^r_M, a^r_M\circ c^{-1}_{M,A}).
\end{equation*}
The functors $F_{\pm}$ are fully faithful monoidal functors. It is easy to observe that $F_+(\C_A) = \C_A^{c_{A,-}}$ and $F_-(\C_A) = \C_A^{c_{-,A}^{-1}}$. Additionally, $\C_A^{\loc} = \C_A^{c_{A,-}} \cap \C_A^{c_{-,A}^{-1}}$.
\end{remark}
Note that if $\C$ is a (finite) abelian category, then so are the categories ${}_A\C_A$, $\C_A^\sigma$ and $\C_A^\loc$ defined above \cite[Exercise~7.8.16]{etingof2016tensor}.


\subsection{Algebras in rigid monoidal categories}
Let $\C$ be a rigid monoidal category and $A$ an algebra in $\C$. For $(M,a^l_M) \in{}_A\C$ and $(N,a^r_N)\in\C_A$, we have that $(M,a_{M^*}^r)\in\C_A$ and $(N,a^l_{{}^*N}) \in {}_A\C$ where 
\begin{equation}\label{eq:dual-modules-actions}
  a_{M^*}^{r} = (\id_{M^*} \otimes \ev_A) \circ ((a_M^{l})^* \otimes \id_A), \quad 
  a_{{}^*N}^l = (\ev_{{}^*A} \otimes \id_{{}^*N} ) (\id_A \otimes {}^*(a^r_N)).
\end{equation}
Let $B$ also be an algebra in $\C$, and let $(M, b^r_M)$ be a right $B$-module. In general, the left dual object $M^*$ is not a left $B$-module. Instead, $M^*$ is a left $B^{**}$-module by the action
\begin{equation}\label{eq:dual-modules-actions-2}
  b^l_{M^*} := (\ev_{B^*} \otimes \id_{M^*}) \circ (\id_{B^{**}} \otimes (b_M^r)^*) : B^{**} \otimes M^* \to M^*.
\end{equation}
Similarly, the right dual of a left $B$-module is a left ${}^{**}B$-module. Furthermore, as has been pointed out in \cite[Lemma 2.4.13]{douglas2018dualizable}, these constructions give antiequivalences
\begin{equation}\label{eq:dual-bimodules}
  (-)^* : {}_A\C_{B} \to {}_{B^{**}}\C_A
  \quad \text{and} \quad
  {}^*(-) : {}_B\C_{A} \to {}_A\C_{{}^{**}\!B}.
\end{equation}

Now suppose that $\mathcal{C}$ is a tensor category. Given two algebras $A_1$ and $A_2$, their direct sum $A_1 \oplus A_2$ is also an algebra in $\mathcal{C}$. An algebra is called \textit{indecomposable} if it is not isomorphic to a direct sum of two algebras. An algebra $A$ is \textit{haploid} if $\Hom_{\mathcal{C}}(\unit, A) \cong \kk$. Clearly, a haploid algebra is indecomposable. An algebra $A$ is \textit{simple} if it does not admit any two-sided ideals other than the trivial ideal and $A$.

\subsubsection*{Frobenius algebras}
A \textit{Frobenius algebra} is a tuple $(A,\mu,\eta,\Delta,\varepsilon)$ where $(A,\mu,\eta)$ is an algebra, $(A,\Delta,\varepsilon)$ is a coalgebra and $(\mu\otimes\id_A)(\id_A\otimes \Delta) = \Delta \mu = (\id_A\otimes \mu)(\Delta\otimes\id_A)$ holds.

Let $\C$ be a pivotal tensor category with pivotal structure $\mathfrak{p}$, and let $(A, \mu, \eta, \Delta, \varepsilon)$ be a Frobenius algebra in $\C$. There is an isomorphism $\phi := (\varepsilon \mu \otimes \id_{A^*}) (\id_A \otimes \coev_A) : A \to A^*$ of right $A$-modules in $\C$ with the inverse given by $\phi^{-1} = (\ev_A \otimes \id_A)(\id_{A^*} \otimes \Delta \eta)$. The \textit{Nakayama automorphism} \cite{fuchs2008frobenius}\footnote{The morphism $\phi$ considered here corresponds to the morphism $\Phi_{\kappa, \mathrm{r}}$ in \cite{fuchs2008frobenius}, and the morphism $\nu_A$ is, strictly speaking, the inverse of the Nakayama automorphism $\mho$ of \cite{fuchs2008frobenius}.} of $A$ is given by
\begin{equation*}
    \nu_A := \mathfrak{p}_A^{-1} \circ (\phi^{-1})^* \circ \phi : A \to A.
\end{equation*}

\begin{lemma}
    With the above notation, we have
    \begin{equation*}
        \mathrm{tr}(\nu_A) \id_{\unit} = \varepsilon \mu \Delta \eta.
    \end{equation*}
\end{lemma}
\begin{proof}
  By the definition of the pivotal trace and $\nu_A$, we have
  \begin{equation*}
    \mathrm{tr}(\nu_A) \id_{\unit}
    = \begin{tikzpicture}[x = 12pt, y = 12pt, baseline = 0]
      \node (B1) [draw] at (0, 0) {\makebox[2em][c]{$\mathstrut\smash{\phi}$}};
      \node (B2) [draw] at (3, 0) {\makebox[2em][c]{$\mathstrut\smash{\phi^{-1}}$}};
      \draw (B1.north) to [out = 90, in = 90, looseness = 2] (B2.north);
      \draw (B1.south) to [out = -90, in = -90, looseness = 2] (B2.south);
      \node at (-.7, 1.5) {$A$};
      \node at (3.7, 1.5) {$A^*$};
      \node at (-.7, -1.5) {$A^*$};
      \node at (3.5, -1.5) {$A$};
    \end{tikzpicture}
    = \begin{tikzpicture}[x = 12pt, y = 12pt, baseline = 0]
      \draw (0,0) to [out = -90, in = -90, looseness = 2]
      coordinate[midway] (M1) ++(1.5,0)
      to [out = 90, in = 90, looseness = 2] ++(1.5,0) coordinate (L1)
      to [out = -90, in = -90, looseness = 2] ++(4.5,0)
      to [out = 90, in = 90, looseness = 2]
      coordinate[midway] (M2) ++(-1.5,0)
      to [out = -90, in = -90, looseness = 2] ++(-1.5,0) coordinate (L2)
      to [out = 90, in = 90, looseness = 2] cycle;
      \node at (M1) {$\bullet$};
      \draw (M1) -- ++(0,-1) node {$\bullet$};
      \node at (M2) {$\bullet$};
      \draw (M2) -- ++(0,1) node {$\bullet$};
      \node at ($(L1)+(-.5,-2)$) {$A^*$};
      \node at ($(L2)+( .5, 2)$) {$A^*$};
    \end{tikzpicture}
    = \varepsilon \mu \Delta \eta,
  \end{equation*}
  where $\Delta$ and $\varepsilon$ are expressed by the upside down diagrams of $\mu$ and $\eta$.
\end{proof}

A Frobenius algebra $(A, \mu, \eta, \Delta, \varepsilon)$ is called \textit{special} if $\varepsilon\circ \eta=\beta_0\id_{\unit}$ and $\mu\circ\Delta=\beta_2\id_A$ for some $\beta_0,\beta_2\in\kk^{\times}$.
By renormalizing $\Delta$ and $\varepsilon$, we may assume that $\beta_0 = 1$. Then, by the above lemma, we have $\beta_2 = \mathrm{tr}(\nu_A)$. We have proved that the trace of $\nu_A$ is non-zero if $A$ is special Frobenius. Under the assumption that $A$ is haploid, the converse is also true:

\begin{lemma}
  \label{lem:special-Frobenius-trace-Nakayama}
  A haploid Frobenius algebra in $\C$ is special if and only if $\mathrm{tr}(\nu_A) \ne 0$.
\end{lemma}
\begin{proof}
    By the above discussion, $A$ cannot be special Frobenius when $\mathrm{tr}(\nu_A) = 0$. We now assume that $\mathrm{tr}(\nu_A) \ne 0$ and show that $A$ is special Frobenius.
    By the assumption that $A$ is haploid, we have $\Hom_{\C_A}(A, A) \cong \Hom_{\C}(\unit, A) \cong \kk$.
    This implies $\mu \Delta = \beta_2 \id_A$ for some scalar $\beta_2 \in \kk$, and therefore we have $\beta_2 \varepsilon \eta = \varepsilon \mu \Delta \eta = \mathrm{tr}(\nu_A) \id_{\unit}\ne 0$.
    Therefore $\beta_2 \ne 0$ and $\varepsilon \eta = \beta_0 \id_{\unit}$ for some $\beta_0 \in \kk^{\times}$. The proof is done.
\end{proof}
We say that a Frobenius algebra $A\in\C$ is \textit{symmetric} if $\nu_A = \id_A$.

\begin{lemma}\cite[Proposition~2.25]{frohlich2006correspondences}
  \label{lem:symmetric-frobenius-theta}
  Let $\C$ be a braided category with a twist $\theta$. Then, a commutative Frobenius algebra $A$ is symmetric if and only if $\theta_A=\id_A$. \qed
\end{lemma}


\section{Rigidity and ribbon structure of \texorpdfstring{$\C_A^{\loc}$}{C-A-loc}}\label{sec:closed}


A monoidal category $\C$ is called \textit{closed} \cite{eilenberg1966closed}, if the functors $X\otimes -$ and $-\otimes X$ admit a right adjoint for all $X\in\C$.
If $\C$ is a closed monoidal category with equalizers and coequalizers, then the monoidal category ${}_A\C_A$ is closed for every algebra $A$ in $\C$; see Section~\ref{subsec:closed-bimodules}.
In Section~\ref{subsec:closed-modules}, we consider the setting of $(A,\sigma)$ being a lax central commutative algebra and show that $\C_A^{\sigma}$ is a monoidal subcategory of ${}_A\C_A$ with the same internal Hom functor as ${}_A\C_A$. In Section~\ref{subsec:closed-rigid-A-C-A}, we note that this observation implies that $\C_A^{\sigma}$ is rigid provided ${}_A\C_A$ is rigid.
When $\C$ is braided and $A$ is a commutative algebra in $\C$, this implies that $\C_A$ and $\C_A^{\loc}$ are rigid provided ${}_A\C_A$ is rigid. When $\C$ is ribbon, in Section~\ref{subsec:ribbon-case}, we give sufficient conditions for $\C_A^{\loc}$ to be a ribbon monoidal category. 


\subsection{Closed monoidal categories}
Let $\C$ be a closed monoidal category, and let $\uHom(X,-):\C\rightarrow\C$ denote the right adjoint of $- \otimes X$. One can extend this to a bifunctor $\uHom : \C^{\op} \times \C \to \C$, called the \emph{internal Hom functor} of $\C$. Thus there is a natural isomorphism:
\begin{equation}
  \label{eq:iHom-adj}
  \Hom_{\C}(W \otimes X, Y)
  \cong \Hom_{\C}(W, \uHom(X, Y))
  \quad (W, X, Y \in \C).
\end{equation}
The unit and the counit of this adjunction are denoted by
\begin{equation*}
  \underline{\coev}_{X,Y}: Y \to \uHom(X, Y \otimes X),
  \quad \underline{\ev}_{X,Y}: \uHom(X, Y) \otimes X \to Y.
\end{equation*}
\begin{remark}
  When working with closed monoidal categories we will use (without proof) various identities of the form $f=g:X\rightarrow \uHom(M,N)$. Such identities can be readily proved by showing that the adjuncts of $f$ and $g$, namely $\underline{\ev}_{M,N}\circ (f\otimes \id_M)$ and $\underline{\ev}_{M,N}\circ (g\otimes \id_M)$ are equal.
\end{remark}
For $W, X, Y \in \C$, we define the morphism
\begin{equation*}
  \underline{\mathord{\otimes} W}_{X,Y}:
  \uHom(X, Y) \to \uHom(X \otimes W, Y \otimes W)
\end{equation*}
to be the adjunct of the morphism
\begin{equation*}
  \uHom(X, Y) \otimes X \otimes W
  \xrightarrow{\quad \underline{\ev}_{X,Y} \otimes \id_W \quad}
  Y \otimes W.
\end{equation*}


\subsection{Explicit description of the internal Hom in \texorpdfstring{{}${}_A\C_A$}{ACA}}\label{subsec:closed-bimodules}

Let $\C$ be a closed monoidal category admitting equalizers and coequalizers, and let $A$ be an algebra in $\C$.
We note that the tensor product of $\C$ preserves coequalizers by the general fact that a functor admitting a right adjoint preserves colimits.
Consider $M \in {}_A \C$. Then, the left action gives rise to a morphism
\begin{equation*}
  \rho_M^{\ell} := \uHom(\id_M, a_M^{\ell}) \circ \underline{\coev}_{M,A} : A \to \uHom(M, M)
\end{equation*}
of algebras in $\C$.
When $M, N \in {}_A \C_A$, $\uHom(M, N)$ has a structure of $A$-bimodule given by
\begin{equation}
  \label{eq:iHom-action-def}
  \begin{aligned}
    a_{\uHom(M, N)}^{\ell} & = \underline{\circ}_{M,N,N} \circ (\rho_N^{\ell} \otimes \id_{\uHom(M, N)}), \\
    a_{\uHom(M, N)}^{r} & = \underline{\circ}_{M,M,N} \circ (\id_{\uHom(M, N)} \otimes \rho_M^{\ell}),
  \end{aligned}
\end{equation}
where $\underline{\circ}_{X,Y,Z} : \uHom(Y,Z) \otimes \uHom(X,Y) \to \uHom(X,Z)$ is the composition for the internal Hom functor. We note that the following equations hold:
\begin{equation}
  \label{eq:iHom-action-def-adjunt}
  \begin{aligned}
    \underline{\ev}_{M,N} \circ (a_{\uHom(M, N)}^{\ell} \otimes \id_M)
    & = a_N^{\ell} \circ (\id_A \otimes \underline{\ev}_{M,N}), \\
    \underline{\ev}_{M,N} \circ (a_{\uHom(M,N)}^{r} \otimes \id_M)
    & = \underline{\ev}_{M,N} \circ (\id_{\uHom(M,N)} \otimes a_M^{r}).
  \end{aligned}
\end{equation}

\begin{definition}
  For $X, Y \in {}_A \C_A$, we define
  \begin{equation}\label{eq:iHom-A-C-A-def}
    \uHom_A(X, Y) = \mathrm{equalizer} \Big(
    \begin{tikzcd}[column sep = 100pt]
      \uHom(X, Y)
      \arrow[r, yshift = .3em, "{\uHom(a_X^{r}, \id_Y)}"]
      \arrow[r, yshift = -.3em, "{\uHom(\id_{X \otimes A}, a_Y^{r}) \circ \underline{\mathord{\otimes}A}_{X,Y}}"']
      & \uHom(X \otimes A, Y)
    \end{tikzcd} \Big)
  \end{equation}
  We make $\uHom_A(X,Y)$ an $A$-bimodule as a subbimodule of $\uHom(X,Y)$.
  Namely, the left and the right action of $A$ are morphisms determined by
  \begin{equation}
    \label{eq:iHom-A-X-Y-action-def}
    \begin{aligned}
      i_{X,Y} \circ a_{\uHom_A(X,Y)}^{\ell}
      & = a^{\ell}_{\uHom(X,Y)} \circ (\id_A \otimes i_{X,Y}), \\
      i_{X,Y} \circ a_{\uHom_A(X,Y)}^{r}
      & = a^{r}_{\uHom(X,Y)} \circ (i_{X,Y} \otimes \id_A),
    \end{aligned}
  \end{equation}
  where $i_{X,Y}: \uHom_A(X, Y) \to \uHom(X, Y)$ is the canonical monomorphism.
\end{definition}

As $i_{X,Y}$ is an equalizer, we have the following equation:
\begin{equation}
  \label{eq:iHom-A-X-Y-equalizer-1}
  \uHom(a_X^r, \id_Y) \circ i_{X,Y}
  = \uHom(\id_{X \otimes A}, a_Y^r) \circ \underline{\mathord{\otimes}A}_{X,Y} \circ i_{X,Y}.
\end{equation}
By taking adjuncts of both sides in \eqref{eq:iHom-A-X-Y-equalizer-1}, we have
\begin{equation}
  \label{eq:iHom-A-X-Y-equalizer-2}
  \underline{\ev}_{X,Y} \circ (i_{X,Y} \otimes a_{X}^r)
  = a_Y^r \circ (\underline{\ev}_{X,Y} \otimes \id_A) \circ (i_{X,Y} \otimes \id_X \otimes \id_A).
\end{equation}

\begin{lemma}\label{lem:bimodClosed}
  The adjunction isomorphism \eqref{eq:iHom-adj} induces a natural isomorphism
  \begin{align*}
    \Hom_{{}_A\C_A}(W \otimes_A X,Y)
    \cong \Hom_{{}_A\C_A}(W, \uHom_A(X,Y)). \qed
  \end{align*}
\end{lemma}
This proves that $-\otimes_A X$ admits a right adjoint, namely $\uHom_A(X,-)$. A similar argument shows that $X\otimes_A -$ admits a right adjoint. Thus, we obtain the following corollary.

\begin{corollary}\label{cor:bimodClosed}
  The category ${}_A\C_A$ is closed monoidal. \qed
\end{corollary}


\subsection{\texorpdfstring{$\C_A^{\sigma}$}{CAsigma} is closed}\label{subsec:closed-modules}
Let $\C$ be a closed monoidal category admitting equalizers and coequalizers, and let $(A,\sigma)$ be a lax central commutative algebra in $\C$.
Then we have the following result.

\begin{theorem}\label{thm:RepA-int-hom}
  If $X, Y \in \C^{\sigma}_A$, then $\uHom_A(X,Y) \in \C^{\sigma}_A$.
\end{theorem}

For our purpose, it suffices to consider the case where $\sigma$ is invertible. Nevertheless, since the invertibility of $\sigma$ does not simplify the proof, we state the theorem in a form that includes the lax setting with an eye toward future applications to some theories where lax central commutative algebras naturally arise, like \cite{bruguieres2011hopf}.

\begin{proof}
  Let $X, Y \in \C_A^{\sigma}$. It suffices to show that the equation
  \begin{equation}
    \label{eq:iHom-A-X-Y-commutativity}
    i_{X,Y} \circ a_{\uHom_A(X,Y)}^{\ell} = i_{X,Y} \circ a_{\uHom_A(X,Y)}^{r} \circ \sigma_{\uHom_A(X,Y)}
  \end{equation}
  holds. We note that both sides are morphisms
  \begin{equation*}
    A \otimes \uHom_A(X,Y) \to \uHom(X,Y).
  \end{equation*}
  We compute the adjuncts of both sides of \eqref{eq:iHom-A-X-Y-commutativity}. The adjunct of the left hand side is:
  \begin{align*}
    & \underline{\ev}_{X,Y} \circ ((\text{the left hand side of \eqref{eq:iHom-A-X-Y-commutativity}}) \otimes \id_X) \\
    {}^{\eqref{eq:iHom-A-X-Y-action-def}}
    & = \underline{\ev}_{X,Y} \circ (a_{\uHom(X,Y)}^{\ell}(\id_A \otimes i_{X,Y}) \otimes \id_X) \\
    {}^{\eqref{eq:iHom-action-def-adjunt}}
    & = a_Y^{\ell} \circ (\id_A \otimes \underline{\ev}_{X,Y}) \circ (\id_A \otimes i_{X,Y} \otimes \id_X) \\
    & = a_Y^r \circ \sigma_Y \circ (\id_A \otimes \underline{\ev}_{X,Y}) \circ (\id_A \otimes i_{X,Y} \otimes \id_X)
      \quad \text{(since $Y \in \C_A^{\sigma}$)} \\
    {}^{\text{(nat.)}}
    & = a_Y^r \circ (\underline{\ev}_{X,Y} \otimes \id_A) \circ (i_{X,Y} \otimes \id_X \otimes \id_A) \circ \sigma_{\uHom_A(X,Y) \otimes X} \\
    {}^{\eqref{eq:iHom-A-X-Y-equalizer-2}}
    & = \underline{\ev}_{X,Y} \circ (i_{X,Y} \otimes a_{X}^r)
      \circ \sigma_{\uHom_A(X,Y) \otimes X},
  \end{align*}
  where the naturality of the half-braiding $\sigma$ is used at (nat.).
  The adjunct of the right hand side is:
  \begin{align*}
    & \underline{\ev}_{X,Y} \circ ((\text{the right hand side of \eqref{eq:iHom-A-X-Y-commutativity}}) \otimes \id_X) \\
    {}^{\eqref{eq:iHom-A-X-Y-action-def}}
    & = \underline{\ev}_{X,Y} \circ (a_{\uHom(X,Y)}^r(i_{X,Y} \otimes \id_A) \otimes \id_X)
      \circ (\sigma_{\uHom_A(X,Y)} \otimes \id_X) \\
    {}^{\eqref{eq:iHom-action-def-adjunt}}
    & = \underline{\ev}_{X,Y} \circ (i_{X,Y} \otimes a_{X}^{\ell})
      \circ (\sigma_{\uHom_A(X,Y)} \otimes \id_X) \\
    & = \underline{\ev}_{X,Y} \circ (i_{X,Y} \otimes a_{X}^{r} \sigma_X)
      \circ (\sigma_{\uHom_A(X,Y)} \otimes \id_X)
      \quad \text{(since $X \in \C_A^{\sigma}$)} \\
    {}^{\eqref{eq:lax-half-braiding}}
    & = \underline{\ev}_{X,Y} \circ (i_{X,Y} \otimes a_{X}^r)
      \circ \sigma_{\uHom_A(X,Y) \otimes X}
  \end{align*}
  Thus \eqref{eq:iHom-A-X-Y-commutativity} is verified.
\end{proof}

It is well-known that $\C^{\sigma}_A$ is a monoidal subcategory of ${}_A\C_A$ with $\unit_{\C_A^{\sigma}}=A$ and tensor product over $A$. Lemma~\ref{lem:bimodClosed} and Theorem~\ref{thm:RepA-int-hom} together imply:

\begin{corollary}\label{cor:C-A-sigma-closed}
  Let $\C$ be a closed monoidal category admitting equalizers and coequalizers, and let $(A, \sigma)$ be a lax central commutative algebra in $\C$. 
  Then the category $\C^{\sigma}_A$ is closed monoidal with the same internal Hom functor as the monoidal category ${}_A\C_A$. \qed
\end{corollary}
Using Remark~\ref{rem:local-int}, the following is immediate.

\begin{corollary}\label{C-A-loc-closed}
  Let $\C$ be a braided closed monoidal category admitting equalizers and coequalizers, and let $A$ be a commutative algebra in $\C$. Then the categories $\C_A$ and $\C_A^{\loc}$ are closed monoidal with the same internal Hom functor as the monoidal category ${}_A\C_A$.\qed
\end{corollary}


\subsection{Duals in \texorpdfstring{$\C_A^{\sigma}$}{CAsigma}}\label{subsec:closed-rigid-A-C-A}

Let $\D$ be a closed monoidal category with internal Hom functor $\uHom_{\D}$. Then there is a natural transformation
\begin{equation*}
  \varsigma_{M,N} : N \otimes \uHom_{\D}(M, \unit) \to \uHom_{\D}(M, N)
  \quad (M, N \in \D)
\end{equation*}
induced by the morphism $\id_N \otimes \underline{\ev}_{M, \unit} : N \otimes \uHom_{\mathcal{D}}(M, \unit) \otimes M \to N$.
It is known that $M \in \D$ is left rigid if and only if $\varsigma_{M,N}$ is an isomorphism for all objects $N \in \D$ \cite[Remark~3.2]{bruguieres2011hopf}. Furthermore, if this is the case, a left dual object of $M$ is given by $M^{\dagger} := \uHom_{\D}(M, \unit)$ with the evaluation and the coevaluation morphisms given respectively by
\begin{equation*}
  \underline{\ev}_{M, \unit} : M^{\dagger} \otimes M \to \unit
  \quad \text{and} \quad
  \varsigma_{M,M}^{-1} \circ \underline{\coev}_{M, \unit} : \unit \to M \otimes M^{\dagger}.
\end{equation*}

\begin{theorem}\label{thm:ACA-rigid-consequences}
  Let $\C$ be a closed monoidal category admitting equalizers and coequalizers.
  \begin{enumerate}
  \item[\textup{(a)}] If $(A,\sigma)$ is a lax central commutative algebra in $\C$ such that ${}_A\C_A$ is rigid, then $\C_A^{\sigma}$ is rigid.
  \item[\textup{(b)}] If $\C$ is braided and $A$ is a commutative algebra in $\C$ such that ${}_A\C_A$ is rigid, then both $\C_A$ and $\C_A^{\loc}$ are rigid. 
  \end{enumerate}
  In all cases, the left duality functor of the category is given by
  \[M \mapsto M^{\dagger} := \uHom_A(M, A).  \]
\end{theorem}
\begin{proof}
By the above argument and Corollary~\ref{cor:C-A-sigma-closed}, as well as their left-right reversed versions, $M \in \C_A^{\sigma}$ has a left (right) dual object in $\C_A^{\sigma}$ if and only if it has a left (right) dual object in ${}_A\C_A$. This proves (a). Using Corollary~\ref{C-A-loc-closed}, part (b) follows.
\end{proof}


\subsection{Rigid case}\label{subsec:rigid-case}
In this section $\C$ is assumed to be a rigid monoidal category with coequalizers and that the tensor product preserves coequalizers. Then, $\C$ is closed monoidal and the internal Hom functor is given by $\uHom(X, Y) = Y \otimes X^*$ for $X, Y \in \C$.
Hence, by \eqref{eq:iHom-A-C-A-def}, the internal Hom functor of ${}_A\C_A$ is given by
\begin{equation}\label{eq:iHom-A-in-ACA}
  \uHom_A(M, N) = \mathrm{equalizer} \big(
  \begin{tikzcd}[column sep = 130pt]
    N\otimes M^*
    \arrow[r, yshift = .3em, "{N\otimes (a^r_M)^*}"]
    \arrow[r, yshift = -.3em, "(a^r_N\otimes A^*\otimes M^*) \circ (N\otimes \coev_A \otimes M^*)"']
    & N\otimes A^*\otimes M^*
  \end{tikzcd} \big).
\end{equation}
The canonical monomorphism for $\uHom_A(M, N)$ as an equalizer \eqref{eq:iHom-A-in-ACA} is given by
\begin{equation*}
  i_{M,N} := (\pi_{M,{}^*N})^* : \uHom_A(M, N)
  \to (M \otimes {}^*N)^* = N \otimes M^*,
\end{equation*}
where $\pi_{X,Y} : X \otimes Y \to X \otimes_A Y$ for $X \in \C_A$ and $Y \in {}_A\C$ is the canonical epimorphism.
From (\ref{eq:iHom-A-in-ACA}) and the definition of $\otimes_A$, we also have
\begin{equation}
  \label{eq:internal-Hom-rigid-case}
  \uHom_A(M,N) = (M\otimes_A {}^*N)^*.
\end{equation}
If $M$ and $N$ are $A$-bimodules, then $N \otimes M^*$ is also an $A$-bimodule by
\begin{equation*}
  a_{N \otimes M^*}^{l} = a_N^{l} \otimes \id_{M^*}
  \quad \text{and}
  \quad a_{N \otimes M^*}^{r} = \id_N \otimes a_{M^*}^{r},
\end{equation*} 
and we can make $\uHom_A(M, N)$ an $A$-bimodule as a subbimodule of $N \otimes M^*$.


\subsection{Technical results on duality}\label{subsec:central-rigid-case}
We assume that $\mathbf{A} = (A, \sigma)$ is a central commutative algebra. As we have remarked, the duality functor on $\C_A^{\sigma}$ is given by $M \mapsto \uHom_A(M, A)$ provided that $\C_A^{\sigma}$ is rigid. However, there is another duality functor $D$ on $\C_A^{\sigma}$ constructed based on the duality functor of $\mathcal{C}$ as follows: If $M$ is an object of $\C_A^{\sigma}$, then $D(M) := M^*$ also belongs to $\C_A^{\sigma}$ by the right action $a^r_{D(M)} := a^r_{M^*}$ given by \eqref{eq:dual-modules-actions} and the left action $a^l_{D(M)} := a^r_{D(M)} \circ \sigma_{M^*}$. The assignment $M \mapsto D(M)$ gives rise to an anti-autoequivalence $D$ on the category $\C_A^{\sigma}$.

The functor $D$ plays an important role when we discuss ribbon structures on $\C_A^{\loc}$. Here we provide technical lemmas concerning the functor $D$. We recall that the algebra $A^{**}$ also acts on $M^{*}$ from the left by the action $b_{M^*}^l : A^{**} \otimes M^* \to M^*$ given by \eqref{eq:dual-modules-actions} with $b_M^r = a_M^r$. The actions $b^l_{M^*}$ and $a^l_{D(M)}$ have the following relation:

\begin{lemma}
\label{lem:GV-duality-on-CA-1}
With the above notation, we have $a^l_{D(M)} = b_{M^*}^l \circ (u \otimes \id_{M^*})$, where $u = \mathfrak{u}_{\mathbf{A}}$ is the component of the Drinfeld isomorphism in $\Z(\C)$ for $\mathbf{A}$.
\end{lemma}
\begin{proof}
  By the definition of actions, we have
\begin{align*}
  \ev_M \circ (a_{D(M)}^{l} \otimes \id_{M^*})
  & = \ev_M \circ (a_{D(M)}^r \otimes \id_M) \circ (\sigma_{M^*} \otimes \id_M) \\
  & = \ev_M \circ (\id_{M^*} \otimes a_M^{l}) \circ (\sigma_{M^*} \otimes \id_M) \\
  & = \ev_M \circ (\id_{M^*} \otimes a_M^{r})
    \circ (\id_{M^*} \otimes \sigma_{M})
    \circ (\sigma_{M^*} \otimes \id_M) \\
  {}^{\eqref{eq:lax-half-braiding}}
  & = \ev_{M \otimes A} \circ ((a_M^{r})^* \otimes \id_{M} \otimes \id_{A})
    \circ \sigma_{M^* \otimes M}.
\end{align*}
Now we set $\xi = (\id_{A^*} \otimes \ev_M) \circ ((a_M^r)^* \otimes \id_{M})$ for simplicity of notation. Since $b^l_{M^*}$ is defined by \eqref{eq:dual-modules-actions} with $b_M^r = a_M^r$, we have
\begin{equation*}
  \ev_{A^*} (\id_{A^{**}} \otimes \xi)
  = (\ev_{A^*} \otimes \ev_M) (\id_{A^{**}} \otimes (a_M^r)^* \otimes \id_{M})
  = \ev_M (b_{M^*}^{l} \otimes \id_{M}),
\end{equation*}
and thus we continue the above computation as follows:
\begin{align*}
  ...
  & = \ev_{A} \otimes (\id_{A^*} \otimes \ev_M \otimes \id_{A})
    \circ ((a_M^{r})^* \otimes \id_{M} \otimes \id_{A})
    \circ \sigma_{M^* \otimes M} \\
  & = \ev_{A} \circ (\xi \otimes \id_{A}) \circ \sigma_{M^* \otimes M} \\
  {}^{\text{(nat.)}}
  & = \ev_{A} \circ \sigma_{A^*} \circ (\id_{A} \otimes \xi) \\
  {}^{\eqref{eq:def-Drinfeld-iso}}
  & = \ev_{A^*} \circ (u \otimes \id_{A^*}) \circ (\id_{A} \otimes \xi) \\
  & = \ev_{M} \circ (b^l_{M^*} \otimes \id_M) \circ (u \otimes \id_{M^*} \otimes \id_{M}),
\end{align*}
where `(nat.)' follows from the naturality of the half-braiding $\sigma$.
The proof is done, since the map $f \mapsto \ev_M (f \otimes \id_{M})$ is injective.
\end{proof}

The above lemma yields a handy description of $D^2$:

\begin{lemma}
\label{lem:GV-duality-on-CA-2}
For $M \in \C_A^{\sigma}$, the right action of $A$ on $D^2(M) = M^{**}$ is given by
\begin{equation*}
  a^r_{D^2(M)} = (a^r_M)^{**} \circ (\id_{M^{**}} \otimes u).
\end{equation*}
\end{lemma}
\begin{proof}
  Let $b^l_{M^*}$ be as in the previous lemma. We recall that $a^r_{D^2(M)}$ is defined by \eqref{eq:dual-bimodules} with $M$ replaced by $D(M)$. By the previous lemma, we have
  \begin{align*}
    a^r_{D^2(M)}
    & = (\id_{M^{**}} \otimes \ev_A) \circ ((a^l_{D(M)})^* \otimes \id_{A}) \\
    & = (\id_{M^{**}} \otimes \ev_A)
      \circ (\id_{M^{**}} \otimes u^* \otimes \id_{A})
      \circ ((b^l_{M^*})^* \otimes \id_{A}) \\
    {}^{\eqref{eq:dual-modules-actions}}
    & = (\id_{M^{**}} \otimes \ev_A)
      \circ (\id_{M^{**}} \otimes u^* \otimes \id_{A}) \\*
    & \quad \circ ((a^r_M)^{**} \otimes \id_{M^{**}} \otimes \id_{A})
      \circ (\id_{M^{**}} \otimes (\ev_{A^{*}})^* \otimes \id_A).
  \end{align*}
  Now the proof is completed by an easy graphical calculus.
\end{proof}


\subsection{Ribbon case}\label{subsec:ribbon-case}
In this section, we assume that $\C$ is a braided monoidal category satisfying the assumptions of \S\ref{subsec:rigid-case}. 
We will now give sufficient conditions for $\C_A^{\loc}$ to be ribbon.
To start, we note the following result.
\begin{theorem}\cite[Theorem~1.17]{kirillov2002q}
  Let $\C$ be a braided category with a twist $\theta$. If $A$ is a commutative algebra in $\C$ such that $\theta_A = \id_A$, then $\theta$ gives a twist on $\C_A^{\loc}$. \qed
\end{theorem}

From now on, $\C$ will be a ribbon category with braiding $c$ and twist $\theta$.
We have introduced an anti-autoequivalence $D$ on $\C_A^{\sigma}$, where $\sigma_X = c_{A,X}$ ($X \in \C$), in \S\ref{subsec:central-rigid-case}. A description of $D^2$ has been given by Lemma \ref{lem:GV-duality-on-CA-2}. We use this lemma to prove:

\begin{proposition}\label{prop:D2-on-C-A-loc}
  If $\theta_A = \id_A$, then there is a natural isomorphism
  \begin{equation*}
      D^2(M) \cong M \quad (M \in \C_A^{\sigma}).
  \end{equation*}
\end{proposition}
\begin{proof}
  Set $\mathfrak{p}_X = \mathfrak{u}_X \theta_X$ for $X \in \C$, where $\mathfrak{u}$ is the Drinfeld isomorphism of $\C$. Then $\mathfrak{p} = \{ \mathfrak{p}_X \}_{X \in \C}$ is a pivotal structure on $\C$. Without assuming $\theta_A = \id_A$, for $M \in \C_A^{\sigma}$, we have
  \begin{gather*}
     a^r_{D^2(M)} \circ (\mathfrak{p}_M \otimes \id_A)
     = (a^r_M)^{**} \circ (\mathfrak{p}_M \otimes \mathfrak{u}_{A})
     = (a^r_M)^{**} \circ (\mathfrak{p}_M \otimes \mathfrak{p}_{A} \theta_A^{-1}) \\
     = (a^r_M)^{**} \circ \mathfrak{p}_{M \otimes A} \circ (\id_M \otimes \theta_A^{-1})
     = \mathfrak{p}_M \circ a^r_M \circ (\id_M \otimes \theta_A^{-1})
  \end{gather*}
  by Lemma \ref{lem:GV-duality-on-CA-2}. In other words, $\mathfrak{p}$ gives an isomorphism between $D^2$ and the functor induced by the algebra automorphism $\theta_A^{-1} : A \to A$. Now the claim is obvious.
\end{proof}

It is natural to ask when $\C_A^{\loc}$ is closed under $D$. We show that this is the case precisely if the equation $\theta_A^2 = \id_A$ holds. More precisely, we have:

\begin{proposition}\label{prop:K-in-loc}
$D(A) \in \C_A^{\loc}$ if and only if $\theta_A^2 = \id_A$.
If these equivalent conditions are satisfied, then we have $D(M) \in \C_A^{\loc}$ for all $M \in \C_A^{\loc}$.
\end{proposition}
\begin{proof}
  Let $M \in \C_A^{\loc}$. We first discuss when $D(M)$ belongs to $\C_A^{\loc}$. By the naturality of $\theta$ and the locality of $M$, we have
  \begin{equation}
    \label{eq:prop:K-in-loc-proof-1}
    \theta_M \circ a^r_M
    = a^r_M \circ \theta_{M \otimes A}
    = a^r_M \circ c_{A,M} \circ c_{M,A} \circ (\theta_{M} \otimes \theta_A)
    = a^r_M \circ (\theta_{M} \otimes \theta_A).
  \end{equation}
  Hence, we compute $a^r_{D(M)} \circ c_{M^*,A}^{-1} \circ c_{A,M^{*}}^{-1}$
  \begin{align*}
    {}^{\eqref{eq:def-twist}}
    & = a^r_{D(M)} \circ \theta_{M^* \otimes A}^{-1} \circ (\theta_{M^*} \otimes \theta_A) \\
    & = \theta_{M^*}^{-1} \circ a^r_{D(M)} \circ (\theta_{M^*} \otimes \theta_A) \\
    {}^{\eqref{eq:dual-modules-actions}}
    & = (\theta_{M}^*)^{-1} \circ (\id_{M^*} \otimes \ev_A) \circ ((a^r_M)^* \otimes \id_{A}) \circ ((\theta_{M})^{*} \otimes \theta_A) \\
    {}^{\eqref{eq:prop:K-in-loc-proof-1}}
    & = (\theta_{M}^*)^{-1} \circ (\id_{M^*} \otimes \ev_A)
    \circ (\theta_M^* \otimes \theta_A^* \otimes \id_{A})
    \circ ((a^r_M)^* \otimes \id_{A})
    \circ (\id_{M^*} \otimes \theta_A) \\
    & = (\theta_{M}^*)^{-1} \circ (\id_{M^*} \otimes \ev_A)
    \circ (\theta_M^* \otimes \id_{A^*} \otimes \theta_{A})
    \circ ((a^r_M)^* \otimes \id_{A})
    \circ (\id_{M^*} \otimes \theta_A) \\
    {}^{\eqref{eq:dual-modules-actions}}
    & = a^r_{D(M)} \circ (\id_{M^*} \otimes \theta_A^2),
  \end{align*}
  where the second to last equation follows from the dinaturality of $\ev_X$ in $X \in \C$.
  Thus we have obtained the following conclusion: $D(M)$ belongs to $\C_A^{\loc}$ if and only if
  \begin{equation}
    \label{eq:prop:K-in-loc-proof-2}
     a^r_{D(M)} = a^r_{D(M)} \circ (\id_{M^*} \otimes \theta_A^2).
  \end{equation}
  Now we show that $D(A) \in \C_A^{\loc}$ if and only if $\theta_A^2 = \id_A$.
  The `if' part is trivial by \eqref{eq:prop:K-in-loc-proof-2}. The `only if' part is obtained by applying the map
  \begin{equation*}
     \Hom_{\C}(A^* \otimes A, A^*) \to \Hom_{\C}(A,A), 
     \quad f \mapsto (\id_A \otimes \eta^*f) (\coev_A \otimes \id_A)
  \end{equation*}
  to both sides of \eqref{eq:prop:K-in-loc-proof-2} for $M = A$, where $\eta$ is the unit of $A$.
  The latter part of the statement of this proposition is obvious from \eqref{eq:prop:K-in-loc-proof-2}.
\end{proof}

\begin{remark}
Let $\C$ be a rigid braided category and $A$ a commutative algebra in $\C$. As we have noted, the rigid (left) dual $M^\dagger$ of an object $M$ in the category $\C_A$ (or $\C_A^\loc$) is typically different from $M^*$. Nevertheless, $\C_A$ (and $\C_A^\loc$, under the additional assumption that $\theta^2=\id$) is closed under the functor $D=(-)^*$. Moreover, $D$ is an antiequivalence. In fact, by \cite[Theorem~3.9]{mcrae2025rigidity}, both $\C_A$ and $\C_A^\loc$ are Grothendieck-Verdier categories \cite{boyarchenko2013duality}, under the assumptions mentioned, with duality functor $D$ and dualizing object $K=A^*$. 
Also, compare Proposition~\ref{prop:K-in-loc} with \cite[Theorem~3.11]{mcrae2025rigidity}.
\end{remark}

Now we prove the following main theorem of this section:

\begin{theorem}
\label{thm:local-module-ribbon}
Suppose that $(A,\mu,\eta)$ is a commutative algebra in $\C$ such that $\C_A^{\loc}$ is rigid.
Suppose moreover that $\theta_A = \id_A$ holds, so that we have the object $D(A) \in \C_A^{\loc}$ of Proposition \ref{prop:K-in-loc}.
Then the following assertions are equivalent:
\begin{enumerate}
\item[\textup{(a)}] $\C_A^{\loc}$ is a ribbon monoidal category with the twist $\theta$.
\item[\textup{(b)}] $D(A)$ belongs to the M\"uger center of $\C_A^{\loc}$.
\end{enumerate}
\end{theorem}
\begin{proof}
We use the setup of the rigid case, see \S\ref{subsec:rigid-case}.
By Proposition \ref{prop:D2-on-C-A-loc}, the $A$-modules $D(A)$ and $D^{-1}(A)$ are isomorphic. Thus, below, we actually show that (a) is equivalent to $D^{-1}(A) \in (\C_A^{\loc})'$. We compute $\theta_{X^{\dagger}}$ ($X^{\dagger} = \uHom_A(X, A)$) for $X \in \C_A^{\loc}$ as follows:
\begin{gather*}
  i_{X,A} \circ \theta_{X^{\dagger}}
  = \theta_{A \otimes X^*} \circ i_{X,A}
  = c_{X^*,A} \circ c_{A,X^*} \circ (\theta_A \otimes \theta_X^*) \circ i_{X,A} \\
  = c_{X^*,A} \circ c_{A,X^*} \circ i_{X,A} \circ \uHom_A(\theta_X, A)
  = c_{X^*,A} \circ c_{A,X^*} \circ i_{X,A} \circ (\theta_X)^{\dagger}.
\end{gather*}
Thus, if the equation $\theta_{X^{\dagger}} = (\theta_X)^{\dagger}$ holds, then we have
\begin{equation}
  \label{eq:Mueger-centrality-of-K-1}
  i_{X,A} = c_{X^*,A} \circ c_{A,X^*} \circ i_{X,A}
\end{equation}
by the invertibility of $\theta_{X^{\dagger}}$. If, conversely, the equation \eqref{eq:Mueger-centrality-of-K-1} holds, then we have $\theta_{X^{\dagger}} = (\theta_X)^{\dagger}$ by the above computation and the monicity of $i_{X,A}$. Summarizing, the equation $\theta_{X^{\dagger}} = (\theta_X)^{\dagger}$ is equivalent to \eqref{eq:Mueger-centrality-of-K-1}. Therefore the condition (a) of this theorem is equivalent to that the equation \eqref{eq:Mueger-centrality-of-K-1} holds for all $X \in \C_A^{\loc}$. By the definition of the braiding of $\C_A^{\loc}$, we have
\begin{gather*}
  {}^{*}(\text{the right hand side of \eqref{eq:Mueger-centrality-of-K-1}})
  = \pi_{X, D^{-1}(A)} \circ {}^*(c_{X^*,A} \circ c_{A,X^*}) \\
  = \pi_{X, D^{-1}(A)} \circ c_{X,{}^*\!A} \circ c_{{}^*\!A,X}
  = c^A_{X, D^{-1}(A)} \circ c^A_{D^{-1}(A),X} \circ \pi_{X,D^{-1}(A)}. 
\end{gather*}
Thus, since $\pi_{X,{}^*\!A}$ is an epimorphism, the equation \eqref{eq:Mueger-centrality-of-K-1} is equivalent to
\begin{equation*}
  c^A_{X, D^{-1}(A)} \circ c^A_{D^{-1}(A),X} = \id_{D^{-1}(A)} \otimes_A \id_{X}.
\end{equation*}
Therefore (a) is equivalent to (b).
\end{proof}

\begin{corollary}
  \label{cor:local-module-ribbon}
  Suppose that $A$ is a commutative algebra such that $\C_A^{\loc}$ is rigid. Suppose moreover that $\theta_A = \id_A$ and $A$ belongs to the M\"uger center of $\C$. Then $\C_A^{\loc}$ ($=\C_A$) is a ribbon monoidal category.
\end{corollary}
\begin{proof}
As $\theta_A=\id_A$, Proposition~\ref{prop:K-in-loc} implies that $D(A) \in\C_A^{\loc}$. By the assumption that $A \in \C'$, we have $A^* \in \C'$. Thus, $D(A) \in (\C_A^{\loc})'$. Now, the claim follows by Theorem~\ref{thm:local-module-ribbon}.
\end{proof}


\section{Pivotality and the double dual in \texorpdfstring{$\C_A^{\sigma}$}{C-A-sigma}}\label{sec:double-dual-rigid}

\subsection{Main result}

Throughout this section, $\mathcal{C}$ is a rigid monoidal category admitting coequalizers and $\mathbf{A} = (A, \sigma)$ is a commutative algebra in the Drinfeld center $\mathcal{Z}(\mathcal{C})$ with multiplication $\mu$, unit $\eta$ and half-braiding $\sigma$.
We assume moreover that ${}_A\C_A$ is rigid.
Then, as we have proved in Theorem~\ref{thm:ACA-rigid-consequences}, $\mathcal{C}_A^{\sigma}$ is a rigid monoidal category with left duality functor $(-)^{\dagger} := \uHom_A(-, A)$.
The aim of this section is to give a convenient description of the double dual functor $(-)^{\dagger\dagger}$ and then use it  to give a sufficient condition for the category $\C_A^{\sigma}$ to be pivotal.

Our result is described as follows: Let $\mathfrak{u}$ denote the Drinfeld isomorphism of $\Z(\C)$.
Since $\mathbf{A}$ is commutative, $u := \mathfrak{u}_{\mathbf{A}} : \mathbf{A} \to \mathbf{A}^{**}$ is in fact an isomorphism of algebras in $\Z(\C)$.
Given a morphism $f: R \to S$ of algebras in $\mathcal{C}$, we denote by $(-)_{(f)} : \mathcal{C}_S \to \mathcal{C}_R$ the functor induced by $f$. 

Recall from \S\ref{subsec:central-rigid-case}, the anti-autoequivalence $D=(-)^*$ on $\C_A^\sigma$. We consider the autoequivalence $G := D^2$ on $\C_A^{\sigma}$ considered in Lemma \ref{lem:GV-duality-on-CA-2}. According to that lemma, it is given by
\begin{equation}
\label{eq:def-auto-equiv-G}
G: \C_A^{\sigma} \to \C_A^{\sigma},
\quad G(M) = (M^{**})_{(u)} \quad (M \in \C_A^{\sigma}).
\end{equation}
We also define
$\mathtt{g}_{0} : A \to G(A)$ and
$\mathtt{g}_{M,N} : G(M) \otimes_A G(N) \to G(M \otimes_A N)$
for $M, N \in \mathcal{C}_A^{\sigma}$ to be the unique morphism in $\mathcal{C}$ such that
\begin{equation}
  \label{eq:CA-double-dual-monoidal-structure-def-1}
  \mathtt{g}_{0} = u \quad \text{and} \quad
  \mathtt{g}_{M,N} \circ \pi_{M^{**}, N^{**}} = \pi_{M,N}^{**},
\end{equation}
where $\pi_{X, Y} : X \otimes Y \to X \otimes_A Y$ for $X \in \C_A$ and $Y \in {}_A\C$ is the canonical epimorphism. 

\begin{lemma}
  \label{lem:CA-double-dual-monoidal-structure}
  The morphisms defined by \eqref{eq:CA-double-dual-monoidal-structure-def-1} make $G$ a strong monoidal functor.
\end{lemma}
\begin{proof}
  By the definition of the action of $A$ on $A^{**}$, $\mathtt{g}_0$ is $A$-linear.
  The balancing property of $\pi_{M,N}$ implies that $\pi_{M,N}^{**}$ is a morphism of $A^{**}$-bimodules satisfying
  \begin{equation*}
    \pi_{M,N}^{**}((a_M^l)^{**} \otimes \id_{N^{**}})
    = \pi_{M,N}^{**}(\id_{M^{**}} \otimes (a_N^r)^{**}).
  \end{equation*}
  From this, we see that $\mathtt{g}_{M,N}$ is a well-defined $A$-linear morphism.
  It is easy to verify that $\mathtt{g}$ and $\mathtt{g}_0$ satisfy the axioms of a monoidal functor.

The map $\mathtt{g}_0$ is an isomorphism because it is a component of the Drinfeld isomorphism.
  As $(-)^{**}$ is an autoequivalence on $\C$, $\pi_{M,N}^{**}$ is also a coequalizer. As $\pi_{M^{**}, N^{**}}$ is also a coequalizer of the same maps, by the uniqueness of coequalizers, $\mathtt{g}_{M,N}$ is an isomorphism. This proves that $G$ is strong monoidal.
\end{proof}

Our main result in this section is:

\begin{theorem}
  \label{thm:CA-double-dual}
  There is an isomorphism of monoidal functors
  \begin{equation*}
    \phi_M : G(M) \to M^{\dagger\dagger}
    \quad (M \in \mathcal{C}_A^{\sigma}).
  \end{equation*}
\end{theorem}

An explicit description of the morphism $\phi_M$ will be given in Lemma~\ref{lem:CA-double-dual-as-functor} in terms of some natural isomorphisms introduced later. We will use this theorem to prove:

\begin{theorem}
  \label{thm:pivotality-C-A-sigma}
  Suppose, moreover, that $\C$ has a pivotal structure $\mathfrak{p}$ such that $\mathfrak{u}_{\mathbf{A}} = \mathfrak{p}_A$.
  Then the monoidal category $\mathcal{C}_A^{\sigma}$ has a pivotal structure given by
  \begin{equation}
    \label{eq:CA-pivotal-structure}
    M \xrightarrow{\quad \mathfrak{p}_M \quad}
    M^{**} \xrightarrow{\quad \phi_M \quad} M^{\dagger\dagger}
    \quad (M \in \mathcal{C}_A^{\sigma}).
  \end{equation}
\end{theorem}
\begin{proof}
  The assumption $\mathfrak{u}_{\mathbf{A}} = \mathfrak{p}_A$ guarantees that $\mathfrak{p}_M : M \to G(M)$ is a morphism in $\mathcal{C}_A^{\sigma}$. Moreover, we have
  $\mathfrak{p}_A = \mathtt{g}_0$
  and $\mathfrak{p}_{M \otimes_A N}
  = \mathtt{g}_{M,N} \circ (\mathfrak{p}_M \otimes_A \mathfrak{p}_N)$
  for all $M, N \in \mathcal{C}_A^{\sigma}$ by \eqref{eq:CA-double-dual-monoidal-structure-def-1} and the definition of a pivotal structure.
  Thus \eqref{eq:CA-pivotal-structure} is an isomorphism of monoidal functors.
\end{proof}

As a special case, we can consider the setting where $\mathcal{C}$ is, in addition, braided with braiding $c$, and $A$ is a commutative algebra in $\mathcal{C}$.
There are two ways to make $A$ a commutative algebra in $\mathcal{Z}(\mathcal{C})$.
For $X \in \mathcal{C}$, we put
\begin{equation*}
  \sigma^{+}_X = c_{A,X} \quad \text{and} \quad \sigma^{-}_X = c_{X,A}^{-1}
\end{equation*}
and then define $\mathbf{A}^{\pm} := (A, \sigma^{\pm})$.

\begin{theorem}
  \label{thm:pivotality-C-A-loc}
  If $\mathfrak{u}_{\mathbf{A}^{+}} = \mathfrak{p}_A$ or $\mathfrak{u}_{\mathbf{A}^{-}} =\mathfrak{p}_A$, then $\mathcal{C}_A^{\mathrm{loc}}$ is a braided pivotal category.
\end{theorem}
\begin{proof}
If $\mathfrak{u}_{\mathbf{A}^{+}} = \mathfrak{p}_A$ (resp., $\mathfrak{u}_{\mathbf{A}^{-}} =\mathfrak{p}_A$), then by Theorem~\ref{thm:pivotality-C-A-sigma}, the category $\C_{A}^{\sigma^+}$ (resp., $\C_{A}^{\sigma^-}$) is pivotal. Thus, its full subcategory $\C_A^{\loc}$ is pivotal as well.
\end{proof}


The rest of this section is organized as follows:
after recalling some basic facts in Sections~\ref{subsec:bimod-duals-basics} and \ref{subsec:module-structure-of-dagger}, we prove Theorem \ref{thm:CA-double-dual} in three steps in the subsequent subsections. In Section~\ref{subsec:step-1}, we establish the isomorphism $\phi : G \to (-)^{\dagger\dagger}$ at the level of functors. Section~\ref{subsec:step-2} explains the monoidal structure of $(-)^{\dagger\dagger}$. Finally, in Section~\ref{subsec:step-3}, we prove Theorem~\ref{thm:CA-double-dual}. 


\subsection{Basics on the internal Hom functor}\label{subsec:bimod-duals-basics}

As we have discussed in Section \ref{subsec:closed-bimodules}, the monoidal category ${}_A\mathcal{C}_A$ is closed with the internal Hom functor $\uHom_A$. The unit and the counit of the adjunction $(-) \otimes_A M \dashv \uHom_A(M, -)$ for $M \in {}_A\mathcal{C}_A$ are morphisms in $\mathcal{C}$ characterized by the following equations: For $N \in {}_A\C_A$,
\begin{gather}
  \label{eq:internal-Hom-def-coev}
  i_{M, N \otimes_A M} \circ \ucoev_{M,N}
  = (\pi_{N, M} \otimes \id_{M^*}) \circ (\id_N \otimes \coev_M), \\
  \label{eq:internal-Hom-def-eval}
  \underline{\ev}_{M,N} \circ \pi_{\uHom_A(M, N), M}
  = (\id_N \otimes \ev_M) \circ (i_{M,N} \otimes \id_M).
\end{gather}

Since the functor $(-) \otimes_A A : {}_A\mathcal{C}_A \to {}_A\mathcal{C}_A$ is isomorphic to the identity functor, its right adjoint $\uHom_A(A, -)$ is also isomorphic to the identity functor.
We note how this isomorphism looks if we realize the internal Hom functor by \eqref{eq:internal-Hom-rigid-case}.

\begin{lemma}
  \label{lem:internal-Hom-iso-alpha}
  For $M \in {}_A\mathcal{C}_A$, there is a natural isomorphism
  \begin{equation*}
    \alpha_M : M \to \uHom_A(A, M)
  \end{equation*}
  of $A$-bimodules in $\mathcal{C}$ determined by
  \begin{equation}
    \label{eq:internal-Hom-iso-alpha}
    i_{A,M} \circ \alpha_M = (a_M^{r} \otimes \id_{A^*}) \circ (\id_M \otimes \coev_A).
  \end{equation}
  The inverse of $\alpha_M$ is given by
  \begin{equation*}
    \alpha_M^{-1} = (\id_M \otimes \eta^*) \circ i_{A,M},
  \end{equation*}
  where $\eta : \unit \to A$ is the unit of A.
\end{lemma}
\begin{proof}
  Straightforward.
\end{proof}

In a rigid monoidal category with left duality functor $(-)^*$, there is a natural isomorphism $(X \otimes Y)^* \cong Y^* \otimes X^*$ for all objects $X$ and $Y$. We consider this natural isomorphism in the rigid monoidal category $({}_A\C_A, \otimes_A, A)$. Namely, we have a natural isomorphism 
\begin{equation}
  \label{eq:gamma-def}
  \gamma_{M,N} : M^{\dagger} \otimes_A N^{\dagger} \to (N \otimes_A M)^{\dagger}
  \quad (M, N \in {}_A\mathcal{C}_A)
\end{equation}
characterized by the equation
\begin{equation}
  \label{eq:internal-Hom-iso-gamma-def}
  \underline{\ev}_{N \otimes_A M, A} \circ (\gamma_{M,N} \otimes_A \id_{N} \otimes_A \id_{M})
  = \underline{\ev}_{M,A} \circ (\id_{M^{\dagger}} \otimes_A \underline{\ev}_{N,A} \otimes_A \id_{M}),
\end{equation}
where we have identified $M^{\dagger} \otimes_A A$ with $M^{\dagger}$ in the right hand side. An alternative description of the natural isomorphism $\gamma$ is:
\begin{lemma}
  \label{lem:internal-Hom-iso-gamma-eq-1}
  For $M, N \in {}_A\mathcal{C}_A$, we have
  \begin{align}
    \label{eq:internal-Hom-iso-gamma-eq-1}
    i_{N \otimes M, A}
    \circ \pi_{M,N}^{\dagger} \circ \gamma_{M, N} \circ \pi_{M^{\dagger}, N^{\dagger}}
    = (\id_A \otimes a^r_{M^*} \otimes \id_{N^*}) \circ (i_{M, A} \otimes i_{N, A}).
  \end{align}
\end{lemma}
\begin{proof}
  Given $X_1, \cdots, X_n \in {}_A\C_A$, we denote by
  \begin{equation*}
    \pi_{X_1, \cdots, X_n}^{(n)} : X_1 \otimes \cdots \otimes X_n \to X_1 \otimes_A \cdots \otimes_A X_n
  \end{equation*}
  the epimorphism obtained by using $\pi_{-,-}$ iteratively (since the tensor product of $\C$ preserves coequalizers, the morphism $\pi_{X_1, \cdots, X_n}^{(n)}$ does not depend on the order of taking iterative quotients). Now let $L$ and $R$ be the left and the right hand side of \eqref{eq:internal-Hom-iso-gamma-eq-1}, respectively.
  We compute
  \begin{align*}
    & (\id_A \otimes \ev_{N \otimes M}) \circ (L \otimes \id_N \otimes \id_M) \\
    {}^{\text{(nat.)}}
    & = (\id_A \otimes \ev_{N \otimes M}) \circ ((\id_A \otimes \pi_{M,N}^*)i_{N \otimes_A M,A} \gamma_{M,N} \pi_{M^{\dagger},N^{\dagger}} \otimes \id_{N \otimes M}) \\
    {}^{\text{(nat.)}}
    & = (\id_A \otimes \ev_{N \otimes_A M}) \circ (i_{N \otimes_A M,A} \gamma_{M,N} \otimes \id_{N \otimes_A M}) \circ (\pi_{M^{\dagger},N^{\dagger}} \otimes \pi_{N,M}) \\
    {}^{\eqref{eq:internal-Hom-def-eval}}
    & = \underline{\ev}_{N \otimes_A M, A}
      \circ \pi_{(N \otimes_A M)^{\dagger}, N \otimes_A M}
      \circ (\gamma_{M,N} \otimes \id_{N \otimes_A M})
      \circ (\pi_{M^{\dagger}, N^{\dagger}} \otimes \pi_{N,M}) \\
    {}^{\text{(nat.)}}
    & = \underline{\ev}_{N \otimes_A M, A}
    \circ (\gamma_{M,N} \otimes_A \id_{N} \otimes_A \id_{M})
      \circ \pi_{M^{\dagger}, N^{\dagger}, N, M}^{(4)} \\
    {}^{\eqref{eq:internal-Hom-iso-gamma-def}}
    & = \underline{\ev}_{M,A}
      \circ (\id_{M^{\dagger}} \otimes_A \underline{\ev}_{N,A} \otimes_A \id_{M})
      \circ \pi_{M^{\dagger}, N^{\dagger}, N, M}^{(4)} \\
    {}^{(*)}
    & = \underline{\ev}_{M,A}
      \circ (\id_{M^{\dagger}} \otimes_A \underline{\ev}_{N,A} \pi_{N^{\dagger},N} \otimes_A \id_{M})
      \circ \pi_{M^{\dagger}, N^{\dagger} \otimes N, M}^{(3)} \\
    {}^{\text{(nat.)}}
    & = \underline{\ev}_{M,A} \circ \pi_{M^{\dagger}, A, M}^{(3)}
      \circ (\id_{M^{\dagger}} \otimes \underline{\ev}_{N,A} \pi_{N^{\dagger},N} \otimes \id_{M})\\
    {}^{\eqref{eq:internal-Hom-def-eval}}
    & = \underline{\ev}_{M,A}
      \circ \pi_{M^{\dagger}, A, M}^{(3)}
      \circ (\id_{M^{\dagger}} \otimes \id_A \otimes \ev_N \otimes \id_M)
      \circ (\id_{M^{\dagger}} \otimes i_{N,A} \otimes \id_N \otimes \id_M),
  \end{align*}
  where `(nat.)' follows from the (di)naturality of $i_{P,A}$, $\ev_X$ and $\pi_{P,Q}$ in $P,Q \in {}_A\C_A$ and $X \in \C$, and $(*)$ follows from the definition of $\pi^{(n)}_{X_1, \cdots, X_n}$. Moreover, the canonical isomorphism $M^{\dagger} \otimes_A A \cong M^{\dagger}$ is implicit in the last four expressions. Noticing that the canonical isomorphisms $A \otimes_A X \cong X \cong X \otimes_A A$ for $X \in {}_A\C_A$ are induced by actions of $A$ on $X$, we have
  \begin{equation}
    \label{eq:internal-Hom-iso-gamma-eq-1-proof-1}
      \pi_{X,Y}(\id_X \otimes a_Y^l)
      = \pi_{X, A, Y}^{(3)}
      = \pi_{X,Y}(a_X^r \otimes \id_Y)
      \quad (X, Y \in {}_A\C_A)
  \end{equation}
  when we identify $X \otimes_A A \otimes_A Y = X \otimes_A Y$. Thus we continue the above computation as follows:
  \begin{align*}
    ...
    & = \underline{\ev}_{M,A}
      \circ \pi_{M^{\dagger}, M}
      \circ (a_{M^{\dagger}}^r \otimes \ev_N \otimes \id_M)
      \circ (\id_{M^{\dagger}} \otimes i_{N,A} \otimes \id_N \otimes \id_M) \\
    {}^{\eqref{eq:internal-Hom-def-eval}}
    & = (\id_A \otimes \ev_M)
      \circ (i_{M,A} a_{M^{\dagger}}^r \otimes \ev_N \otimes \id_M)
      \circ (\id_{M^{\dagger}} \otimes i_{N,A} \otimes \id_N \otimes \id_M) \\
    {}^{(*)}
    & = (\id_A \otimes \ev_M)
      \circ ((\id_A \otimes a_{M^*}^r)i_{M,A} \otimes \ev_N \otimes \id_M)
      \circ (\id_{M^{\dagger}} \otimes i_{N,A} \otimes \id_N \otimes \id_M) \\
    & = (\id_A \otimes \ev_{N \otimes M}) \circ (R \otimes \id_N \otimes \id_M),
  \end{align*}
  where the equality $(*)$ follows from that $i_{M,A} : M^{\dagger} \to A \otimes M^*$ is a morphism of right $A$-modules in $\mathcal{C}$. Therefore \eqref{eq:internal-Hom-iso-gamma-eq-1} is verified.
\end{proof}


\subsection{Module structure of the duality functor}
\label{subsec:module-structure-of-dagger}

We note that $\C_A^{\sigma}$ is a $\C$-bimodule category by the actions as stated by the following lemma:

\begin{lemma}
  \label{lem:C-A-sigma-bimodule-structure}
  Let $X \in \mathcal{C}$ and $M \in \mathcal{C}_A^{\sigma}$ be arbitrary objects.
  \begin{enumerate}
  \item $X \otimes M \in \mathcal{C}_A^{\sigma}$ by the actions
    \begin{equation*}
      a_{X \otimes M}^{l} = (\id_X \otimes a_M^{l}) \circ (\sigma_X \otimes \id_M), \quad
      a_{X \otimes M}^{r} = \id_X \otimes a_M^r.
    \end{equation*}
  \item $M \otimes X \in \mathcal{C}_A^{\sigma}$ by the actions
    \begin{equation*}
      a_{M \otimes X}^{l} = a_M^{l} \otimes \id_X, \quad
      a_{M \otimes X}^{r} = (a_M^{r} \otimes \id_X) \circ (\id_M \otimes \sigma_X^{-1}).
    \end{equation*}
  \end{enumerate}
\end{lemma}
\begin{proof}
  Straightforward.
\end{proof}

The actions of $\mathcal{C}$ on $\mathcal{C}_A^{\sigma}$ are compatible with the internal Hom functor in the following sense:

\begin{lemma}
  \label{lem:internal-Hom-module-structure-1}
  For $X, Y \in \mathcal{C}$ and $M, N \in \mathcal{C}_A^{\sigma}$, we have
  \begin{equation}
    \label{eq:internal-Hom-module-structure-1}  
    X \otimes \uHom_A(M, N) \otimes Y^*
    = \uHom_A(Y \otimes M, X \otimes N)
  \end{equation}
  as right $A$-modules.
\end{lemma}
\begin{proof}
  By the definition of the right action of $A$ on $Y \otimes M$ and $X \otimes N$ given in Lemma~\ref{lem:C-A-sigma-bimodule-structure}, it is easy to see that \eqref{eq:internal-Hom-module-structure-1} holds as an equation of subobjects of $X \otimes N \otimes M^* \otimes Y^*$. Again by the definition of the action given by Lemma~\ref{lem:C-A-sigma-bimodule-structure}, one can verify that the actions of $A$ on the left and the right hand side of \eqref{eq:internal-Hom-module-structure-1} agree. The proof is done.
\end{proof}

For $X \in \mathcal{C}$ and $M \in \mathcal{C}_A^{\sigma}$, we define
\begin{equation*}
  \xi_{X,M} : X^* \otimes M^{\dagger} \to (M \otimes X)^{\dagger}
\end{equation*}
to be the unique morphism in $\mathcal{C}$ such that
\begin{equation}
  \label{eq:internal-Hom-iso-xi-def}
  i_{X \otimes M, A} \circ \xi_{X,M} = (\sigma_{X^*}^{-1} \otimes \id_{M^*}) \circ (\id_{X^*} \otimes i_{M,A}).
\end{equation}
The following lemma says that the left duality functor $(-)^{\dagger} : \C_A^{\sigma} \to (\C_A^{\sigma})^{\op}$ is a right $\C$-module functor if we view the target as a right $\C$-module category by the action $\triangleright$ given by $X \triangleright M = X^* \otimes M$.

\begin{lemma}
  \label{lem:internal-Hom-module-structure-2}
  $\xi$ is a natural isomorphism of right $A$-modules satisfying
  \begin{equation}
    \label{eq:internal-Hom-iso-xi-eq-1}
    \xi_{\unit, M} = \id_{M^{\dagger}}
    \quad \text{and} \quad
    \xi_{X \otimes Y, M} = \xi_{Y, M \otimes X} \circ (\id_{Y^*} \otimes \xi_{X,M})
  \end{equation}
  for all $X, Y \in \mathcal{C}$ and $M \in \mathcal{C}_A^{\sigma}$.
\end{lemma}
\begin{proof}
  To prove this lemma, we introduce the following notation:
  \begin{equation*}
      \tau_{\mathbf{Z}, V} := z_V : Z \otimes V \to V \otimes Z
      \quad (\mathbf{Z} = (Z, z) \in \Z(\C), V \in \C)
  \end{equation*}
  Thus, for $\mathbf{A} = (A, \sigma)$, we have $\tau_{\mathbf{A},V} = \sigma_V$ for $V \in \C$.
  According to \cite[Section 7.13]{etingof2016tensor}, a left dual object of $\mathbf{A}$ has $A^*$ as its underlying object, and the half-braiding is given by $\tau_{\mathbf{A}^*, V} = (\sigma_{{}^*V}^{-1})^*$ for $V \in \C$. Thus we have
  \begin{equation}
    \label{eq:lem:internal-Hom-module-structure-2-proof-1}
      (a_{M \otimes X}^r)^*
      = (\id_M \otimes \sigma_X^{-1})^* \circ (a_M^r \otimes \id_X)^*
      = (\tau_{\mathbf{A}^*, X^*} \otimes \id_{M^*}) \circ (\id_{X^*} \otimes (a_M^r)^*),
  \end{equation}
  see Lemma \ref{lem:C-A-sigma-bimodule-structure} (b) for the definition of $a_{M \otimes X}^r$.

  We first verify that $\xi_{X,M}$ is well-defined. Let $R$ be the right hand side of \eqref{eq:internal-Hom-iso-xi-def}, and set
  \begin{equation}
    \label{eq:def-delta-A}
    \delta_A = (\mu \otimes \id_{A^*}) \circ (\id_A \otimes \coev_A) : A \to A \otimes A^*,
  \end{equation}
  where $\mu$ is the multiplication of $A$.
  By the expression \eqref{eq:iHom-A-in-ACA} of $\uHom_A$ as an equalizer, to show that $\xi_{X,M}$ is well-defined, we shall show:
  \begin{equation*}
      (\delta_A \otimes \id_{X^*} \otimes \id_{M^*}) \circ R
      = (\id_A \otimes (a_{M \otimes X}^r)^*) \circ R.
  \end{equation*}
  Since $\mu$ and $\coev_A$ are morphisms in $\Z(\C)$, $\delta_A$ is also a morphism in $\Z(\C)$. Thus we have
  \begin{align*}
    & (\delta_A \otimes \id_{X^*} \otimes \id_{M^*}) \circ R \\
    & = (\delta_A \otimes \id_{X^*} \otimes \id_{M^*})
      \circ (\tau_{\mathbf{A}, X^*}^{-1} \otimes \id_{M^*}) \circ (\id_{X^*} \otimes i_{M,A}) \\
    & = (\tau_{\mathbf{A} \otimes \mathbf{A}^*, X^*}^{-1} \otimes \id_{M^*})
      \circ (\id_{X^*} \otimes \delta_A \otimes \id_{M^*})
      \circ (\id_{X^*} \otimes i_{M,A}) \\
    {}^{\eqref{eq:internal-Hom-iso-xi-def}}
    & = (\tau_{\mathbf{A} \otimes \mathbf{A}^*, X^*}^{-1} \otimes \id_{M^*})
      \circ (\id_{X^*} \otimes \id_A \otimes (a_M^r)^*)
      \circ (\id_{X^*} \otimes i_{M,A}) \\
    & = ((\id_A \otimes \tau_{\mathbf{A}^*, X^*}^{-1}) (\tau_{\mathbf{A}, X^*}^{-1} \otimes \id_{A}) \otimes \id_{M^*})
      \circ (\id_{X^*} \otimes \id_A \otimes (a_M^r)^*)
      \circ (\id_{X^*} \otimes i_{M,A}) \\
    & = (\id_A \otimes \tau_{\mathbf{A}^*, X^*}^{-1} \otimes \id_{M^*})
      \circ (\id_A \otimes \id_{X^*} \otimes (a_M^r)^*)
      \circ (\tau_{\mathbf{A},X^*}^{-1} \otimes \id_{M^*})
      \circ (\id_{X^*} \otimes i_{M,A}) \\
    {}^{\eqref{eq:lem:internal-Hom-module-structure-2-proof-1}}
    & = (\id_A \otimes (a_{M \otimes X}^r)^*) \circ R,
  \end{align*}
  where the second equality follows from that $\delta_A$ is a morphism in $\Z(\C)$ and the third from the equalizer property of $i_{M,A}$.
  Therefore the morphism $\xi_{X, M}$ is well-defined.

  Since the right hand side of \eqref{eq:internal-Hom-iso-xi-def} is $A$-linear, so is $\xi_{X,M}$.
  The inverse of $\xi_{X,M}$ is given as the unique morphism satisfying
  \begin{equation*}
    (\id_{X^*} \otimes i_{M,A}) \circ \xi_{X,M}^{-1}
    = (\sigma_{X^*} \otimes \id_{M^*}) \circ i_{M \otimes X, A}.
  \end{equation*}
  Equation \eqref{eq:internal-Hom-iso-xi-eq-1} follows from the axiom of half-braidings. The proof is done.
\end{proof}


\subsection{Step 1: isomorphism of functors}\label{subsec:step-1}
Suppose that $R$ and $S$ are algebras in $\mathcal{C}$. Then the assignment $M \mapsto (-) \otimes_R M$ gives an equivalence from ${}_R\C_S$ to the category of left $\C$-module functors $\C_R \to \C_S$ preserving coequalizers \cite[Theorem~4.2]{pareigis1977non}. We refer to this fact as {\em the Eilenberg-Watts theorem}. We note that given an object $F : \C_R \to \C_S$ in the latter category, the corresponding $R$-$S$-bimodule is the right $S$-module $F(R) \in \C_S$ equipped with the left $R$-module structure given by $F(\mu_R) \circ f_{R,R} : R \otimes F(R) \to F(R)$, where $f_{X,M} : X \otimes F(M) \to F(X \otimes M)$ is the module structure of $F$ and $\mu_R : R \otimes R \to R$ is the multiplication of $R$. 

We construct a natural isomorphism $M^{\dagger\dagger} \cong G(M)$ for $M \in \C_A^{\sigma}$, where $G$ is the functor given by \eqref{eq:def-auto-equiv-G}. The Eilenberg-Watts theorem is essential for this purpose.
By Lemma~\ref{lem:internal-Hom-module-structure-1}, we have
\begin{equation*} 
  M^{\dagger} \otimes X^* = \uHom_A(M, A) \otimes X^* = \uHom_A(X \otimes M, A) = (X \otimes M)^{\dagger} 
\end{equation*}
for $X \in \mathcal{C}$ and $M \in \mathcal{C}_A^{\sigma}$.
Hence we have a natural isomorphism
\begin{equation*}
  \xi_{X^{*}, M^{\dagger}} : X^{**} \otimes M^{\dagger\dagger}
  \to (M^{\dagger} \otimes X^*)^{\dagger} = (X \otimes M)^{\dagger\dagger}
\end{equation*}
for $X \in \mathcal{C}$ and $M \in \mathcal{C}_A^{\sigma}$, where $\xi$ is the natural isomorphism discussed in Lemma~\ref{lem:internal-Hom-module-structure-2}.
Equation \eqref{eq:internal-Hom-iso-xi-eq-1} means that the double dual functor $(-)^{\dagger\dagger}$ on $\mathcal{C}_A^{\sigma}$ is a `$(-)^{**}$-twisted' left $\mathcal{C}$-module autoequivalence on $\mathcal{C}_A^{\sigma}$. To `untwist' the double dual functor $(-)^{**}$, we consider the functor
\begin{equation*}
  F: \mathcal{C}_{A^{**}}
  \xrightarrow{\quad {}^{**}(-) \quad}
  \mathcal{C}_A^{\sigma}
  \xrightarrow{\quad (-)^{\dagger\dagger} \quad}
  \mathcal{C}_A^{\sigma},
\end{equation*}
which is an equivalence of $\mathcal{C}$-module categories with the module structure given by
\begin{equation}
\label{eq:isomorphism-xi-twisted}
  \xi_{X,M}' : X \otimes F(M)
  = X \otimes ({}^{**}M)^{\dagger\dagger}
  \xrightarrow{\quad \xi_{{}^* \! X, ({}^{**}M)^{\dagger}} \quad}
  ({}^{**}X \otimes {}^{**}M)^{\dagger\dagger}
  = F(X \otimes M)
\end{equation}
for $X \in \C$ and $M \in \C_{A^{**}}$.
Thus, by the Eilenberg-Watts theorem, $F(A^{**})$ has a canonical structure of $A^{**}$-$A$-bimodule and there is a natural isomorphism
\begin{equation*}
  ({}^{**}M)^{\dagger\dagger} = F(M) \cong M \otimes_{A^{**}} F(A^{**})
  \quad (M \in \mathcal{C}_{A^{**}}).
\end{equation*}
Thus we obtain the following description of the double dual functor:
\begin{equation*}
  M^{\dagger\dagger} = F(M^{**}) \cong M^{**} \otimes_{A^{**}} F(A^{**})
  \quad (M \in \mathcal{C}_A^{\sigma}).
\end{equation*}
The problem is what the bimodule $F(A^{**})$ is. We define
\begin{equation*}
  \tilde{\alpha} := \uHom_A(\alpha_A^{-1}, A) \circ \alpha_A
  : A \to \uHom_A(\uHom_A(A,A),A) = F(A^{**}),
\end{equation*}
where $\alpha_M$ for $M \in {}_A\mathcal{C}_A$ is the natural isomorphism discussed in Lemma \ref{lem:internal-Hom-iso-alpha}.
We transport the $A^{**}$-$A$-bimodule structure of $F(A^{**})$ to $A$ via the isomorphism $\tilde{\alpha}$ and let $\tilde{a}{}_A^{l}$ and $\tilde{a}{}_A^r$ be the resulting left and right action of $A^{**}$ and $A$ on $A$, respectively.

\begin{lemma}
  \label{lem:actions-on-F(A**)}
  The actions on $A$ are given by
  \begin{equation*}
    \tilde{a}{}_A^{l} = \mu \circ (u^{-1} \otimes \id_A)
    \quad \text{and} \quad
    \tilde{a}{}_A^r = \mu,
  \end{equation*}
  where $u = \mathfrak{u}_{\mathbf{A}}$ is the component of the Drinfeld isomorphism for $\mathbf{A} = (A, \sigma)$.
\end{lemma}
\begin{proof}
  The formula for $\tilde{a}{}_A^r$ follows from that $\alpha_M$ for $M \in \mathcal{C}_A^{\sigma}$ is actually an isomorphism of right $A$-modules. We prove the formula for $\tilde{a}{}_A^l$.
  In the proof of this lemma, we write $[M,N] = \uHom_A(M, N)$ to save space. By Lemma~\ref{lem:internal-Hom-iso-alpha}, we have
  \begin{equation}
    \label{eq:lem:actions-on-F(A**)-proof-eq-1}
    \begin{gathered}
      i_{[A,A],A} \circ \tilde{\alpha}
      = i_{[A,A],A} \circ [\alpha_A^{-1}, A] \circ \alpha_A
      = (\id_A \otimes (\alpha_A^{-1})^*) \circ i_{A,A} \circ \alpha_A \\
      = (\id_A \otimes i_{A,A}^*) \circ (\id_A \otimes \eta^{**} \otimes \id_{A^*}) \circ \delta_A,
    \end{gathered}
  \end{equation}
  where $\eta : \unit \to A$ is the unit of $A$ and $\delta_A$ is given by \eqref{eq:def-delta-A}. 
  Now we consider the commutative diagram of Figure~\ref{fig:lem:actions-on-F(A**)-1}.
  By the unit law of the action, the right column is equal to $\sigma_{M^*}^{-1} \otimes \id_{A^*}$.
  By definition, the morphism $\tilde{\alpha} \circ \tilde{a}{}_A^{l}$ is equal to the left column of the diagram with $M = A$.
  Thus we obtain
  \begin{gather}
    \label{eq:lem:actions-on-F(A**)-proof-eq-2}
    i_{[A, A], A} \circ \tilde{\alpha} \circ \tilde{a}{}_A^{l}
    = (\id_A \otimes i_{A,A}^*) \circ (\sigma_{A^{**}}^{-1} \otimes \id_{A^*}) \circ (\id_{A^{**}} \otimes \delta_A).
  \end{gather}
  Let $b := (\ev_{A^*} \otimes \id_{A^*}) \circ (\id_{A^{**}} \otimes \mu^{*}) : A^{**} \otimes A^* \to A^*$ be the left action of $A^{**}$ on $A^*$.
  By the balancing property of $\pi$ and the unit law,
  \begin{equation}
    \label{eq:lem:actions-on-F(A**)-proof-eq-3}
    i_{A,A}^*
    = i_{A,A}^* \circ (\mu^{**} \otimes \id_{A^{*}}) \circ (\eta^{**} \otimes \id_{A^{**}} \otimes \id_{A^*}) 
    = i_{A,A}^* \circ (\eta^{**} \otimes b).    
  \end{equation}
  By the graphical calculus as in Figure \ref{fig:lem:actions-on-F(A**)-2}, we have
  \begin{equation}
    \label{eq:lem:actions-on-F(A**)-proof-eq-4}
    (\id_A \otimes b) \circ (\sigma_{A^{**}}^{-1} \otimes \id_{A^*}) \circ (\id_{A^{**}} \otimes \delta_A)
    = \delta_A \circ \mu \circ (\overline{u} \otimes \id_A),
  \end{equation}
  where $\overline{u} = (\id_A \otimes \ev_{A^*}) (\sigma_{A^{**}}^{-1} \otimes \id_{A^*})(\id_{A^{**}} \otimes \coev_A) : A^{**} \to A$ (the second equality in the figure follows from the associativity and the commutativity of $\mu$).
  One can verify that $\overline{u}$ is the inverse of $u$.
  Hence, $i_{[A, A], A} \circ \tilde{\alpha} \circ \tilde{a}{}_A^{l}$
  \begin{align*}
    {}^{\eqref{eq:lem:actions-on-F(A**)-proof-eq-2},
    \eqref{eq:lem:actions-on-F(A**)-proof-eq-3}}
    & = (\id_A \otimes i_{A,A}^*) \circ
      (\id_A \otimes \eta^{**} \otimes b)
      \circ (\sigma_{A^{**}}^{-1} \otimes \id_{A^*}) \circ (\id_{A^{**}} \otimes \delta_A) \\
    {}^{\eqref{eq:lem:actions-on-F(A**)-proof-eq-4}}
    & = (\id_A \otimes i_{A,A}^*) \circ
      (\id_A \otimes \eta^{**} \otimes \id_{A^{*}})
      \circ \delta_A \circ \mu \circ (u^{-1} \otimes \id_A) \\
    {}^{\eqref{eq:lem:actions-on-F(A**)-proof-eq-1}}
    & =  i_{[A, A], A} \circ \tilde{\alpha} \circ \mu \circ (u^{-1} \otimes \id_A).
  \end{align*}
  The formula as stated follows since $i_{[A, A], A} \circ \tilde{\alpha}$ is monic.
\end{proof}

\begin{figure}
  \begin{equation*}
    \begin{tikzcd}[column sep=24pt, row sep=32pt]
      M^{**} \otimes A
      \arrow[rr, "{\id \otimes \delta_A}"]
      \arrow[d, "{\id \otimes \tilde{\alpha}}"']
      \arrow[rrd, phantom, "{\eqref{eq:lem:actions-on-F(A**)-proof-eq-1}}"]
      & & M^{**} \otimes A \otimes A^*
      \arrow[d, "{\id \otimes \id \otimes \eta^{**} \otimes \id}"]
      \\ 
      M^{**} \otimes [[A, A], A]
      \arrow[d, "{\xi_{M^*, [A,A]^{\dagger}}}"']
      \arrow[r, hook, "{\id \otimes i_{[A,A], A}}" {yshift = 5pt}]
      \arrow[rd, phantom, "{\scriptsize \text{(definition of $\xi$)}}"]
      & M^{**} \otimes A \otimes [A, A]^*
      \arrow[d, "{\sigma_{M^{*}}^{-1} \otimes \id \otimes \id}" {yshift = 5pt}]
      & M^{**} \otimes A \otimes A^{**} \otimes A^*
      \arrow[l, two heads, "{\id \otimes \id \otimes i_{A,A}^*}"' {yshift = 5pt}]
      \arrow[d, "{\sigma^{-1}_{M^{**}} \otimes \id \otimes \id \otimes \id}"]
      \\ 
      {} [[A, A] \otimes M^*, A]
      \arrow[d, equal]
      \arrow[r, hook, "{i_{[A,A] \otimes M^*, A}}" {yshift = 5pt}]
      & A \otimes M^{**} \otimes [A, A]^*
      \arrow[d, equal]
      & A \otimes M^{**} \otimes A^{**} \otimes A^*
      \arrow[l, two heads, "{\id \otimes \id \otimes i_{A,A}^*}"' {yshift = 5pt}]
      \arrow[d, equal]
      \\[-15pt] 
      {} [[M \otimes A, A], A]
      \arrow[d, "{[[a_M^r, A], A]}"']
      \arrow[r, hook, "{i_{[M \otimes A,A], A}}" {yshift = 5pt}]
      & A \otimes [M \otimes A, A]^*
      \arrow[d, "{\id \otimes [a_M^r, A]^*}" {yshift = 0pt}]
      & A \otimes M^{**} \otimes A^{**} \otimes A^{*}
      \arrow[l, two heads, "{\id \otimes i_{M \otimes A,A}^*}"' {yshift = 5pt}]
      \arrow[d, "{\id \otimes (a_M^r)^{**} \otimes \id}"]
      \\ 
      {} [[M, A], A]
      \arrow[r, hook, "{i_{[M,A], A}}" {yshift = 0pt}]
      & A \otimes [M, A]^*
      & A \otimes M^{**} \otimes A^*
      \arrow[l, two heads, "{\id \otimes i_{M, A}^*}"' {yshift = 0pt}]
    \end{tikzcd}
  \end{equation*}
  \caption{Proof of Equation~\eqref{eq:lem:actions-on-F(A**)-proof-eq-1}}
  \label{fig:lem:actions-on-F(A**)-1}

  \bigskip
    \begin{equation*}
    \tikzset{overstrand/.style={preaction={-,draw=white,line width=8pt}}}
    \begin{tikzpicture}[x = 11pt, y = 11pt, baseline = 0]
      \coordinate (S1) at (.5,3); \node at (S1) [above] {$A^{**}$};
      \coordinate (S2) at (2,3); \node at (S2) [above] {$A$};
      \coordinate (T1) at (1,-3); \node at (T1) [below] {$A$};
      \coordinate (T2) at (4,-3); \node at (T2) [below] {$A^{*}$};
      \draw (S2) -- ++(0, -1)
      to [out = -90, in = -90, looseness = 1.5] coordinate[midway](P1) ++(1,0)
      to [out = +90, in = +90, looseness = 1.5] ++(1,0)
      -- ++(0, -1) coordinate (P2);
      \draw (P1) [out = -90, in = 90] to (T1);
      \node at (P1) {$\bullet$};
      \node (B1) at (P2) [below, draw] {$\mu^*$};
      \draw ($(B1.south) + (+.25,0)$) [out = -90, in = +90] to (T2);
      \draw ($(B1.south) + (-.25,0)$) -- ++(0,-.1) coordinate (P2);
      \draw [overstrand] (P2) to [out = -90, in = -90, looseness = 1.5] (S1);
    \end{tikzpicture}
    = \begin{tikzpicture}[x = 11pt, y = 11pt, baseline = 0]
      \coordinate (S1) at (0,3); \node at (S1) [above] {$A^{**}$};
      \coordinate (S2) at (1.5,3); \node at (S2) [above] {$A$};
      \coordinate (T1) at (0,-3); \node at (T1) [below] {$A$};
      \coordinate (T2) at (4,-3); \node at (T2) [below] {$A^{*}$};
      \draw (T2) -- ++(0,4)
      to [out = +90, in = +90, looseness = 1.5] ++(-2,0)
      to [out = -90, in = -90, looseness = 1.8] coordinate[midway] (P1) ++(.66,0)
      to [out = +90, in = +90, looseness = 1.8] ++(.66,0) -- ++(0,-1.5) coordinate (P2);
      \draw (P1) to [out = -90, in = -90, looseness = 1.8] coordinate[midway] (P3) ++(-1,0)
      to [out = +90, in = -90, looseness = 1.8] (S2);
      \draw (P3) to [out = -90, in = +90, looseness = 1] (T1);
      \node at (P1) {$\bullet$};
      \node at (P3) {$\bullet$};
      \draw [overstrand] (P2)
      to [out = -90, in = -90, looseness = 1.4] ($(S1)-(0,2)$) -- (S1);
    \end{tikzpicture}
    = \begin{tikzpicture}[x = 11pt, y = 11pt, baseline = 0]
      \coordinate (S1) at (.5,3); \node at (S1) [above] {$A^{**}$};
      \coordinate (S2) at (2,3); \node at (S2) [above] {$A$};
      \coordinate (T1) at (1,-3); \node at (T1) [below] {$A$};
      \coordinate (T2) at (4.7,-3); \node at (T2) [below] {$A^{*}$};
      \draw (S2) to [out = -90, in = +90, looseness = .8] ++(1,-2)
      to [out = -90, in = -90, looseness = 1.5] coordinate[midway] (P1) ++(-1,0) coordinate (P5);
      \draw (P1) to [out = -90, in = -90, looseness = 1.8] coordinate[midway] (P2) ++(1,0)
      -- ++(0,2) coordinate (P3)
      to [out = +90, in = +90, looseness = 1.8] ++(2,0) coordinate (P4)
      to [out = -90, in = 90] (T2);
      \node at (P1) {$\bullet$};
      \node at (P2) {$\bullet$};
      \draw (P2) to [out = -90, in = 90] (T1);
      \draw [overstrand] (P5) to [out = 90, in = -90] ($(P3)+(.66,0)$) coordinate (P6);
      \draw (P6) to [out = 90, in = +90, looseness = 1.8] ($(P4)-(.66,0)$) coordinate (P7);
      \draw [overstrand] (P7) to [out = -90, in = -90, looseness = 3.5] (S1);
    \end{tikzpicture}
    = \begin{tikzpicture}[x = 11pt, y = 11pt, baseline = 0]
      \coordinate (S1) at (0,3); \node at (S1) [above] {$A^{**}$};
      \coordinate (S2) at (2,3); \node at (S2) [above] {$A$};
      \coordinate (T1) at (1.75,-3); \node at (T1) [below] {$A$};
      \coordinate (T2) at (4,-3); \node at (T2) [below] {$A^{*}$};
      \draw (S2) -- ++(0, -2)
      to [out = -90, in = -90, looseness = 2] coordinate[midway] (P1) ++(-2,0)
      to [out = 90, in = 90, looseness=2] ++(1,.7) coordinate (P2);
      \draw [overstrand] (S1) -- ++(0,-.7)
      to [out = -90, in = -90, looseness=2] (P2);
      \draw (P1) -- ++(0,-.75)
      to [out = -90, in = -90, looseness = 2] coordinate[midway] (P3) ++(1.5,0)
      to [out = +90, in = +90, looseness = 2] ++(1.5,0)
      to [out = -90, in = +90] (T2);
      \draw (P3)
      to [out = -90, in = +90] (T1);
      \node at (P1) {$\bullet$};
      \node at (P3) {$\bullet$};
    \end{tikzpicture}
  \end{equation*}
  \caption{Proof of Equation~\eqref{eq:lem:actions-on-F(A**)-proof-eq-3}}
  \label{fig:lem:actions-on-F(A**)-2}
\end{figure}

\begin{lemma}
  \label{lem:CA-double-dual-as-functor}
  There is a natural isomorphism
  \begin{equation}
    \label{eq:CA-double-dual-iso-phi-def}
    \phi_M := (a_M^r)^{\dagger\dagger} \circ \xi_{M^*, A^{\dagger}} \circ (\id_{M^{**}} \otimes \tilde{\alpha} \eta) : G(M) \to M^{\dagger\dagger}
    \quad (M \in \mathcal{C}_A^{\sigma})
  \end{equation}
  of right $A$-modules in $\C$.
\end{lemma}
\begin{proof} 
  By Lemma~\ref{lem:actions-on-F(A**)}, the morphism $\tilde{\alpha} \circ u^{-1} : (A^{**})_{(u)} \to F(A^{**})$ is an isomorphism of $A^{**}$-$A$-modules. Thus we have a natural isomorphism
  \begin{align*}
    \phi_M' : (M^{**})_{(u)}
    \xrightarrow{\cong} M^{**} \otimes_{A^{**}} (A^{**})_{(u)}
    \xrightarrow{\cong} M^{**} \otimes_{A^{**}} F(A^{**})
    \xrightarrow{\cong} F(M^{**}) = M^{\dagger\dagger}
  \end{align*}
  of right $A$-modules. Specifically, the first arrow is the composition of $\id_{M^{**}} \otimes \eta^{**}$ and the quotient morphism $M^{**} \otimes A^{**} \to M^{**} \otimes_{A^{**}} A^{**}$, the second one is induced by $\tilde{\alpha} \circ u^{-1}$, and the third one is obtained by the Eilenberg-Watts theorem. We recall that the module structure $\xi'$ of $F$ is given by \eqref{eq:isomorphism-xi-twisted}. Thus, by the proof of the Eilenberg-Watts equivalence \cite[Theorem~4.2]{pareigis1977non}, we find that the third arrow is induced by
  \begin{equation*}
      M^{**} \otimes F(A^{**})
      \xrightarrow{\quad \xi_{M^{**}, A^{**}}' \quad}
      F(M^{**} \otimes A^{**})
      \xrightarrow{\quad F((a_{M}^r)^{**}) \quad} F(M^{**}).
  \end{equation*}
  Thus the isomorphism $\phi_M'$ is equal to the following composition:
\begin{equation*}
  \begin{gathered}
    M^{**} \xrightarrow{\quad \id \otimes \eta^{**} \quad}
    M^{**} \otimes A^{**} \xrightarrow{\quad \id \otimes u^{-1} \quad}
    M^{**} \otimes A \\ \xrightarrow{\quad \id \otimes \tilde{\alpha} \quad}
    M^{**} \otimes A^{\dagger\dagger}
    \xrightarrow{\quad \xi_{M^*,A^{\dagger}} \quad}
    (M \otimes A)^{\dagger\dagger}
    \xrightarrow{\quad (a_M^r)^{\dagger\dagger} \quad} M^{\dagger\dagger}.
  \end{gathered}
\end{equation*}
Since $u$ is a component of the Drinfeld isomorphism, we have $u^{-1} \circ \eta^{**} = \eta$.
Therefore $\phi_M'$ is equal to the morphism $\phi_M$ in the statement of this lemma. The proof is done.
\end{proof}


\subsection{Step 2: monoidal structure}\label{subsec:step-2}

The double dual functor of a rigid monoidal category has a canonical structure of a monoidal functor. In our setting, the monoidal structure of $(-)^{\dagger\dagger} : {}_A\mathcal{C}_A \to {}_A\mathcal{C}_A$ is given explicitly by
\begin{gather*}
  (\alpha_A^{-1})^{\dagger} \circ \alpha_A : A \to A^{\dagger\dagger}, \quad
  (\gamma_{N,M}^{-1})^{\dagger} \circ \gamma_{M^{\dagger}, N^{\dagger}}:
  M^{\dagger\dagger} \otimes_A N^{\dagger\dagger}
  \to (M \otimes_A N)^{\dagger\dagger}
\end{gather*}
for $M, N \in {}_A\mathcal{C}_A$, where $\alpha$ and $\gamma$ are natural transformations defined by \eqref{eq:internal-Hom-iso-alpha} and \eqref{eq:gamma-def}, respectively.
We describe how the evaluation morphism $\underline{\ev}_{M^{\dagger}, A} : M^{\dagger\dagger} \otimes_A M^{\dagger} \to A$ looks like if we identify $M^{\dagger\dagger}$ with $G(M) = M^{**}$ via the natural isomorphism $\phi_M$ given in Lemma~\ref{lem:CA-double-dual-as-functor}.
For $P \in {}_A\C_A$, we define
\begin{equation}
  \label{eq:thm:CA-double-dual:epsilon-def}
  \varepsilon_P := (\ev_{P^*} \otimes \id_A)
  \circ (\id_{P^{**}} \otimes \sigma_{P^*})
  \circ (\id_{P^{**}} \otimes i_{P,A}):
  P^{**} \otimes P^{\dagger} \to A.
\end{equation}
By the (di)naturality of the evaluation and the half-braiding, we have
\begin{equation}
  \label{eq:thm:CA-double-dual:epsilon-dinat}
  \varepsilon_P \circ (f^{**} \otimes \id_{P^{\dagger}})
  = \varepsilon_Q \circ (\id_{Q^{**}} \otimes f^{\dagger})
\end{equation}
for any morphism $f: P \to Q$ in ${}_A\mathcal{C}_A$.

\begin{lemma}
  \label{lem:CA-double-dual-eval}
  For $M \in \mathcal{C}_A^{\sigma}$, we have
  \begin{equation}
    \label{eq:CA-double-dual-eval}
      \underline{\ev}_{M^{\dagger},A} \circ (\phi_M \otimes_A \id_{M^{\dagger}}) \circ \pi_{M^{**}, M^{\dagger}} = \varepsilon_M.
  \end{equation}
\end{lemma}
\begin{proof}
  We consider the commutative diagram given as Figure~\ref{fig:lem:CA-double-dual-eval}, where $\delta_A$ is given by \eqref{eq:def-delta-A}.
  We note that the left column is equal to $(\phi_M \otimes_A \id_{M^{\dagger}}) \circ \pi_{M^{**}, M^{\dagger}}$. Thus we have
  \begin{align*}
    & \underline{\ev}_{M^{\dagger}, A} \circ (\phi_M \otimes_A \id_{M^{\dagger}}) \circ \pi_{M^{**}, M^{\dagger}} \\
    & = (\id_A \otimes \ev_{A \otimes M^*})
      \circ (\id_A \otimes \id_{M^{**}} \otimes \id_{A^*} \otimes i_{M,A}) \\
    & \qquad \qquad
      \circ (\sigma_{M^{**}}^{-1} \otimes \id_{A^*} \otimes \id_{M^{\dagger}})
      \circ (\id_{M^{**}} \otimes \coev_A \otimes \id_{M^{\dagger}}) \\
    & = (\id_A \otimes \ev_M) \circ (\sigma_{M^{**}}^{-1} \otimes \id_{M^*}) \circ (\id_{M^{**}} \otimes i_{M,A}) \\
    {}^{\eqref{eq:inverse-of-half-br}}
    & = (\ev_{M^*} \otimes \id_A)
  \circ (\id_{M^{**}} \otimes \sigma_{M^*})
  \circ (\id_{M^{**}} \otimes i_{M,A})
    = \varepsilon_M. \qedhere
  \end{align*}
\end{proof}

\begin{figure}
  \begin{equation*}
    \begin{tikzcd}[column sep=24pt, row sep=24pt]
      M^{**} \otimes M^{\dagger}
      \arrow[d, "{\id \otimes \eta \otimes \id}"']
      \arrow[rr, equal]
      & & M^{**} \otimes M^{\dagger}
      \arrow[d, "{\id \otimes \coev_A \otimes \id}"]
      \\ 
      M^{**} \otimes A \otimes M^{\dagger}
      \arrow[rr, "{\id \otimes \delta_A}"]
      \arrow[d, "{}"']
      \arrow[rrd, phantom, "{\text{(Figure \ref{fig:lem:actions-on-F(A**)-1})} \otimes M^{\dagger}}"]
      & & M^{**} \otimes A \otimes A^* \otimes M^{\dagger}
      \arrow[d, "{}"]
      \\[16pt] 
      {} M^{\dagger\dagger} \otimes M^{\dagger}
      \arrow[r, hook, "{i_{[M,A], A} \otimes \id}" {yshift = 5pt}]
      \arrow[d, "{\pi_{M^{\dagger\dagger}, M^{\dagger}}}"']
      & A \otimes (M^{\dagger})^* \otimes M^{\dagger}
      \arrow[d, "{\id \otimes \ev_{M^{\dagger}, A}}"]
      & A \otimes M^{**} \otimes A^* \otimes M^{\dagger}
      \arrow[l, two heads, "{\id \otimes i_{M, A}^* \otimes \id}"' {yshift = 5pt}]
      \arrow[d, "{\id \otimes \id \otimes \id \otimes i_{M, A}}"]
      \\ 
      M^{\dagger\dagger} \otimes_A M^{\dagger}
      \arrow[r, "{\underline{\ev}_{M^{\dagger}, A}}"]
      & A &
      A \otimes M^{**} \otimes A^* \otimes A \otimes M^*
      \arrow[l, "{\id \otimes \ev_{A \otimes M^*}}"']
    \end{tikzcd}
  \end{equation*}
  \caption{Proof of Lemma~\ref{lem:CA-double-dual-eval}}
  \label{fig:lem:CA-double-dual-eval}
\end{figure}


\subsection{Step 3: finishing the proof}\label{subsec:step-3}

Now we have all the ingredients needed to prove the main result.

\begin{proof}[Proof of Theorem~\ref{thm:CA-double-dual}]
We prove that $\phi_M : G(M) \to M^{\dagger\dagger}$ ($M \in \C_A^{\sigma}$) is an isomorphism of monoidal functors, that is, the equations
\begin{equation*}
\phi_A \circ \mathtt{g}_0
= (\alpha_A^{-1})^{\dagger} \circ \alpha_A
\quad \text{and} \quad
(\gamma_{N, M}^{-1})^{\dagger} \circ \gamma_{M^{\dagger}, N^{\dagger}} \circ (\phi_M \otimes_A \phi_N)
= \phi_{M \otimes_A N} \circ \mathtt{g}_{M,N}
\end{equation*}
hold for all $M, N \in \C_A^{\sigma}$.
Here, $\mathtt{g}_0$ and $\mathtt{g}_{M,N}$ ($M, N \in \C_A^{\sigma}$) are the monoidal structure of $G$ given by \eqref{eq:CA-double-dual-monoidal-structure-def-1}. Also, $\alpha_A$ and $\gamma_{M,N}$ are the monoidal structure of $(-)^{\dagger}$ given in \eqref{eq:internal-Hom-iso-alpha} and \eqref{eq:gamma-def}, respectively. The isomorphism $\phi_M$ is defined by \eqref{eq:CA-double-dual-iso-phi-def}.

  We first verify the equation $\phi_A \circ \mathtt{g}_0 = \tilde{\alpha}$, where $\tilde{\alpha} := (\alpha_A^{-1})^{\dagger} \circ \alpha_A$. We compute
  \begin{align*}
    \phi_A \circ \mathtt{g}_0 \circ \eta
    & = \mu^{\dagger\dagger} \circ \xi_{A^*, A^{\dagger}} \circ (\eta^{**} \otimes \tilde{\alpha} \eta) \\
    & = (\mu \circ (\eta \otimes \id_A))^{\dagger\dagger} \circ \xi_{{}^{*}\unit, A^{\dagger}} \circ (\id_{\unit} \otimes \tilde{\alpha} \eta)
      = \tilde{\alpha} \circ \eta,
  \end{align*}
  where the first equality follows from \eqref{eq:CA-double-dual-iso-phi-def}, the  naturality of $\xi$ is used at the second one, and the third one follows from the unit axiom.
  This implies $\phi_A \circ \mathtt{g}_0 = \tilde{\alpha}$, since both sides are morphisms of $A$-modules from $A$.

  To complete the proof, we fix objects $M, N \in \mathcal{C}_A^{\sigma}$ and verify the equation
  \begin{equation*}
    \gamma_{M^{\dagger}, N^{\dagger}} \circ (\phi_M \otimes_A \phi_N)
    = (\gamma_{N, M})^{\dagger} \circ \phi_{M \otimes_A N} \circ \mathtt{g}_{M,N}.
  \end{equation*}
  Let $L$ and $R$ be the left and the right hand side of this equation.
  Using notation in Lemmas \ref{lem:internal-Hom-iso-gamma-eq-1} and \ref{lem:CA-double-dual-eval}, we compute: \allowdisplaybreaks[4]
  \begin{align*}
    & \underline{\ev}_{N^{\dagger} \otimes_A M^{\dagger}, A}
      \circ (L \otimes_A \id_{N^{\dagger}} \otimes_A \id_{M^{\dagger}})
      \circ \pi_{M^{**},N^{**},N^{\dagger},M^{\dagger}}^{(4)} \\
    {}^{\eqref{eq:internal-Hom-iso-gamma-def}}
    & = \underline{\ev}_{M^{\dagger}, A}
      \circ (\id_{M^{\dagger\dagger}} \otimes_A \underline{\ev}_{N^{\dagger}, A} \otimes_A \id_{M^{\dagger}})
      \circ (\phi_M \otimes_A \phi_N \otimes_A \id_{N^{\dagger}} \otimes_A \id_{M^{\dagger}}) \\*
    & \qquad \circ (\id_{M^{**}} \otimes_A \pi_{N^{**},N^{\dagger}} \otimes_A \id_{M^{\dagger}})
    \circ \pi_{M^{**}, N^{**} \otimes N^{\dagger}, M^{\dagger}}^{(3)} \\*
    & \qquad \text{(where the isomorphism $A \otimes_A M^{\dagger} \cong M^{\dagger}$ is implicit)} \\
    {}^{\eqref{eq:CA-double-dual-eval}}
    & = \underline{\ev}_{M^{\dagger}, A} \circ (\phi_M \otimes_A \varepsilon_N \otimes_A \id_{M^{\dagger}})
    \circ \pi_{M^{**}, N^{**} \otimes N^{\dagger}, M^{\dagger}}^{(3)} \\    
    {}^{\text{(nat.)}}
    & = \underline{\ev}_{M^{\dagger}, A} \circ (\phi_M \otimes_A \id_{M^{\dagger}})
    \circ \pi_{M^{**}, A, M^{\dagger}}^{(3)}
      \circ (\id_{M^{**}} \otimes \varepsilon_N \otimes \id_{M^{\dagger}}) \\
    {}^{\eqref{eq:internal-Hom-iso-gamma-eq-1-proof-1}}
    & = \underline{\ev}_{M^{\dagger}, A} \circ (\phi_M \otimes_A \id_{M^{\dagger}})
    \circ \pi_{M^{**}, M^{\dagger}}^{(3)}
      \circ (\id_{M^{**}} \otimes a_{M^{\dagger}}^{l} (\varepsilon_N \otimes \id_{M^{\dagger}})) \\
    {}^{\eqref{eq:CA-double-dual-eval}}
    & = \varepsilon_M \circ (\id_{M^{**}} \otimes a_{M^{\dagger}}^{l} (\varepsilon_N \otimes \id_{M^{\dagger}})) \\
    {}^{\eqref{eq:thm:CA-double-dual:epsilon-def}}
    & = (\ev_{M^*} \otimes \id_A)
      \circ (\id_{M^{**}} \otimes \sigma_{M^*} i_{M,A} a_{M^{\dagger}}^{l} (\varepsilon_N \otimes \id_{M^{\dagger}})) \\
    {}^{\text{($*$)}}
    & = (\ev_{M^*} \otimes \id_A)
      \circ (\id_{M^{**}} \otimes \sigma_{M^*} (\mu \otimes \id_{M^*}) (\varepsilon_N \otimes i_{M,A})) \\
    {}^{\eqref{eq:thm:CA-double-dual:epsilon-def}}
    & = (\ev_{M^*} \otimes \id_A)
      \circ (\id_{M^{**}} \otimes \ev_{N^*} \otimes \sigma_{M^*}) \\*
    & \qquad \circ (\id_{M^{**} \otimes N^{**}} \otimes \mu \otimes \id_{M^*})
      \circ (\id_{M^{**} \otimes N^{**}} \otimes \sigma_{N^*} i_{N,A} \otimes i_{M,A}) \\
    & = (\ev_{(M \otimes N)^*} \otimes \id_A)
      \circ (\id_{M^{**}} \otimes \id_{N^{**}} \otimes L'(i_{N,A} \otimes i_{M,A})),
  \end{align*}
  where `(nat.)' means that the equality at that point follows from the naturality of $\pi^{(n)}_{X_1, \cdots, X_n}$,
  $(*)$ follows from that $i_{M,A} : M^{\dagger} \to A \otimes M^*$ is a morphism of left $A$-modules, and $L' : A \otimes N^* \otimes A \otimes M^* \to N^* \otimes M^* \otimes A$ is the morphism given by
  \begin{equation*}
    L' = (\id_{N^*} \otimes \sigma_{M^*})
    \circ (\id_{N^*} \otimes \mu \otimes \id_{M^*})
    \circ (\sigma_{N^*} \otimes \id_{A} \otimes \id_{M^*}).
  \end{equation*}

  We also compute $\underline{\ev}_{N^{\dagger} \otimes_A M^{\dagger}, A} \circ (R \otimes_A \id_{N^{\dagger}} \otimes_A \id_{M^{\dagger}}) \circ \pi_{M^{**},N^{**},N^{\dagger},M^{\dagger}}^{(4)}$.
  For this purpose, we note that the following equation holds:
  \begin{equation}
    \label{eq:thm:CA-double-dual:proof-1}
    (\id_A \otimes a_{N^*}^r) \circ (i_{N,A} \otimes \id_A)
    = (\mu \otimes \id_{N^*}) \circ (\id_A \otimes \sigma_{N^*}^{-1}) \circ (i_{N,A} \otimes \id_A).
  \end{equation}
  Indeed, we have $(\id_A \otimes a_{N^*}^r) (i_{N,A} \otimes \id_A)$
  \begin{align*}
    & = i_{N,A} \circ a_{N^{\dagger}}^r
    = i_{N,A} \circ a_{N^{\dagger}}^l \circ \sigma_{N^{\dagger}}^{-1} \\
    & = (\mu \otimes \id_{N^*}) \circ (\id_A \otimes i_{N,A}) \circ \sigma_{N^{\dagger}}^{-1} \\
    & = (\mu \otimes \id_{N^*}) \circ \sigma_{A \otimes N^*}^{-1} \circ (i_{N,A} \otimes \id_A) \\
    & = (\mu \otimes \id_{N^*}) \circ (\sigma_A^{-1} \otimes \id_{N^*})
    \circ (\id_A \otimes \sigma_{N^*}^{-1}) \circ (i_{N,A} \otimes \id_A) \\
    & = (\mu \otimes \id_{N^*})
    \circ (\id_A \otimes \sigma_{N^*}^{-1}) \circ (i_{N,A} \otimes \id_A),
  \end{align*}
  where the first and the third equalities follow from that $i_{N,A} : N^{\dagger} \to N \otimes A^*$ is a morphism of $A$-bimodules, the second equality holds since $N^{\dagger} \in \C_A^{\sigma}$, the fourth one follows from the naturality of the half-braiding, and the last one follows from the commutativity of $A$. Thus we compute
  \begin{align*}
    & \underline{\ev}_{N^{\dagger} \otimes_A M^{\dagger}, A}
      \circ (R \otimes_A \id_{N^{\dagger}} \otimes_A \id_{M^{\dagger}})
      \circ \pi_{M^{**}, N^{**}, N^{\dagger}, M^{\dagger}}^{(4)} \\
    {}^{\text{(nat.)}}
    & = \underline{\ev}_{(M \otimes_A N)^{\dagger}, A}
      \circ (\phi_{M \otimes_A N} \mathtt{g}_{M,N} \otimes_A \gamma_{N,M})
      \circ \pi_{M^{**}, N^{**}, N^{\dagger}, M^{\dagger}}^{(4)} \\
    {}^{\text{(nat.)}}
    & = \underline{\ev}_{(M \otimes_A N)^{\dagger}, A}
      \circ (\phi_{M \otimes_A N} \otimes \id_{(M \otimes_A N)^{\dagger}}) \\*
    & \qquad \circ \pi_{(M \otimes_A N)^{**}, (M \otimes_A N)^{\dagger}}
      \circ (\mathtt{g}_{M,N} \pi_{M^{**}, N^{**}} \otimes \gamma_{N,M} \pi_{N^{\dagger}, M^{\dagger}}) \\
    {}^{\eqref{eq:CA-double-dual-monoidal-structure-def-1}, \eqref{eq:CA-double-dual-eval}}
    & = \varepsilon_{M \otimes_A N}
      \circ (\pi_{M, N}^{**} \otimes \gamma_{N,M} \pi_{N^{\dagger}, M^{\dagger}}) \\
    {}^{\eqref{eq:thm:CA-double-dual:epsilon-dinat}}
    & = \varepsilon_{M \otimes N}
      \circ (\id_{M^{**}} \otimes \id_{N^{**}}
      \otimes \pi_{M, N}^{\dagger} \gamma_{N,M} \pi_{N^{\dagger}, M^{\dagger}}) \\
    {}^{\eqref{eq:thm:CA-double-dual:epsilon-def}}
    & = (\ev_{(M \otimes N)^*} \otimes \id_A)
      \circ (\id_{(M \otimes N)^{**}} \otimes \sigma_{(M \otimes N)^*} i_{M \otimes N, A} \pi_{M, N}^{\dagger} \gamma_{N,M} \pi_{N^{\dagger}, M^{\dagger}}) \\
    {}^{\eqref{eq:internal-Hom-iso-gamma-eq-1}}
    & = (\ev_{(M \otimes N)^*} \otimes \id_A)
      \circ (\id_{M^{**}} \otimes \id_{N^{**}} \otimes \sigma_{(M \otimes N)^*}) \\*
    & \qquad \circ (\id_{M^{**}} \otimes \id_{N^{**}} \otimes
      (\id_{A} \otimes a_{N^*}^r \otimes \id_{M^*}) (i_{N,A} \otimes i_{M,A})), \\
    {}^{\eqref{eq:thm:CA-double-dual:proof-1}}
    & = (\ev_{(M \otimes N)^*} \otimes \id_A)
      \circ (\id_{M^{**}} \otimes \id_{N^{**}} \otimes \sigma_{(M \otimes N)^*}) \\*
    & \qquad \circ (\id_{M^{**}} \otimes \id_{N^{**}} \otimes
      ((\mu \otimes \id_{N^*})(\id_A \otimes \sigma_{N^*}^{-1}) \otimes \id_{M^*}) (i_{N,A} \otimes i_{M,A})) \\
    & = (\ev_{(M \otimes N)^*} \otimes \id_A)
      \circ (\id_{M^{**}} \otimes \id_{N^{**}} \otimes R'(i_{N,A} \otimes i_{M,A})),
  \end{align*}
  where (nat.) follows from the (di)naturality of $\underline{\ev}_{X,A}$ and $\pi_{X,Y}$ in $X$ and $Y$, and $R' : A \otimes N^* \otimes A \otimes M^* \to N^* \otimes M^* \otimes A$ is the morphism given by
  \begin{equation*}
    R' = \sigma_{(M \otimes N)^*}
    \circ (\mu \otimes \id_{N^*} \otimes \id_{M^*})
    \circ (\id_{A} \otimes \sigma_{N^*}^{-1} \otimes \id_{M^*}).
  \end{equation*}
  By our assumption, $\mu$ is a morphism in $\mathcal{Z}(\mathcal{C})$.
  An easy graphical calculus shows $L' = R'$. Hence we have obtained the equation
  \begin{align*}
    & \underline{\ev}_{N^{\dagger} \otimes_A M^{\dagger}, A}
      \circ (L \otimes_A \id_{N^{\dagger}} \otimes_A \id_{M^{\dagger}})
      \circ \pi_{M^{**},N^{**},N^{\dagger},M^{\dagger}} \\
    = \mathrel{} \mbox{}
    & \underline{\ev}_{N^{\dagger} \otimes_A M^{\dagger}, A}
      \circ (R \otimes_A \id_{N^{\dagger}} \otimes_A \id_{M^{\dagger}})
      \circ \pi_{M^{**},N^{**},N^{\dagger},M^{\dagger}},
  \end{align*}
  which implies $L = R$. The proof is done.
\end{proof}


\section{(Braided) finite tensor categories}\label{sec:finite-cats}
This section focuses on braided FTCs $\mathcal{C}$ over an algebraically closed field $\kk$. In particular, we apply the results of Sections~\ref{sec:closed} and \ref{sec:double-dual-rigid} to an algebra $A$ in $\mathcal{C}$. Understanding the rigidity of ${}_A\mathcal{C}_A$ is the first step, and this occurs if and only if $A$ is exact \cite{etingof2004finite}. We recall some known results about exact algebras in Section~\ref{subsec:exact-algebras}.

For commutative, haploid and exact algebras $A$, we demonstrate that $\mathcal{C}_A^{\mathrm{loc}}$ is a braided FTC. This is detailed in Section~\ref{subsec:comm-exact-algebras}, where we prove this result and examine additional properties of commutative exact algebras. The core of this section is Section~\ref{subsec:MTC-main-result}, where we establish sufficient conditions for $\mathcal{C}_A^{\mathrm{loc}}$ to be pivotal, ribbon and modular.
Finally, in Section~\ref{subsec:relation-to-recent-work}, we compare our results with previous works \cite{kirillov2002q, davydov2013witt, laugwitz2023constructing}.


\subsection{Exact algebras}\label{subsec:exact-algebras}

Let $\C$ be a finite tensor category. An algebra $A\in\C$ is called \textit{exact} if for $P\in\C$ projective and any $M\in\C_A$, $P\otimes M\in\C_A$ is projective. 
We start by recalling various reformulations of the exactness property:

\begin{theorem}\label{thm:exact-characterization}
  Let $\C$ be a finite tensor category and $A$ an indecomposable algebra in $\C$. Then, the following are equivalent:
  \begin{enumerate}
    \item[\textup{(a)}] $A$ is exact.
    \item[\textup{(b)}] $-\otimes_A - : \C_A\times {}_A\C\rightarrow \C$ is a biexact functor.
    \item[\textup{(c)}] $A$ is simple and there exists $X\in\C$ with an injective map ${}^*A\hookrightarrow A\otimes X$ in ${}_A\C$.
    \item[\textup{(d)}] ${}_A\C_A$ is rigid.
  \end{enumerate}
\end{theorem}
\begin{proof}
The equivalence (a)$\iff$(b) follows from \cite[Example~7.9.8]{etingof2016tensor}. The equivalence (a)$\iff$(c) is \cite[Theorem~B.1]{etingof2021frobenius}. For (a) $\implies$(d), see \cite[\S3.3]{etingof2004finite}.

To complete the proof, we prove (d)$\implies$(b). We assume that ${}_A\C_A$ is rigid. Take $M\in\C_A$ and $N\in{}_A\C$. Note that $A\otimes M,\, N\otimes A \in {}_A\C_A$ using the regular action of $A$. Using the description of internal Hom of ${}_A\C_A$ provided in \eqref{eq:internal-Hom-rigid-case}, we have that
\begin{equation*}
  (A \otimes M) \otimes_A (N \otimes A)
  \cong A \otimes {}^*\uHom_A(M, N^*) \otimes A.
\end{equation*}
As ${}_A\C_A$ is rigid, the tensor product $\otimes_A$ in ${}_A\C_A$ is biexact.
As tensoring reflects exact sequences (see Remark~\ref{rem:tensor-reflect-exact}) it follows that the functor ${}^*\uHom_A(-,(-)^*) = (-\otimes_A -):\C_A\times {}_A\C\rightarrow\C$ is biexact, that is, (b) holds. The proof is done.
\end{proof}

\begin{remark}\label{rem:exact-Deligne}
  Given two FTCs $\C$ and $\D$, and algebras $A\in\C$ and $B\in\D$, we have an equivalence ${}_{A\boxtimes B}(\C\boxtimes\D)_{A\boxtimes B} \cong {}_A\C_A\boxtimes{}_B\D_B$ of monoidal categories. Thus, if $A$ and $B$ are exact, then $A\boxtimes B$ is exact in $\C\boxtimes\D$. See \cite[Theorem~2.3.7]{douglas2018dualizable} for a proof.
\end{remark}

An algebra $(A,\mu,\eta)\in\C$ is called \textit{separable} if there exists an $A$-bimodule map $\Delta:A\rightarrow A\otimes A$ that satisfies $\Delta\circ \mu=\id_A$.

\begin{lemma}
  \label{lem:separable-is-exact}
  Separable algebras are exact.
\end{lemma} 
\begin{proof}
Let $\C$ be a finite tensor category, and let $A$ be an algebra in $\C$ with multiplication $\mu$ and unit $\eta$.
By definition, there is a morphism $\tilde{e} : A \to A \otimes A$ of $A$-bimodules in $\C$ such that $\mu \circ \tilde{e} = \id_A$. We define
\begin{equation*}
    e_{X,Y} = (a_X^r \otimes a_Y^l) \circ (\id_X \otimes \tilde{e} \eta \otimes \id_Y) : X \otimes Y \to X \otimes Y
\end{equation*}
for $X \in \C_A$ and $Y \in {}_A\C$. The morphism $e_{X,Y}$ satisfies
\begin{equation*}
    e_{X,Y}^2 = e_{X,Y}, \quad
    \pi_{X,Y} \circ e_{X,Y} = \pi_{X,Y} \quad \text{and} \quad
    e_{X,Y} \circ (a_X^r \otimes \id_Y) = e_{X,Y} \circ (\id_X \otimes a_Y^l).
\end{equation*}
By these equations, we see that $X \otimes_A Y$ is isomorphic to the image of $e_{X,Y}$, which is a direct summand of $X \otimes Y$. Hence, for all $X \in \C$, the functor $X \otimes_A (-): {}_A \C \to \C$ is exact as a direct summand of the exact functor $X \otimes (-)$. In a similar way, $(-) \otimes_A Y$ is exact for all $Y \in {}_A\C$. Therefore $A$ is an exact algebra. The proof is done.
\end{proof}


\subsection{Commutative exact algebras}
\label{subsec:comm-exact-algebras}

\begin{theorem}\label{thm:C-A-rigid}
Let $\C$ be an FTC, and let $A$ be an algebra in $\C$.
\begin{enumerate}
  \item[\textup{(a)}] If $A$ is a commutative central algebra with half-braiding $\sigma$ then, $A$ is exact if and only if $\C_A^{\sigma}$ is rigid.
  \item[\textup{(b)}] If $\C$ is braided and $A$ is commutative in $\C$, then the monoidal category $\C_A$ (or ${}_A\C$) is rigid if and only if $A$ is exact. 
\end{enumerate}
\end{theorem}
\begin{proof}
If $A$ is exact, then ${}_A\C_A$ is rigid and hence $\C_A^{\sigma}$ is rigid by Theorem~\ref{thm:ACA-rigid-consequences}. Conversely, if $\C_A^{\sigma}$ is rigid, then the tensor product $\otimes_A : \C_A \times {}_A \C \to \C$ is exact in each variable. Hence $A$ is exact and part (a) is proved.

Part (b) for $\C_A$ is proved by applying part (a) to $A$ with half-braiding $\sigma = c_{A,-}$. As $(-)^*:{}_A\C\rightarrow\C_A$ is a tensor equivalence, the result for ${}_A\C$ also follows.
\end{proof}

By the same argument as the above lemma, we get the following result.
\begin{theorem}\label{thm:C-A-loc-rigid}
Let $\C$ be a braided finite tensor category and $A$ a commutative exact algebra in $\C$. Then the categories $\C_A,\; {}_A\C$ and $\C_A^{\loc}$ are finite multi-tensor categories. If $A$ is haploid, then these are finite tensor categories. \qed
\end{theorem}

The following is a basic way of obtaining new commutative exact algebras from known ones.
\begin{lemma}\label{lem:exact-preserve}
  Let $A$ be a haploid, commutative exact algebra in a braided FTC $\C$. If $F:\C\rightarrow\D$ is a fully faithful braided tensor functor, then $F(A)$ is a haploid, commutative and exact.
\end{lemma}
\begin{proof}
  By \cite[Corollary~B.5]{etingof2021frobenius}, $F(A)$ is exact. As $F$ is braided, $F(A)$ is commutative. Also, $F$ being fully faithful implies that $F(A)$ is haploid.
\end{proof}

\begin{remark}
    In a braided monoidal category, the tensor product of two algebras $A$ and $B$ is also an algebra. If in addition, $c_{A,B}\circ c_{B,A} = \id_{A\otimes B}$, then $A\otimes B$ is commutative. However, the tensor product of exact algebras may not be exact \cite[Example~7.6]{coulembier2025simple}. This is unlike separable algebras where the tensor product of commutative separable algebras is commutative separable.
\end{remark}


\subsubsection{Haploid vs indecomposable}
When working with exact algebras, we consider indecomposable algebras. For commutative exact algebras, we focus on haploid algebras because the haploid property is easier to check in examples. Additionally, for commutative exact algebras, haploid and indecomposable are equivalent notions, as shown in the following lemma.
\begin{lemma}\label{lem:comm-indecomposable-haploid}
  A commutative exact algebra is haploid if and only if it is indecomposable.
\end{lemma}
\begin{proof}
  The `only if' part is obvious. We prove the converse.
  Suppose that a commutative exact algebra $A$ in a braided FTC $\C$ is indecomposable. As $A$ is exact, $\C_A$ is rigid monoidal. In fact, it is a multi-tensor category. Now, by the same argument as in \cite[Corollary~4.22]{laugwitz2023constructing}, it follows that $A$ is haploid.
\end{proof}

\begin{lemma}\label{lem:comm-simple}
A commutative simple algebra is haploid and indecomposable.
\end{lemma}
\begin{proof}
One can check that $\Hom_{\C}(\unit,A)\cong \Hom_{\C_A}(A,A) \cong \Hom_{{}_A\C_A}(A,A) \cong\kk$. Here the second last isomorphism holds because $A$ is commutative and the last isomorphism because $A$ is simple. Thus, $A$ is haploid. As haploid algebras are indecomposable, the proof is done.
\end{proof}


\subsubsection{The free functor}\label{subsubsec:FA'}
Fix $\C$ to be a braided FTC and $A$ a commutative, exact and haploid algebra in it. Then we have the free functor 
\[F_A:\C\rightarrow\C_A, \qquad X \mapsto (X\otimes A, \id_X\otimes m). \]
The functor $F_A$ is a surjective tensor functor and its right adjoint is the forgetful functor $F_A^{\ra}=U_A:\C_A\rightarrow\C$ given by $U_A(M,r_M) = M$ \cite[Lemma~4.5]{laugwitz2023constructing}. Furthermore, $F_A$ is a central functor, namely we have a faithful tensor functor $F_A':\C\rightarrow\Z(\C_A)$ such that $F_A = U\circ F_A'$ \cite[Proposition~8.8.10]{etingof2016tensor}.

\begin{lemma}
  \label{lem:free-functor-braided}
  Consider the above setting. Then, $A\in\C'$ if and only if there is a braiding on $\C_A$ making the free functor $F_A$ braided.
\end{lemma}
\begin{proof}
($\Rightarrow$) Suppose that $\C_A$ has a braiding $c^A$ such that $F_A$ is braided. We will show below that every right $A$-module is local. Then in particular, the free module $X \otimes A$ is local for all $X \in \C$. This means that $A$ belongs to the M\"uger center of $\C$ \cite[Lemma~3.15]{davydov2013witt}. 

Since the forgetful functor $U_A: \C_A \to \C$ is right adjoint to
$F_A$, it is a lax braided monoidal functor \cite{kelly1974doctrinal}. Thus the braiding of $\C_A$ satisfies the equation $c^A_{M, N} \circ \pi_{M, N} = \pi_{N, M} \circ c_{M, N}$ for all $M, N \in \C_A$. Since the left hand side is $A$-balanced, so is the right hand side. Hence we have
\begin{align*}
    \pi_{N, M} \circ c_{M, N} \circ (a_M^{r} \otimes \id_N)
    & = \pi_{N, M} \circ c_{M, N} \circ (\id_M \otimes a_N^l) \\
    & = \pi_{N, M} \circ c_{M, N} \circ (\id_M \otimes a_N^r)
    \circ (\id_M \otimes c_{A, N}) \\
    {}^{\text{(nat.)}}
    & = \pi_{N, M} \circ (a_N^r \otimes \id_M) \circ c_{M, N \otimes A} \circ (\id_M \otimes c_{A,N}) \\
    {}^{\text{(bal.)}}
    & = \pi_{N, M} \circ (\id_N \otimes a_M^l) \circ c_{M, N \otimes A} \circ (\id_M \otimes c_{A,N}) \\
    & = \pi_{N, M} \circ (\id_N \otimes a_M^r c_{A,M} c_{M,A}) \circ (c_{M, N} \otimes \id_{A}) \circ (\id_M \otimes c_{A,N}) \\
    {}^{\text{(nat.)}}
    & = \pi_{N, M} \circ c_{M,N} \circ (a_M^r c_{A,M} c_{M,A} \otimes \id_N),
\end{align*}
where (nat.) and (bal.) follows from the naturality of the braiding and the balancing property of $\pi$, respectively.
By letting $A = N$ and composing $\id_M \otimes \eta$ to both sides, we obtain $a_M^r = a_M^r c_{A,M} c_{M,A}$. Namely, $M$ is local.

($\Leftarrow$) Suppose that $A$ is in the M\"uger center of $\C$. Then $\C_A\cong \C_A^{\loc}$, so it is braided. For a proof that $F_A$ is braided, see for example \cite[Theorem~2.67]{creutzig2017tensor}.
\end{proof}


\subsubsection{Further properties}

The following lemma is a generalization of \cite[Corollary~3.26]{davydov2013witt}.
\begin{lemma}\label{lem:nondg-fully-faithful}
  Let $\C,\D$ be braided FTCs with $\C$ non-degenerate. Let $F:\C\rightarrow\D$ be a braided tensor functor. Then, $F$ is full, hence fully faithful.
\end{lemma}
\begin{proof}
By replacing $\D$ with the image of $F$, we may assume that $F$ is surjective.
We consider the commutative algebra $A = F^{\ra}(\unit)$ in $\mathcal{C}$.
There is an equivalence $\D \simeq \C_A$ induced by $F^{\ra}$ and, if we identify $\D$ with $\C_A$, then $F$ is identified with the free module functor $F_A : \C \to \C_A$ given by $F(X) = X \otimes A$ for $X \in \C$ \cite[Proposition~6.1]{bruguieres2011exact}.
By Lemma~\ref{lem:free-functor-braided}, $A$ belongs to the M\"uger center of $\C$. Since $\C$ is assumed to be non-degenerate, and since $\Hom_{\C}(\unit, A) \cong \Hom_{\D}(\unit, \unit) \cong \kk$, the algebra $A$ is actually isomorphic to $\unit$. Thus $F$ is, in fact, the identity functor.
The proof is done.
\end{proof}

\begin{remark}
  A generalization of the above result, where $\C$ is not necessarily non-degenerate, has appeared in \cite[Proposition~3.57]{decoppet2024verlinde}.
\end{remark}

Suppose that $\B$ and $\C$ are finite tensor categories with $\B$ braided. Let $G:\B\rightarrow\Z(\C)$ be a faithful braided tensor functor and set $\Z_{\B}(\C) = \Z_{(2)}(G(\B)\subset \Z(\C))$. In \cite{laugwitz2022relative}, the notation $\Z_{\B}(\C)$ is used to denote $\Z_{(2)}(\B\subset \Z(\C))$, which is different from our convention. 

Given a commutative algebra $(A,\sigma)\in\Z_{\B}(\C)$, \cite[Proposition~5.13]{laugwitz2023constructing} constructs a faithful braided tensor functor $G_A:\B\rightarrow\Z(\C_A^{\sigma})$. Using this, \cite[Theorem~5.15]{laugwitz2023constructing} proves the following equivalence of braided monoidal categories:
\begin{equation}\label{eq:LW-local}
  \Z_{\B}(\C)_{(A,\sigma)}^{\loc} \simeq \Z_{\B}(\C_A).
\end{equation}

One special case of this is when $\B=\Vect$. In this case, for a central algebra $(A,\sigma)\in\Z(\C)$, we have the following equivalence of braided monoidal categories \cite{schauenburg2001monoidal}:
\begin{equation}\label{eq:Schauenburg}
  \Z(\C)_{(A,\sigma)}^{\loc} \simeq \Z(\C_A^{\sigma}).
\end{equation}

The other case is when $\B=\C$. For this, recall the functor $F_A'$ from \S\ref{subsubsec:FA'}.

\begin{theorem}\label{thm:LW-local}
Let $\C$ be a finite tensor category and $A$ a commutative algebra in $\C$. Then $\C_A^{\loc} \simeq \Z_{(2)}(F_A'(\C)\subset \Z(\C_A))$.
\end{theorem}
\begin{proof}
This follows from \cite[Theorem~5.15]{laugwitz2023constructing} by taking $\B=\C$ and the functor $G=i_+:\C\rightarrow \Z(\C)$. In this setting, $\Z_{\B}(\C) \simeq \C$ and the functor $G_A$ is equal to $F_A': \C\rightarrow \Z(\C_A)$ . So that claim follows from (\ref{eq:LW-local}).
\end{proof}

\subsubsection{Non-degeneracy of \texorpdfstring{$\C_A^{\loc}$}{CA-loc}}
Throughout this subsection, $\C$ is a braided FTC and $A$ is a haploid, commutative and exact algebra in $\C$.

A particularly important case of Theorem \ref{thm:LW-local} is when $\C$ is non-degenerate. In this case, the functor $F_A'$ is fully faithful by Lemma~\ref{lem:nondg-fully-faithful}. So $\mathrm{Im}(F'_A) \simeq \C$ and we get 
\begin{equation*}
  \C_A^{\loc} \simeq \Z_{(2)}(\C\subset \Z(\C_A)).
\end{equation*}

\begin{corollary}\label{cor:ndBFTC-Z(C-A)}
If $\C$ is non-degenerate, then:
\begin{enumerate}
  \item $\Z(\C_A) \simeq \C_A^{\loc} \boxtimes \, \oC$ as braided tensor categories.
  \item $\C_A^{\loc}$ is a non-degenerate braided FTC. 
\end{enumerate} 
\end{corollary}
\begin{proof}
  As $\C_A^{\loc}$ is the centralizer of a nondegenerate category, part (a) follows from \cite[Theorem~4.17]{laugwitz2022relative}. Part (b) follows from the same argument as in \cite[Proposition~4.24]{laugwitz2023constructing}.
\end{proof}

\subsubsection{Frobenius-Perron dimensions}
Lastly, we record some consequences for FP dimensions. We use the same notations as in \cite{etingof2016tensor}.

\begin{lemma}\label{lem:FPdim-C_A}
Let $A$ be a haploid, commutative and exact algebra in $\C$. Then $\FPdim(\C_A)=\frac{\FPdim(\C)}{\FPdim_{\C}(A)}$.
\end{lemma}
\begin{proof}
By applying \cite[Lemma~6.2.4]{etingof2016tensor} to the functor $F_A$, it follows that 
\[\FPdim(\C_A) = \FPdim(\C)\frac{\FPdim_{\C_A}(A)}{\FPdim_{\C}(F_A^{\ra}(A))} = \frac{\FPdim(\C)}{\FPdim_{\C}(A)}.  \qedhere \]
\end{proof}

\begin{lemma}
  Let $\C$ be an integral braided FTC and $A$ a haploid, commutative and exact algebra in $\C$.
  Then $\C_A^{\loc}$ is an integral braided FTC and $\C_A$ is an integral FTC.
\end{lemma}

Hence, by \cite[Corollary 6.1.15]{etingof2016tensor}, there is a finite-dimensional quasi-Hopf algebra $Q$ and its quotient $\overline{Q}$ such that $\C_A \simeq \Rep(Q)$ and $\C_A^{\loc} \simeq \Rep(\overline{Q})$ as tensor categories.

\begin{proof}
  As $F_A:\C\rightarrow\C_A$ is a surjective tensor functor, the claim about $\C_A$ follows from \cite[Corollary~6.2.5]{etingof2016tensor}. Thus $\C_A$ has a quasi-fiber functor.
  As $\C_A^{\loc}$ is a tensor full subcategory of $\C_A$, by composing this inclusion with the quasi-fiber functor of $\C_A$, we get a quasi-fiber functor from $\C_A^{\loc}$. Thus, it is also integral.
\end{proof}

\begin{lemma}\label{lem:FPdim-C-A-loc}
Let $A$ be a haploid, commutative and exact algebra in a braided FTC $\C$. Then,
\begin{equation*}
\FPdim(\C_A^{\loc}) 
= \frac{\FPdim(\C_A)^2}{\FPdim(\mathrm{Im}(F_A'))} 
= \frac{\FPdim(\C)}{\FPdim_{\C}(A)^2} \cdot \frac{\FPdim(\C)}{\FPdim(\mathrm{Im}(F_A'))}
\end{equation*}
When $\C$ is nondegenerate, $\FPdim(\C_A^{\loc}) = \frac{\FPdim(\C)}{\FPdim(A)^2}$ and $\FPdim_{\C}(A)\leq \sqrt{\FPdim(\C)}$.
\end{lemma}
\begin{proof}
By Theorem~\ref{thm:LW-local}, we have that $\C_A^{\loc}\simeq \Z_{(2)}(\mathrm{Im}(F_A') \subset\Z(\C_A))$. Thus,
the claim follows by the formula for Frobenius-Perron dimension of M\"uger centralizers \cite[Theorem~4.9]{shimizu2019non}. When $\C$ is nondegenerate, $F_A'$ is an equivalence onto its image. This implies the second claim.
\end{proof}

If a haploid, commutative exact algebra $A$ satisfies $\FPdim(A)^2=\FPdim(\C)$, we call it \textit{Lagrangian}.

\begin{corollary}\label{cor:lagrangian}
Let $A$ be a haploid, commutative, exact algebra in a non-degenerate braided FTC $\C$. Then $A$ is a Lagrangian algebra if and only if $\C_A^{\loc} \simeq \Vect$.
\end{corollary}
\begin{proof}
If $A$ is Lagrangian, $\FPdim(\C_A^{\loc}) = 1$. Thus, $\C_A^\loc \simeq \Vect$. 
The converse is clear.
\end{proof}


\subsection{Main result}\label{subsec:MTC-main-result}
In this section, we will focus on ribbon braided FTCs and modular tensor categories. Recall that the base field $\Bbbk$ is algebraically closed.

\begin{proposition}
  \label{prop:local-module-sFrob-MTC}
  Suppose that $A$ is a haploid commutative exact algebra in an MTC $\C$.
  Then the following assertions are equivalent:
  \begin{enumerate}
  \item[\textup{(a)}] $\C_A^{\loc}$ is a ribbon category with the same twist as $\C$.
  \item[\textup{(b)}] $\theta_A = \id_A$ and $A^*\in(\C_A^{\loc})'$. 
  \item[\textup{(c)}] $A$ is a symmetric Frobenius algebra.
  \end{enumerate}
  If these equivalent conditions are satisfied, then $\C_A^{\loc}$ is a modular tensor category.
\end{proposition}
\begin{proof}
  Suppose that (a) holds. Then $\theta_A = \id_A$ since the twist must be the identity on the unit object. By Theorem \ref{thm:local-module-ribbon}, the object $A^*$ belongs to the M\"uger center of $\C_A^{\loc}$. Thus, (b) holds. 

  Now suppose that (b) holds. Then, by the non-degeneracy of $\C_A^{\loc}$ (Corollary~\ref{cor:ndBFTC-Z(C-A)}(b)), $A^*$ must be isomorphic to the direct sum of finitely many copies of $A$. By comparing the FP dimensions, we conclude that $A \cong A^*$. Namely, $A$ is a Frobenius algebra. Lastly by Lemma~\ref{lem:symmetric-frobenius-theta}, we conclude that $A$ is symmetric. Thus, (c) is proved.

  Suppose that $A$ is a symmetric Frobenius algebra. Then, by Lemma~\ref{lem:symmetric-frobenius-theta}, we have $\theta_A = \id_A$.
  Moreover, the object $A^*$ is isomorphic to $A$. Thus, by Theorem \ref{thm:local-module-ribbon}, $\C_A^{\loc}$ is a ribbon category with the twist $\theta$. This proves (c)$\Rightarrow$(a) and the proof is done.
\end{proof}

\begin{remark}
  By \cite{shimizu2024selfinjective}, condition (c) is proven to be equivalent to  $\C_A$ being pivotal as a left $\C$-module category. It would be interesting to examine if such a condition is related to the unimodularity and/or the sphericality of the tensor category $(\C_A, \otimes_A)$.
\end{remark}

We are now ready to state our main result.

\begin{theorem}\label{thm:FTC-main-result}
  Let $\C$ be a braided finite tensor category and $A$ a haploid, commutative and exact algebra in $\C$. Then $\C_A$ and $\C_A^{\loc}$ are finite tensor categories.
  \begin{enumerate}
    \item[\textup{(a)}] If $\C$ is pivotal and $\fu_A=\fp_A$, then $\C_A$ and $\C_A^{\loc}$ are pivotal FTCs.
    \item[\textup{(b)}] If $\C$ is ribbon, $\theta_A=\id_A$ and $A^*\in(\C_A^{\loc})'$, then $\C_A^{\loc}$ is a ribbon FTC. 
    \item[\textup{(c)}] If $\C$ is ribbon, $\theta_A=\id_A$ and $A\in\C'$, then $\C_A$ is a ribbon FTC.
    \item[\textup{(d)}] If $\C$ is an MTC and $A$ is symmetric Frobenius, then $\C_A^{\loc}$ is an MTC.
  \end{enumerate} 
\end{theorem}
\begin{proof}
For (a), recall that $\bA=(A,c_{A,-})$ is central algebra and $\C_A^{c_{A,-}} = \C_A$. In addition, $\fu_{\bA} = \fu_A$. Thus, by Theorem~\ref{thm:pivotality-C-A-sigma}, $\fu_A=\fp_A$ implies that $\C_A$ is a pivotal finite tensor category. By Theorem~\ref{thm:pivotality-C-A-loc}, $\C_A^{\loc}$ is also a pivotal FTC. This proves (a).
 
As $\C_A^{\loc}$ is rigid, part (b) follows from Theorem~\ref{thm:local-module-ribbon}. 
Part (c) is a direct consequence of Corollary~\ref{cor:local-module-ribbon}.
Lastly, part (d) is a consequence of Proposition~\ref{prop:local-module-sFrob-MTC}.
\end{proof}

\begin{remark}\label{rem:CSZ-simple}
  In \cite[Conjecture~B.6]{etingof2021frobenius}, it was conjectured that an indecomposable algebra $A$ is exact if and only if it is simple. 
  This conjecture was proven in \cite{coulembier2025simple}. Consequently, Theorem~\ref{thm:FTC-main-result} can be rephrased (using Lemma~\ref{lem:comm-simple}) by replacing `haploid and exact' in the assumptions by `simple'.   
\end{remark}
We also obtain the following corollary.

\begin{corollary}
Let $\C$ be an MTC and $A$ a simple commutative symmetric Frobenius algebra in $\C$. Then $\C_A^{\loc}$ is an MTC.
\end{corollary}
\begin{proof}
By Lemma~\ref{lem:comm-simple}, $A$ is indecomposable. Because $A$ is Frobenius, $A\cong A^*$ in $\C_A$. Thus, by Theorem~\ref{thm:exact-characterization}, $A$ is exact. Now the claim follows using Theorem~\ref{thm:FTC-main-result}(d).
\end{proof}


\subsubsection{Remarks about the conditions in Theorem~\ref{thm:FTC-main-result}}
In this section we discuss examples which show that the conditions used in Theorem~\ref{thm:FTC-main-result} are independent of each other.

\underline{$\theta_A \neq \id_A$ and $A^*\notin \C_A^{\loc}$}: Let $\C$ be a finite tensor category, and consider the commutative algebra $A = R(\unit)$, where $R$ is a right adjoint of the forgetful functor $\Z(\C) \to \C$.
It is known that $A^* \cong R(D)$, where $D$ is the distinguished invertible object of $\mathcal{B}$ (see \cite{shimizu2016unimodular}). Thus $A \cong A^*$ as objects of $\mathcal{Z}(\C)$ if and only if $D \cong \unit$. Indeed, the `if' part is obvious. The converse follows from
\begin{equation*}
  \Hom_{\mathcal{Z}(\mathcal{B})}(\unit, A) \cong \Hom_{\mathcal{B}}(\unit, \unit) \cong \Bbbk,
  \quad \Hom_{\mathcal{Z}(\mathcal{B})}(\unit, A^*) \cong \Hom_{\mathcal{B}}(\unit, D).
\end{equation*}
Assume $\mathcal{B}$ is not unimodular and $\C = \mathcal{Z}(\mathcal{B})$ has a ribbon structure (such an example is given by the category of finite-dimensional modules over the Taft algebra at a root of unity of odd order).
By Lemma~\ref{lem:FPdim-C-A-loc}, since $\FPdim(\C_A^{\loc}) = 1$, it is equivalent to $\Vect$. If $A^* \in \C_A^{\loc}$, then $A^*$ is isomorphic to a finite direct sum of $A$, but this is a contradiction. Thus $A^*$ cannot have a structure of a local module over $A$. Therefore $\theta_A \ne \id_A$. This example also shows that $\theta_A = \id_A$ is not a necessary condition for $\C_A^{\loc}$ to be ribbon. In this example, $\C_A^{\loc}$ is equivalent to $\Vect$, which is a ribbon category.

\underline{$\theta_A\neq\id$ and $A^*\in \C_A^{\loc}$}: 
From Proposition~\ref{prop:K-in-loc}, recall that $A^*\in \C_A^{\loc}$ if and only if $\theta_A^2=\id_A$. Here we show an example where $\theta_A\neq\id_A$ and $A^*\in \C_A^{\loc}$. To obtain examples, one can take any symmetric finite tensor category $\C$ with non-trivial twist and any commutative algebra $A$ in $\C$ with $\theta_A\neq\id$. 
For example, let $G$ be a finite abelian group, let $t : G \to \kk^\times$ be a character of order two, and let $H$ be a subgroup of $G$ such that $t|_H$ is non-trivial (e.g., $H = G$). Then the category $\Vect_G$ of finite-dimensional $G$-graded vector spaces is a symmetric finite tensor category with trivial associator and braiding but with non-trivial ribbon structure given by $t$. The group algebra $A = \kk H$ is a commutative algebra in $\Vect_G$ such that $\theta_A \ne \id_A$.

\begin{remark}
  If one knows an example of a unimodular Hopf algebra $H$ that does not admit a pivotal structure but its double $D(H)$ is spherical/ribbon, then a more interesting example can be constructed. In this setting, $\C$ is unimodular but not pivotal. On the other hand, $\Z(\C)$ is spherical (and ribbon). Thus, taking $A=R(\unit)$ in $\Z(\C)$, we get that $A$ will be Frobenius but it cannot be symmetric Frobenius (if so then $\Z(\C)_A \cong \C$ will be pivotal). Hence, $\theta_A\neq \id$. But $A^*\in\C_A^{\loc}$ because $A^*\cong A$ as $A$-modules. This would provide an example of a commutative exact Frobenius algebra with $\theta_A \neq \id$ such that $A^* \in \C_A^{\loc}$.
\end{remark}

\underline{$\theta_A\neq\id_A$ and $A^*\in (\C_A^{\loc})'$}: Such an example is provided in Subsection~\ref{subsec:ex1-even-order}.

\underline{$\theta_A=\id_A$ and $A^*\notin (\C_A^{\loc})'$}: Such an example is provided in Subsection~\ref{subsec:uqsl2}. 

These two examples show that the conditions of $\theta_A=\id_A$ and $A^*\in (\C_A^{\loc})'$ are independent of each other.


\subsubsection{Functors preserving haploid, commutative, exact algebras}
In light of our main result, Theorem~\ref{thm:FTC-main-result}, we need examples of haploid, commutative, exact algebras in braided FTCs with additional properties. Next we introduce certain properties of tensor functors that allow us to transport such algebras from one category to another.

Let $\C$ and $\D$ be finite tensor categories and $F:\C\rightarrow\D$ a tensor functor.
\begin{itemize}
  \item If $\C$ and $\D$ are pivotal, then for all $X\in\C$, there is natural isomorphism $\zeta_X:F(X^*)\rightarrow F(X)^*$. We call $F$ \textit{pivotal} \cite{ng2007higher} if it satisfies the following for all $X\in\C$:
  \[ \fp^{\D}_{F(X)} = ((\zeta_X)^{-1})^* \circ \zeta_{X^*} \circ F(\fp^{\C}_X) . \]
  \item If $\C$ and $\D$ are ribbon and $F$ is braided, we call $F$ \textit{ribbon} if it satisfies $F(\theta_X) = \theta_{F(X)}$ for all $X\in\C$. 
\end{itemize}
Note that a braided tensor functor is pivotal if and only if it is ribbon \cite[Proposition~4.4]{mulevicius2022condensation}.

\begin{lemma}\label{lem:F-preserves}
Let $F:\C\rightarrow\D$ be a fully faithful braided tensor functor between braided FTCs and $A$ a haploid, commutative and exact algebra in $\C$. 
\begin{enumerate}
  \item[\textup{(a)}] If $\C$ and $\D$ are pivotal, $F$ is pivotal and $\fu^{\C}_A=\fp^{\C}_A$, then $F(A)$ satisfies $\fu^{\D}_{F(A)}=\fp^{\D}_{F(A)}$.
  
  \item[\textup{(b)}] If $\C$ and $\D$ are ribbon, $F$ is ribbon and $\theta^{\C}_A=\id^{\C}_A$, then $F(A)$ satisfies $\theta^{\D}_{F(A)}=\id^{\D}_{F(A)}$.
  
  \item[\textup{(c)}] If $\C$ and $\D$ are MTCs, $F$ is ribbon/pivotal and $A$ is symmetric Frobenius, then so is $F(A)$.
\end{enumerate}
\end{lemma}
\begin{proof}
By Lemma~\ref{lem:exact-preserve}, $F(A)$ is exact, haploid and commutative. 

(a) Both $\C$ and $\D$ are braided and pivotal, they both admit a twist $\theta$ which satisfies $\fp=u\theta$. Thus, $u^{\C}_A=\fp^{\C}_A$ implies that $\theta_A^{\C}=\id_A^{\C}$. 
Using \cite[Proposition~4.4]{mulevicius2022condensation}, we have that $F$ is ribbon. Thus, $\theta_{F(A)}^{\D}=\id_{F(A)}$. Using the relation $\fp=\fu\theta$ for the twist $\theta$ of $\D$, the claim follows.

(b) Straightforward.

(c) This follows from the general fact that pivotal functors preserve symmetric Frobenius algebras \cite[Proposition~4.8]{mulevicius2022condensation}.
\end{proof}


\subsection{Relation to prior work}\label{subsec:relation-to-recent-work}
In \cite{kirillov2002q,davydov2013witt}, one gets that if $\C$ is a semisimple MTC and $A$ is a haploid, commutative and separable algebra with trivial twist, then $\C_A^{\loc}$ is a semisimple MTC. 
In \cite[Theorem~4.16]{laugwitz2023constructing}, this result was generalized to get rid of the semisimplicity assumption on $\C$. They used rigid Frobenius algebras $A$ in an MTC $\C$ to construct new MTCs $\C_A^{\loc}$. A rigid Frobenius algebra is a haploid, commutative special Frobenius algebra. 
Special Frobenius algebras are separable and hence exact by Lemma~\ref{lem:separable-is-exact}. Additionally, rigid Frobenius algebras satisfy $\theta_A=\id_A$ which in turn implies that they are symmetric. Thus, a rigid Frobenius algebra satisfies the conditions of Theorem~\ref{thm:FTC-main-result}(d).

In Section~\ref{subsec:comm-alg-in-ZC} we provide a class of examples to illustrate that Theorem~\ref{thm:FTC-main-result}(d) generalizes the result of \cite{laugwitz2023constructing} by removing the assumption that $A$ is special Frobenius and by replacing it by the weaker assumption that $A$ is symmetric Frobenius and exact. 

A different generalization of \cite{kirillov2002q}, namely to monoidal supercategories, was discussed in \cite{creutzig2017tensor}. We don't discuss this in this article, but this aspect of our work is discussed in \cite[\S4.3]{mcrae2025rigidity}.


\section{Two constructions of commutative exact algebras}\label{sec:commutative-algebras}
In this section we provide two constructions of commutative exact algebras in a braided FTC $\C$ over an algebraically closed field $\kk$. The first construction is based on the notion of a simple current algebra. The second construction is based on using adjoints of tensor functors.


\subsection{Simple current algebras}\label{subsec:simple-currents}
In this section, we use direct sums of simple invertible objects to get examples of commutative exact algebras. This is a well-known technique in the semisimple setting, see \cite{fuchs2004tft} for example. We show that this technique can be generalized to the non-semisimple setting and demonstrate it in a few examples.

\subsubsection{Basic setup}
Let $\C$ be a braided finite tensor category over $\kk$. Denote the set of invertible simple objects in $\C$ as $\Inv(\C)$. As $\C$ is braided, $G:=\Inv(\C)$ is an abelian group under the tensor product. For $g\in G$, denote the corresponding object in $\C$ as $X_g$. Then the semisimple subcategory of $\C$ spanned by $\Inv(\C)$ will be a pointed fusion category of the form $\C(G,q)$ where $q:G\rightarrow \kk^{\times}$ is a quadratic form determined by $c_{X_g,X_g} = q(g)\id_{X_{g^2}}$ for all $g\in G$ \cite[\S8.4]{etingof2016tensor}. 

\begin{lemma}
Let $H\subset \Inv(\C)$ be a subgroup and consider the object $A_H=\oplus_{h\in H}X_h$.
\begin{enumerate}
  \item If $c_{X_h,X_h}=\id$ for all $h\in H$, then $A_H$ is a haploid, commutative and exact algebra in $\C$. Hence, $\C_A^{\loc}$ is a braided FTC.
  \item If, in addition, $\C$ is ribbon FTC (or MTC) and $\theta_{X_h}=\id_{X_h}$ for all $h\in H$, then $A_H$ is, in addition, symmetric Frobenius. Hence, $\C_A^{\loc}$ is a ribbon FTC (or MTC).
\end{enumerate} 
\end{lemma}
\begin{proof}
By the discussion in Appendix~\ref{app:pointed}, $A_H$ is a haploid, commutative and exact Frobenius algebra in the category $\C(\Inv(\C),q)$. As $\C(\Inv(\C),q)\hookrightarrow \C$ is a fully faithful braided tensor functor, by Lemma~\ref{lem:F-preserves}, the resulting algebra in $\C$ is also exact Frobenius, commutative and haploid. Thus, part (a) is proved. 

For (b), the assumption implies that $\theta_A=\id$. Also, as $A$ is Frobenius, $A^*\cong A\in(\C_A^{\loc})'$. Thus, by Theorem~\ref{thm:FTC-main-result}, the claim follows.
\end{proof}

We call such an algebra $A$ that is a direct sum of simple invertible objects in $\C$ a \textit{simple current algebra}. 
We note one particular case: if $X\in\Inv(\C)$ satisfies $X^{\otimes n}\cong\unit$ and $c_{X,X}=\id$, then $A=\unit\oplus X \oplus \ldots \oplus X^{n-1}$ is a haploid commutative exact algebra in $\C$. If $\theta_X=\id_X$, then $\theta_A=\id_A$ and $\C_A^{\loc}$ is a ribbon FTC.


\subsubsection{Hopf algebras setting}
Let $(H,R)$ be a finite-dimensional quasi-triangular Hopf algebra with $R$-matrix $R=R_i\otimes R^i$. Let $\C$ denote the finite tensor category $\Rep(H)$. For $V,W\in\Rep(H)$, the braiding is given by
\[ c_{V,W}: V\otimes W\rightarrow W\otimes V, \quad (v\otimes w) \mapsto (R^i\cdot w) \otimes (R_i\cdot v). \]
Given a character $\alpha:H\rightarrow\kk$, let $\kk_{\alpha}$ denote the $H$-module which is $\kk$ as a vector space and $H$-action is given by $h\cdot c = \alpha(h)c$ for $h\in H,c\in\kk$.
The abelian group of invertible objects is 
\[ \Inv(\C) = \{ \kk_{\alpha} \mid \alpha \; \text{is a character of } \; H\}.\]
Note the $\Inv(\C)$ is isomorphic to the group of grouplike elements $G(H^*)$ of the dual Hopf algebra $H^*$.
Using the above mentioned formula for the braiding, we get that quadratic form $q$ on $\Inv(\C)$ is given by 
\[ q(\kk_{\alpha}) = \alpha(R_i)\alpha(R^i) . \]
Thus any subgroup $H\subset \Inv(\C)$ such that $q(\kk_{\alpha})=1$ for all $\kk_{\alpha}\in H$ gives a commutative exact algebra.

One can also take the dual point of view as follows. Let $(H,r)$ be a (possibly infinite-dimensional) coquasi-triangular Hopf algebra with bijective antipode. Here $r:H\otimes H\rightarrow \kk$ is a map called the $r$-form that makes the category $\C=\mathrm{Corep}(H)$, that is the category of finite-dimensional left $H$-comodules, a braided tensor category \cite[\S8.3]{etingof2016tensor}. Given $U,V\in\mathrm{Corep}(H)$, the braiding is given by
\[c_{U, V}: U \otimes V \rightarrow V \otimes U, \quad u \otimes v \mapsto \sum r\left(u_{(1)}, v_{(1)}\right) v_{(0)} \otimes u_{(0)}.\]
The abelian group of invertible objects $\Inv(\C)$ coincides with the group $G(H)$ of grouplike elements of $H$. The quadratic form $q$ on $\Inv(\C)$ is given by $q(g) = r(g,g)$ for $g\in G(H)$. Thus, any subgroup $H\subset G(H)$ such that $r(h,h)=1$ for all $h\in H$ yields a commutative, haploid and exact algebra in $\C$.

\begin{remark}
We note that a similar generalization also works for quasi-Hopf algebras and we leave it for the interested reader.
\end{remark}


\subsubsection{Example 1: $u_q^{\phi}(\fsl_2)$ at $q$ even order}\label{subsec:ex1-even-order}
In this subsection, $\kk=\mathbb{C}$. Recall the description of the modified small quantum group at even order root of unity $q=e^{\pi \iota/p}$ from \cite{creutzig2020quasi}. It is a quasi-Hopf algebra with $3$ generators $E,F,K$. 
It has two invertible objects, the unit object $\unit$ and another object $\psi$ which is $\kk$ as a vector space with action of the generators given by $E,F,K$ given by $0,0,-1$ respectively. Also, $\psi\otimes \psi\cong\unit$, $c_{\psi,\psi}=\iota^p\, \id_{\unit}$ and $\theta_{\psi}=-\iota^p \, \id_{\psi}$. See Appendix~\ref{app:even-order-uqsl2} for definitions and proofs. 

When $p$ is a multiple of $4$, $c_{\psi,\psi}=\id_{\unit}$ and $\theta_{\psi}=-\id_{\psi}$. Thus, the algebra $A=\unit\oplus \psi$ is a haploid, commutative and exact algebra in $\C$. It is Frobenius (thus $A^*\in(\C_A^{\loc})')$ but not symmetric Frobenius because $\theta_A\neq\id_A$. Consequently, one may expect that the braided finite tensor category $\C_A^{\loc}$ is not pivotal. This is indeed the case, see \cite[Example~5.2]{etingof2020frobenius}.

We next discuss how the category $\C_A^{\loc}$ considered above is the same as a category previously studied in \cite[\S10.4]{negron2021log} and \cite[Example~5.13]{etingof2020frobenius}.

\subsubsection{Another interpretation}
Let $\C$ be a non-degenerate braided finite tensor category and $\D\subset\C$ be a full braided tensor subcategory. 
Then, by \cite[Theorem~4.3]{laugwitz2023constructing}, we know that $\Z_{(2)}(\D\subset\C)' = \D'$. Let $A$ be a haploid, commutative and exact algebra in $\D'$. Using the embedding
\[\Z_{(2)}(\D\subset\C)' \subset \Z_{(2)}(\D\subset\C) \subset \C ,\] 
$A$ is also a haploid, commutative, exact algebra in $\Z_{(2)}(\D\subset\C)$ and $\C$ (by Lemma~\ref{lem:exact-preserve}). Thus, we can form the category $\Z_{(2)}(\D\subset\C)_A$. Observe that we have a natural embedding
\begin{equation*}
 \Z_{(2)}(\D\subset\C)_A = \Z_{(2)}(\D\subset\C)_{A}^{\loc} \subset \C_A^{\loc}.
\end{equation*}

\begin{lemma}\label{lem:EN-example}
  $\Z_{(2)}(\D\subset\C)_A \cong \C_A^{\loc} \;\; \iff \;\; \FPdim(\D) = \FPdim(A)$. 
\end{lemma}
\begin{proof}
We have that 
\[ \FPdim \left(\Z_{(2)}(\D\subset\C)_A \right) 
= \frac{\FPdim\left(\Z_{(2)}(\D\subset\C)\right)}{\FPdim(A)} 
= \frac{\FPdim(\C)}{\FPdim(\D) \FPdim(A)}. \]
Here the first equality uses Lemma~\ref{lem:FPdim-C_A} and the second \cite[Theorem~4.9]{shimizu2019non}.
On the other hand, by Lemma~\ref{lem:FPdim-C-A-loc}, $\FPdim(\C_A^{\loc}) = \FPdim(\C)/\FPdim(A)^2$. As $\Z_{(2)}(\D\subset\C)_A$ is a tensor subcategory of $\C_A^{\loc}$, they are equivalent if and only if their Frobenius-Perron dimensions match. Now the claim follows immediately.
\end{proof}

A special case of this consideration is when the symmetric tensor category $\D$ is \textit{Tannakian}, namely, $\D\cong \Rep(G)$ for some finite group $G$. In this case, $\D'=\D$ and one can take $A$ to be the function algebra $\mathcal{O}(G) \in \D'$. Then $\FPdim(A)=\FPdim(\D)$ and by Lemma~\ref{lem:EN-example}, it follows that $\Z_{(2)}(\D\subset\C)_A \cong \C_A^{\loc}$. 

Now the coincidence of our example and that of \cite{negron2021log}, \cite{etingof2020frobenius} is clear. Indeed, while we consider the category $\C_A^{\loc}$, they consider the category $\Z_{(2)}(\D\subset\C)_A$ where $\D$ is the category of representations of $\mathbb{Z}_2$.


\subsubsection{Example 2: Double of Taft algebra}
In this subsection, $\kk=\mathbb{C}$.
Fix a positive integer $n\geq 2$ and take $q=e^{2\pi \iota/n}$.
Let $T_n$ denote the Taft algebra and $D_n$ its Drinfeld double. Then, $D_n$ is a quasitriangular Hopf algebra with generators $a,b,c,d$ subject to the following relations (see \cite{chen1999class}):
\begin{align*}
  b a=q a b, \quad d b=q b d, \quad b c=c b, &\quad c a=q a c, \quad d c=q c d, \quad d a-q a d=1-b c, \\
  a^n=0=d^n, &\quad b^n=1=c^n .
\end{align*}
In this example, we will consider the category $\C=\Rep(D_n)$. 
The category $\C$ is a nondegenerate braided FTC with $\FPdim(\C)=n^4$.

From \cite[\S2]{benkart2022tensor}, we recall that $\Inv(\C)=\{ V(1,s) \mid s\in\mathbb{Z}_n\}$, where $V(1,s)$ is the vector space $\kk$ such that $a$, $b$, $c$ and $d$ act on $V(1, s)$ by the scalars $0$, $q^s$, $q^{-s}$ and $0$, respectively. Also, $V(1,s)\otimes V(1,r)\cong V(1,s+r)$. Thus, $\Inv(\C)\cong\mathbb{Z}_n$ with generator $X:=V(1,1)$.
The $R$-matrix of $D_n$ is $R=\frac{1}{n} \sum_{m, s, t=0}^{n-1} \frac{q^{-t m}}{[s]!} a^s b^t \otimes c^m d^s$ \cite[(3.11)]{benkart2022tensor}.
Using this, we get that 
\begin{align*}
  c_{X,X} & = 
\left[ \frac{1}{n} \sum_{m, s, t=0}^{n-1} \frac{q^{-t m}}{[s]!} [(a^s b^t)\cdot 1] [(c^m d^s)\cdot 1] \right]\, \id 
\stackrel{(s=0)}{=} \left[ \frac{1}{n} \sum_{m, t=0}^{n-1} q^{-t m} (b^t\cdot 1) (c^m\cdot 1)\right] \id \\
& = \left[\frac{1}{n} \sum_{m, t=0}^{n-1} q^{-t m} q^t q^{-m} \right] \id  
= \left[\frac{1}{n} \sum_{m=0}^{n-1} q^{-m} \left(\sum_{t=0}^{n-1} q^{t(1-m)}\right) \right] \id  
\end{align*}
Since the sum in the inner brackets is $n$ if $m=1$ and $0$ otherwise (by Lemma~\ref{lem:sum-roots-of-unity}), we get that $c_{X,X} = q^{-1}\id$. This implies that $c_{X^m,X^m} = (q^{-1})^{m^2}\id$. Also note that $\C$ admits a ribbon structure if and only if $n$ is odd \cite[Prop.~7]{kauffman1993necessary} (which makes it an MTC). The ribbon structure has been explicitly described in \cite[Theorem~3.6(b)]{benkart2022tensor}. Using it, one can calculate that $\theta_X=q^{-1} \id_X$. Using induction, one can prove that $\theta_{X^m}=q^{-m^2}\id_{X^m}$ for all $m\geq 1$.
Next we discuss a special case.

\underline{$n=m^2$ for some integer $m>1$}: 
As $c_{X^m,X^m} = (q^{-1})^{m^2}\id = (q^{-1})^{n}\id = \id$, it follows that $A = \bigoplus_{k=0}^{m-1} X^{mk}$ is a haploid, commutative and exact algebra in $\C$. As $\FPdim(A)=m$, we have that $\FPdim(\C_A^{\loc})=n^4/m^2 = m^6$. Now suppose that $n$ is odd. As $\theta_{X^{mk}} = q^{-m^2k^2}\id = q^{-n k^2}\id =\id$ for all $k$, we get that $\theta_A=\id$ and hence it is symmetric Frobenius. Thus, $\C_A^{\loc}$ is a MTC with $\FPdim(\C_A^{\loc})=m^6$.


\subsubsection{Example~3: Tensor products of \texorpdfstring{$u_q^{\phi}(\fsl_2)$}{uqsl2-phi}}\label{subsec:ex3-Deligne}
In this subsection, $\kk=\mathbb{C}$.
Take $n$ integers $p_1,\cdots,p_n$ where each $p_i\geq 2$. Set $\D_{p_i} = \Rep(u_q^{\phi}(\fsl_2))$ for $q=e^{{\pi\iota}/{p_i}}$ and consider $\C=\D_{p_1} \boxtimes \cdots \boxtimes \D_{p_n}$. From Section~\ref{subsec:ex1-even-order}, recall that $\Inv(\D_{p_i})=\{\unit,\psi\} \cong\mathbb{Z}_2$. 
Thus, the category $\C$ satisfies 
\[\Inv(\C)\cong \mathbb{Z}_2^{n} = \{ \psi^{i_1} \boxtimes \cdots \boxtimes \psi^{i_n} | i_k\in\{0,1\} \; \forall \; 1\leq k \leq n \}. \]
Using the formula for self-braiding and twist of $\psi$ proved in Appendix~\ref{app:even-order-uqsl2}, we get that the quadratic form and twist on objects of $\Inv(\C)$ are  
\[ q(\psi^{i_1} \boxtimes \cdots \boxtimes \psi^{i_n}) = \iota^{i_1p_1+\cdots+i_np_n} ,
\qquad
\theta_{\psi^{i_1} \boxtimes \cdots \boxtimes \psi^{i_n}} = (-1)^{i_1+ \cdots+i_n} \iota^{i_1p_1+\cdots+i_np_n}.\] 
Now, consider the subset $H'\subset \Inv(\C)$ given by
\[ H' = \{(i_1,\ldots,i_n)\in \mathbb{Z}_2^{n} \, | \,2\; \text{divides}\; (i_1+\ldots+i_n) \; \text{and} \; 4 \; \text{divides} \; i_1p_1+\cdots+i_np_n  \} . \]
Clearly, the quadratic form is identically $1$ on $H'$ and the twist is $\id$. Thus, for any subgroup $H\leq \mathbb{Z}_2^n$ such that $H\subset H'$, $A_H$ is a haploid, commutative, exact symmetric Frobenius algebra in $\C$ and $\C_{A_H}^{\loc}$ is an MTC.

\begin{remark}
  Similar examples have also been considered in the VOA literature, see for example \cite[\S4.1.2]{creutzig2017tensor}.
\end{remark}

Next, we discuss a special case of the above example. 
Let $p=4k+2$ for some integer $k\geq 0$ and fix $\D=\Rep(u_q^{\phi}(\fsl_2))$ where $q=e^{\pi \iota/p}$.
Then, $c_{\psi,\psi}=\iota^{4k+2}\id_{\unit} = -\id_{\unit}$ and $\theta_{\psi}=-\iota^{4k+2}\id_{\psi} = \id_{\psi}$. 

Consider the MTC $\C = \D^{\boxtimes d}$ for some $d\geq 2$. Then, $\Inv(\C)\cong \mathbb{Z}_2^{d}$ and the quadratic form is $q(\psi^{i_1} \boxtimes \cdots \boxtimes \psi^{i_d}) = (-1)^{i_1+\ldots+i_d}$ and $\theta_{\psi^{i_1} \boxtimes \cdots \boxtimes \psi^{i_d}} = \id_{\psi^{i_1} \boxtimes \cdots \boxtimes \psi^{i_d}}$.
Consider the subgroup $H\subset \mathbb{Z}_2^{d}$ given by
\[ H = \{ (i_1,\ldots,i_d)\in \mathbb{Z}_2^{d} \, | \, i_1+i_2+\ldots+i_d \;\text{is even} \}. \]
Since $q|_H=1$ and $\theta$ is identity for all elements in $H$, we get that $A_H := \oplus_{i_1+\ldots+i_d\; \text{even}} \psi^{i_1} \boxtimes \cdots \boxtimes \psi^{i_d}$ is a haploid, commutative exact symmetric Frobenius algebra in $\C$. Thus, $\C_{A_H}^{\loc}$ is an MTC. 


\subsection{Commutative exact algebras via adjoints}\label{subsec:comm-alg-via-adjoint}

In this section, we discuss a way to construct commutative exact algebras via adjoints.
We first recall a construction of a central commutative algebra established in the theory of Hopf monads \cite{bruguieres2011hopf}.
Let $F : \C \to \D$ be a strong monoidal functor between rigid monoidal categories $\C$ and $\D$ with monoidal structure $\phi^2 : F(-) \otimes F(-) \to F(- \otimes -)$ and $\phi^0 : \unit \to F(\unit)$.
We assume that $F$ has a right adjoint $R$ and let $\eta$ and $\varepsilon$ be the unit and the counit of the adjunction $F \dashv R$, respectively.
Then we have a natural transformation
\begin{equation}
  \label{eq:Hopf-operator-r}
  \mathbb{H}^r_{X,W} := R((\id_{F(X)} \otimes \varepsilon_W) \circ (\phi^2_{X, R(W)})^{-1}) \circ \eta_{X \otimes R(W)} : X \otimes R(W) \to R(F(X) \otimes W)
\end{equation}
for $X \in \C$ and $W \in \D$, which is identical to the right Hopf operator of the adjunction $R^{\op} \dashv F^{\op}$ in the sense of \cite[2.8]{bruguieres2011hopf}. Similarly, we define
\begin{equation}
  \label{eq:Hopf-operator-l}
  \mathbb{H}^l_{W,X} := R( (\varepsilon_W \otimes \id_{F(X)}) \circ (\phi^2_{R(W),X})^{-1} ) \circ \eta_{R(W) \otimes X}
\end{equation}
for $X \in \C$ and $W \in \D$, which is identical to the left Hopf operator of $R^{\op} \dashv F^{\op}$.
Since $\C$ and $\D$ are assumed to be rigid, $R^{\op} \dashv F^{\op}$ is a Hopf adjunction \cite[Proposition 3.5]{bruguieres2011hopf}.
Thus, the natural transformations $\mathbb{H}^{l}$ and $\mathbb{H}^r$ are invertible.

The functor $R$ is lax monoidal as a right adjoint of $F$. Thus, in particular, $A := R(\unit)$ is an algebra in $\C$. By applying \cite[Theorem 6.6]{bruguieres2011hopf} to the Hopf adjunction $R^{\op} \dashv F^{\op}$ (see also \cite[Theorem 6.1]{bruguieres2011exact}), we find that the algebra $A$ is in fact a central commutative algebra in $\C$ together with the half-braiding given by
\begin{equation}
  \label{eq:central-comm-alg-half-braiding}
  \sigma_{X} = (\mathbb{H}^r_{X,\unit})^{-1} \circ \mathbb{H}^l_{\unit,X} : A \otimes X \to X \otimes A \quad (X \in \C).
\end{equation}

Now we assume that $\C$ has a braiding. In general, $A$ is not a commutative algebra in $\C$. Under an additional assumption that $F$ is central, one can show that the half-braiding of $A$ comes from the braiding of $\C$ and hence $A$ is a commutative algebra in $\C$. The precise statement is as follows:

\begin{proposition}
  \label{prop:central-comm-alg-half-braiding}
  Let $F: \C \to \D$ be a strong monoidal functor between rigid monoidal categories $\C$ and $\D$. We assume that $\C$ has a braiding $c$ and $F$ is central. Then the half-braiding \eqref{eq:central-comm-alg-half-braiding} of the algebra $A = F^{\ra}(\unit)$ is given by
  \begin{equation*}
    \sigma_X = c_{X,A}^{-1}
  \end{equation*}
  for all $X \in \C$. Hence, $A$ is in fact a commutative algebra in $\C$.
\end{proposition}

The commutativity of $A$ can be proved in the same way as \cite[Proposition 8.8.8]{etingof2016tensor}.
The actual purpose of this proposition is to provide a way to connect the results on central commutative (co)algebras in \cite{bruguieres2011hopf} and the results on commutative algebras in braided finite tensor categories given in \cite{etingof2016tensor}.

\begin{proof}
  We fix $\widetilde{F} : \C \to \Z(\D)$ be a braided tensor functor such that $U_{\D} \circ \tilde{F} = F$ as tensor functors. Since $U \circ \widetilde{F} = F$, we may write
  \begin{equation*}
    \widetilde{F}(X) = (F(X), \ \widetilde{c}_{F(X)} : F(X) \otimes (-) \to (-) \otimes F(X))
  \end{equation*}
  for $X \in \C$. Since $\widetilde{F}$ is braided, we have
  \begin{equation}
    \label{eq:central-comm-alg-F-central}
    F(c_{X,Y}) \circ \phi^2_{X,Y} = \phi^2_{Y,X} \circ \widetilde{c}_{F(X), F(Y)}
    \quad (X, Y \in \C),
  \end{equation}
  where $\phi^2$ is the monoidal structure of $F$. Now we compute
  \begin{align*}
    & \mathbb{H}^r_{X, \unit} \circ \sigma_X \circ c_{X,A} \\
    {}^{\eqref{eq:Hopf-operator-l}, \eqref{eq:central-comm-alg-half-braiding}}
    & = F^{\ra}((\varepsilon_{\unit} \otimes \id_{F(X)}) \circ (\phi^2_{A,X})^{-1}) \circ \eta_{A \otimes X} \circ c_{X, A} \\
    {}^{\text{(nat.)}}
    & = F^{\ra}((\varepsilon_{\unit} \otimes \id_{F(X)})
      \circ (\phi^2_{A,X})^{-1} \circ F(c_{X,A})) \circ \eta_{A \otimes X} \\
    {}^{\eqref{eq:central-comm-alg-F-central}}
    & = F^{\ra}((\varepsilon_{\unit} \otimes \id_{F(X)})
      \circ \widetilde{c}_{F(X), F(A)} \circ (\phi^2_{X,A})^{-1}) \circ \eta_{X \otimes A} \\
    {}^{\text{(half-br.)}}
    & = F^{\ra}(\widetilde{c}_{F(X), \unit} \circ (\id_{F(X)} \otimes \varepsilon_{\unit})
      \circ (\phi^2_{X,A})^{-1}) \circ \eta_{X \otimes A} \\
    {}^{\text{(half-br.)}}
    & = F^{\ra}((\id_{F(X)} \otimes \varepsilon_{\unit})
      \circ (\phi^2_{X,A})^{-1}) \circ \eta_{X \otimes A}
      \mathop{=}^{\eqref{eq:Hopf-operator-r}}
      \mathbb{H}^r_{X,\unit}
  \end{align*}
  for all $X \in \C$, where `(nat.)' and `(half-br.)' follow from the naturality of $\eta$ and the axiom of half-braiding, respectively. Hence we have $\sigma_X = c_{X,A}^{-1}$, as desired.
\end{proof}

\begin{proposition}
  \label{prop:comm-alg-via-adjoint-1}
  Let $F: \C \to \D$ be a tensor functor between finite tensor categories $\C$ and $\D$. Then the central commutative algebra $A := F^{\ra}(\unit)$ enjoys the following properties:
  \begin{enumerate}
  \item[\textup{(a)}] $A$ is an indecomposable exact algebra in $\C$.
  \item[\textup{(b)}] $A$ is a central commutative algebra in $\C$ with the half-braiding $\sigma$ given by \eqref{eq:central-comm-alg-half-braiding}.
  \item[\textup{(c)}] The functor $F^{\ra}$ induces an equivalence $K: \mathrm{Im}(F) \to \C_A^{\sigma}$ of tensor categories such that $K \circ F \cong F_A$ as tensor functors, where $F_A : \C \to \C_A^{\sigma}$ is the free module functor.
  \item[\textup{(d)}] $\FPdim_{\C}(A) = \FPdim(\C) / \FPdim(\mathrm{Im}(F))$.
  \end{enumerate}
\end{proposition}
\begin{proof}
This proposition is obtained as a combination of results found in \cite{bruguieres2011hopf,etingof2016tensor}. Here are the details:
let $f : \C \to \mathrm{Im}(F)$ be the corestriction of $F$, and let $r$ be the restriction of $F^{\ra}$ to $\mathrm{Im}(F)$. It is obvious that $f$ is a surjective tensor functor and $r$ is right adjoint to $f$ such that $A = r(\unit)$.
By Lemma~\ref{lem:dominant-and-surjective}, the tensor functor $f$ is dominant.
Now (b) and (c) immediately follow from \cite[Proposition 6.1]{bruguieres2011hopf}.
Since $\mathrm{Im}(F)$ is a finite tensor category, so is $\C_A^{\sigma}$.
In particular, $A$ is a simple object of $\C_A^{\sigma}$. This means that $A$ is indecomposable. Theorem~\ref{thm:C-A-rigid} implies that $A$ is also exact. Hence we have proved (a).
Part (d) is \cite[Lemma 6.2.4]{etingof2016tensor}.
\end{proof}

As a corollary to the above results, we obtain the following:

\begin{corollary}\label{cor:comm-alg-via-adjoint-2}
  Let $F: \C \to \D$ be a central tensor functor between finite tensor categories where $\C$ braided.
  Then $A := F^{\ra}(\unit)$ is a haploid commutative exact algebra in $\C$.
  Moreover, there is an equivalence $K: \mathrm{Im}(F) \to \C_A^{\sigma}$ (for $\sigma=c_{-,A}^{-1}$) induced by $F^{\ra}$ such that $K \circ F \cong F_A$, where $F_A : \C \to \C_A^{\sigma}$ is the free module functor. \qed
\end{corollary}

Next we discuss the unimodular case. 

\begin{proposition}\label{prop:comm-alg-via-adjoint-Frob}
  Let $F: \C \to \D$ be a tensor functor between finite tensor categories $\C$ and $\D$ that preserves projectives. If $\C$ and $\D$ are unimodular, then $A:= F^{\ra}(\unit)$ is Frobenius.
\end{proposition}
\begin{proof}
First note that $A$ is Frobenius iff the forgetful functor $ \C_A^{\sigma}\rightarrow \C$ is Frobenius \cite[\S2.2]{shimizu2016unimodular}. As the forgetful functor is right adjoint to $F_A$ and $F_A \cong K\circ F$ with $K$ an equivalence, it is enough to show that $r=F^{\ra}|_{\mathrm{Im}(F)}:\mathrm{Im}(F)\rightarrow\C$ is Frobenius. For this, it suffices to show that $F^{\ra}$ is Frobenius, which we show next. 
As $\C,\D$ are finite and $F$ preserves projectives, $F^{\ra}$ is exact. Now using \cite[Proposition~4.19]{fuchs2020eilenberg}, it follows that $F^{\rra} \cong F$. In other words, $F^{\ra}$ is Frobenius.
\end{proof}

\newcommand{\Corep}{\mathrm{Corep}}

\subsubsection{Hopf subalgebras}
Let $H$ be a f.d. quasitriangular Hopf algebra and $i:H'\hookrightarrow H$ be a Hopf subalgebra such that the inclusion $i$ factors through the Double $D(H')$ of $H$. This yields a surjective tensor functor $F:\Rep(H)\rightarrow\Rep(H')$, which preserves projectives and is central. Thus, it yields a commutative exact algebra $A=F^{\ra}(\kk)$ in $\Rep(H)$. We discuss a particular example of this below.

Dually, if we have a surjective map of Hopf algebras $K\rightarrow H$, then we get a surjective tensor functor $F:\Corep(K)\rightarrow\Corep(H)$. Suppose furthermore that, $K$ is a \textit{braided central Hopf subalgebra} of $H$ as defined in \cite[Definition~3.1]{angiono2014equivariantization}. Then, by Theorem~3.4 of \cite{angiono2014equivariantization}, $F$ is central. Thus, we get a commutative exact algebra $A=F^{\ra}(\kk)$ in $\Corep(H)$.

\begin{example}
Let $G$ be a finite group scheme over a field $\kk$. Equivalently we can work with the corresponding commutative Hopf algebra $\cO(G)$ of functions on $G$. From this, we get a symmetric finite tensor category $\Rep(G)\cong \Corep(\cO(G))$.
Let $H$ be a subscheme of $G$ for which the induction functor from $H$ to $G$ is exact and faithful (by \cite{cline1977induced}, this happens if and only if the quotient scheme $G/H$ is affine). This induces a map $\cO(H)\rightarrow\cO(G)$ of commutative Hopf algebras.
Then, $A=F^{\ra}(\kk)$ is a commutative, exact algebra in $\Rep(G)$. 
\end{example}


\subsubsection{Lagrangian commutative algebras in the Drinfeld center}\label{subsec:comm-alg-in-ZC}
Let $\C$ be a finite tensor category, and let $\M$ be an indecomposable exact left $\C$-module category.
Then the dual tensor category $\C_{\M}^*$ is defined.
There is a canonical equivalence $\Z(\C) \simeq \Z(\C_{\M}^*)$ of braided tensor categories.
By composing this equivalence and the forgetful functor $\Z(\C_{\M}^*) \to \C_{\M}^*$, we obtain a central surjective tensor functor $\Psi : \Z(\C) \to \C_{\M}^*$.
By Proposition~\ref{prop:comm-alg-via-adjoint-1}, $\mathbf{A}_{\M} := \Psi^{\ra}(\id_{\M})$ is a haploid, commutative and exact algebra in $\Z(\C)$. The algebra $\mathbf{A}_{\M}$ is Lagrangian. Indeed, again by Proposition~\ref{prop:comm-alg-via-adjoint-1}, we have
\begin{equation*}
\FPdim_{\Z(\C)}(\mathbf{A}_{\M})
= \frac{\FPdim(\Z(\C))}{\FPdim(\C_{\M}^*)}
= \frac{\FPdim(\C)^2}{\FPdim(\C)}
= \FPdim(\C),
\end{equation*}
where the second equality follows from the formula of $\FPdim(\Z(\C))$ and the Morita invariance of the FP dimension \cite[Corollary~7.16.7]{etingof2016tensor}.

\begin{remark}\label{rem:Z(C)-lagrangian}
    When $\M=\C$, the functor $\Psi:\Z(\C)\to\C^*_\M$ is the forgetful functor $U_{\C}:\Z(\C)\to\C$. Denoting its right adjoint as $R_{\C}$, we obtain that $R_\C(\unit)$ is a Lagrangian algebra in $\Z(\C)$. Using Corollary~\ref{cor:lagrangian}, $\Z(\C)_{R_\C(\unit)}^\loc\simeq \Vect$.
\end{remark}

\begin{example}\label{eg:uqg-Nsquare-order}
  By \cite{etingof2009small}, when $q$ is a root of unity of order $N^2$ for $N$ odd, we have a braided tensor equivalence $F: \Rep(u_q(\fg)) \xrightarrow{\sim} \Z(\D)$ for some finite tensor category $\D$. Thus, in this case, $\Rep(u_q(\fg))$ admits a Lagrangian algebra $\mathbf{A}'_{\M}$ for every exact indecomposable $\D$-module category $\M$, namely $\mathbf{A}_{\M}' = F^{-1}(\mathbf{A}_{\M})$. 
\end{example}

The algebra $\mathbf{A}_{\M}$ is not a special Frobenius algebra in $\Z(\C)$ in general. Indeed, for $\M = \C$, we can prove:

\begin{lemma}
  Let $\C$ be a pivotal unimodular finite tensor category, and let $R$ denote the right adjoint of the forgetful functor $U:\Z(\C)\rightarrow\C$.
  The algebra $\mathbf{A} := R(\unit)$ in $\Z(\C)$ is a symmetric Frobenius, haploid, commutative and exact algebra in $\Z(\C)$.
  The algebra $\mathbf{A}$ is a special Frobenius algebra if and only if $\C$ is semisimple and $\dim(\C) \ne 0$.
\end{lemma}

Here, $\dim(\C)$ is defined by $\dim(\C) = \sum_{i = 1}^m \dim(V_i^*) \dim(V_i)$, where $V_1, \cdots, V_m$ are the set of representatives of the isomorphism classes of simple objects of $\C$ and $\dim(X)$ for an object $X \in \C$ means the pivotal dimension of $X$.

\begin{proof}
By the above discussion, the algebra $\mathbf{A}$ is haploid, exact and commutative.
By the unimodularity, $\mathbf{A}$ is Frobenius \cite[Theorem 5.6]{shimizu2016unimodular}.
It is also known that $\mathfrak{p}_{\mathbf{A}} = \mathfrak{u}_{\mathbf{A}}$ \cite[Theorem 6.1 and Proposition 6.4]{shimizu2015pivotal}.
Hence, by Lemma \ref{lem:symmetric-frobenius-theta}, $\mathbf{A}$ is a symmetric Frobenius algebra.

To prove the latter part, we assume that $\mathbf{A}$ is a special Frobenius algebra. Let $\varepsilon : U R \to \id_{\C}$ be the counit of the adjunction, and set $A = U(\mathbf{A})$. The morphism $\varepsilon_{\unit} : A \to \unit$ is in fact a morphism of algebras in ${}_A\C_A$. If we view $\unit$ as an $A$-bimodule in $\C$ via $\varepsilon_{\unit}$, then $\varepsilon_{\unit} : A \to \unit$ is a morphism of $A$-bimodules in $\C$. It is obvious that $\varepsilon_{\unit}$ splits as a morphism in $\C$. Since $A$ is a separable algebra, $\varepsilon_{\unit}$ splits as a morphism in ${}_A\C_A$ \cite[Theorem 1.30]{ardizzoni2007hochschild}.
Now let $\Lambda : \unit \to A$ be a morphism of $A$-bimodules in $\C$ such that $\varepsilon_{\unit} \circ \Lambda = \id_{\unit}$, which exists by the above argument. The morphism $\Lambda$ is a categorical integral introduced in \cite{shimizu2019integrals}. By an analogue of Maschke's theorem for finite tensor categories \cite[Theorem~4.12]{shimizu2019integrals}, we conclude that $\C$ is semisimple. The semisimplicity implies that there is an isomorphism $A \cong \bigoplus_{i = 1}^m V_i^* \otimes V_i$ in $\mathcal{C}$ (see \cite[\S6.1]{shimizu2017monoidal}).
By Lemma \ref{lem:special-Frobenius-trace-Nakayama} and that $\mathbf{A}$ is symmetric Frobenius, we have $\dim(\mathcal{C}) = \dim(A) = \dim(\mathbf{A}) \ne 0$. We have proved the `only if' part.

We shall prove the converse.
We assume that $\C$ is semisimple and $\dim(\C) \ne 0$.
As we have mentioned in the above discussion, $\dim(\C) = \dim(\mathbf{A})$.
Hence, by Lemma \ref{lem:special-Frobenius-trace-Nakayama}, $\mathbf{A}$ is special Frobenius. The proof is done.
\end{proof}


\subsubsection{A commutative algebra in the relative Drinfeld center}\label{subsec:Yetter-Drinfeld-relative}
Let $\B$ be a braided finite tensor category, and let $\C$ be a finite tensor category such that $\Z(\C)$ contains $\B$ as a braided tensor full subcategory. Then the category $\Z_{\B}(\C)=\Z_{(2)}(\B\subset \Z(\C))$, called the ($\B$-)relative Drinfeld center of $\C$ is defined. According to \cite{laugwitz2022relative}, $\Z_{\B}(\C)$ is a tensor subcategory of $\Z(\C)$ with FP dimension
\begin{equation*}
    \FPdim(\Z_{\B}(\C)) = \frac{\FPdim(\C)^2}{\FPdim(\B)}.
\end{equation*}
The forgetful functor $U: \Z_{\B}(\C) \to \C$ is a central tensor functor and therefore, $A := U^{\ra}(\unit)$ is an indecomposable exact commutative algebra in $\Z_{\B}(\C)$.
If $U$ is surjective, then 
\begin{equation}\label{eq:FPdim-ZBC}
    \Z_{\B}(\C)_A \simeq \C \quad \text{and} \quad \FPdim_{\Z_{\B}(\C)}(A) = \frac{\FPdim(\C)}{\FPdim(\B)}.
\end{equation}
Furthermore, we have that
\[\Z_{\B}(\C)_A^{\loc} \simeq \Z_{(2)}(\Z_{\B}(\C) \subset \Z(\Z_{\B}(\C)_A)) \simeq \Z_{(2)}(\Z_{\B}(\C) \subset \Z(\C)) \simeq \B.\]
as braided tensor categories. Here the first equivalence is by Theorem~\ref{thm:LW-local}, the second by \eqref{eq:FPdim-ZBC}, and the third by \cite[Theorem~4.9]{shimizu2019non}.

Next we discuss the case of Yetter-Drinfeld modules.


\subsubsection{A commutative algebra in the Yetter-Drinfeld category}

Let $\B$ be a braided finite tensor category, and let $K$ be a Hopf algebra in $\B$.
We take $\C$ to be the category of left $K$-modules in $\B$, which is a finite tensor category of FP dimension $\FPdim(\B) \cdot \FPdim_{\B}(K)$.
The category $\B$ is identified with the tensor full subcategory of $\C$ consisting of all $K$-modules where the action is given through the counit of $K$. The $\B$-relative Drinfeld center $\Z_{\B}(\C)$ is equivalent to the category ${}^K_K\YD(\B)$ of Yetter-Drinfeld modules over $K$ in $\B$ \cite{laugwitz2020relative}.
Thus, as explained above, we obtain an indecomposable exact commutative algebra $A = F^{\ra}(\unit)$ in ${}^K_K\YD(\B)$ by using a right adjoint of the forgetful functor $F: {}^K_K\YD(\B) \to \C$. The FP dimension of $A$ is
\begin{equation*}
    \FPdim_{{}^K_K\YD(\B)}(A) = \frac{\FPdim(\C)}{\FPdim(\B)} = \FPdim_{\B}(K).
\end{equation*}
Suppose for a while that the base field is of characteristic zero.
Let $\mathfrak{g}$ be a simple Lie algebra, and let $q$ be a root of unity of order $N > 1$.
We assume that $N$ is odd and, moreover, $N$ is not divisible by $3$ when $\mathfrak{g}$ is of type $G_2$.
Then the category $\Rep(u_q(\mathfrak{g}))$ is equivalent to ${}^K_K\YD(\B)$ for some $K$ and $\B$ \cite[Example 5.17]{laugwitz2021braided}. The above discussion gives a non-trivial exact commutative algebra in $\Rep(u_q(\mathfrak{g}))$. 

We next provide an explicit description of the algebra $A$ in the case of $\Rep(u_q(\mathfrak{sl}_2))$.


\subsubsection{Example: Small quantum $\fsl_2$ at odd order root of unity}\label{subsec:uqsl2}
Suppose that $\kk=\mathbb{C}$. Let $q = e^{2\pi \iota/N}$ be the root of unity of odd order $N=2h+1$. Then, $u_q(\fsl_2)$ is a Hopf algebra with generators $E,F,K$ and subject to certain relations (see Appendix~\ref{app:odd-order-uqsl2} for details).
The category $\C=\Rep(u_q(\fsl_2))$ is a non-degenerate finite ribbon category \cite{lyubashenko1995invariants}. In this section, we consider the algebra $A=\kk[x]/(x^N)$. It admits the following $u_q(\fsl_2)$-action making it an algebra object in $\Rep(u_q(\fsl_2))$
\begin{equation*}
  K\cdot x = q^{-2} x, \quad E\cdot x =  1, \quad F\cdot x = -q x^2 .
\end{equation*}

\begin{lemma}
$A$ is a haploid, commutative, exact algebra in $\Rep(u_q(\fsl_2))$ with $\theta_A=\id_A$ and $A^*\notin (\C_A^{\loc})'$.
\end{lemma}
\begin{proof}
Observe that $\Hom_{\C}(\unit,A) = \{ f:\kk\rightarrow A| \varepsilon(h)f(1)= h \cdot f(1) \;\forall h\in u_q(\fsl_2) \}$. This implies that $f(1)\in\kk$. Hence, $\Hom_{\C}(\unit,A)\cong \kk$, so $A$ is haploid. 

Next, we prove that $A$ is simple, that is, it does not contain any two-sided ideal in $\C$. Suppose there exists a non-trivial ideal $I$. Consider an element $f\in I$. Then, as a polynomial in $x$, suppose that the leading term of $f$ is $ax^n$ for some scalar $a$ and $n<N$. Consider the element $y=E^{\deg(f)}\cdot f\in I$. As acting by $E$ reduces the degree by $1$, $y$ will be a scalar.
In fact, $y$ will be a nonzero scalar. Thus, $I=A$, thereby proving that $A$ is simple. As $A$ is simple and indecomposable, by \cite[Proposition~B.7]{etingof2021frobenius}, it follows that $A$ is exact.

Next we show that $A$ is commutative. From \cite{laugwitz2022relative}, recall that we can realize $\Rep(u_q(\fsl_2))$ as the category ${}^H_H\YD(\D)$ where $H=\kk[x]/(x^N)$ and $\D=\Rep(\mathbb{Z}_N)$ with a chosen braiding. 
Additionally, by \cite[Example~5.18]{laugwitz2023constructing}, this algebra $H\in{}^H_H\YD(\D)$ is commutative. Thus, $A$ is a commutative algebra in $\Rep(u_q(\fsl_2))$ as well.

For a proof that $\theta_A=\id_A$, see Lemma~\ref{lem:theta-id-odd-order-sl2}. Lastly, we show that $A^*\notin (\C_A^{\loc})'$. As $\C$ is nondegenerate, if $A^*\in (\C_A^{\loc})'$, then $A$ would be isomorphic to $A^*$, thereby making it a Frobenius algebra. In particular, it would admit a nonzero map $\varepsilon:A\rightarrow\unit$. 
But if $f\in\Hom_{\C}(A,\unit)$ then $f(1) = f(E\cdot x) = \varepsilon(E)f(x) = 0$. In a similar manner, one checks that $f(x^i)=0$ for all $i$. Thus, we have obtained a contradiction. Hence, $A^*\notin (\C_A^{\loc})'$.
\end{proof}


\section{Witt equivalence and VOA extensions}\label{sec:Witt-equivalence}
In this section, we discuss the Witt equivalence of braided tensor categories and its relation to extensions of vertex operator algebras (VOAs). The ground field is assumed to be algebraically closed.


\subsection{Definitions and basic results}\label{subsec:Witt-definition}
We will frequently use the fact that for tensor categories $\C$ and $\D$, we have a tensor equivalence $\C\boxtimes\D\simeq \D\boxtimes\C$ given by $X \boxtimes Y \mapsto Y \boxtimes X$. If, furthermore, $\C$ and $\D$ are braided, it is a braided tensor equivalence. Also, we use the notation $R_{\C}:\C\to\Z(\C)$ to denote the right adjoint to the forgetful functor $\Z(\C)\to\C$.

\begin{theorem}\label{thm:Witt-equivalence}
Let $\C$ and $\D$ be non-degenerate braided finite tensor categories. Then, the following are equivalent:
\begin{enumerate}
    \item[\textup{(a)}] There exist finite tensor categories $\cA$ and $\B$ such that $\C\boxtimes\Z(\cA) \simeq \D\boxtimes\Z(\B)$ as braided tensor categories.
    \item[\textup{(b)}] There exists a non-degenerate braided finite tensor category $\E$ with haploid commutative exact algebras $A_1$ and $A_2$ such that $\E_{A_1}^{\loc}\simeq\C$ and $\E_{A_2}^{\loc}\simeq\D$.
    \item[\textup{(c)}] There exist nondegenerate braided finite tensor categories $\D_i$, $1\leq i\leq n$, and finite tensor categories $\cA_i$, $0\leq i\leq n$, for some $n\in\mathbb{Z}_{\geq 0}$ such that 
    \[ \D_i\boxtimes\oD_{i+1}\simeq \Z(\cA_i), \quad \forall\; 0\leq i\leq n\]
    with $\D_0\simeq \C$ and $\D_{n+1}\simeq \D$.
  \end{enumerate}
\end{theorem}
\begin{proof}
(a) $\Rightarrow$ (b) Take $\E = \C\boxtimes \Z(\cA)$, and the exact algebra (by Remark~\ref{rem:exact-Deligne}) $A_1=\unit_{\C}\boxtimes R_{\cA}(\unit_{\cA})$. Observe that 
\[\E_{A_1}^{\loc}
\simeq (\C\boxtimes \Z(\cA))_{\unit_{\C}\boxtimes R_{\cA}(\unit_{\cA})}^{\loc} 
\simeq \C_{\unit_{\C}}^{\loc}\boxtimes \Z(\cA)_{R_{\cA}(\unit_{\cA})}^{\loc} 
\stackrel{(\heartsuit)}{\simeq} \C\boxtimes\Vect
\simeq \C.\]
The equivalence ($\heartsuit$) uses Remark~\ref{rem:Z(C)-lagrangian}. 
As $\E\simeq \D\boxtimes \Z(\B)$, by taking $A_2=\unit_{\D}\boxtimes R_{\B}(\unit_{\B})$ we get that $\E_{A_2}^{\loc}\simeq \D$ in a similar manner. It is clear that $A_1$ and $A_2$ are haploid as well.

(b) $\Rightarrow$ (c) 
Note the equivalences, 
\[ \C\boxtimes\oE \simeq \oE \boxtimes \E_{A_1}^\loc \stackrel{(\spadesuit)}{\simeq} \Z(\E_{A_1}) 
\quad \text{and} \quad 
\E \boxtimes \oD \simeq \overline{\oE\boxtimes \D} \simeq \overline{\oE\boxtimes \E_{A_2}^\loc} \stackrel{(\spadesuit)}{\simeq} \overline{\Z(\E_{A_2})} \stackrel{(\diamondsuit)}{\simeq} \Z((\E_{A_2})^{\op}) \]
Here the equivalences ($\spadesuit$) use Corollary~\ref{cor:ndBFTC-Z(C-A)}(a) and ($\diamondsuit$) uses \cite[Exercise~8.5.2]{etingof2016tensor}. Since $\E_{A_1}$ and $(\E_{A_2})^{\op}$ are finite tensor categories, the claim follows.

(c) $\Rightarrow$ (a) This follows by the following sequence of braided tensor equivalences:   
\begin{align*}
  \D_0\boxtimes \Z(\D_1\boxtimes \ldots\boxtimes \D_{n+1}) 
  \simeq & \D_0 \boxtimes (\D_1\boxtimes \oD_1) \boxtimes\ldots \boxtimes (\D_{n+1}\boxtimes\oD_{n+1}) \\ 
  \simeq & (\D_0\boxtimes \oD_1) \boxtimes (\D_1\boxtimes \oD_2) \boxtimes \ldots (\D_n\boxtimes\oD_{n+1}) \boxtimes \D_{n+1} \\
  \simeq & \Z(\cA_0) \boxtimes \Z(\cA_1) \boxtimes \ldots \boxtimes \Z(\cA_n) \boxtimes \D_{n+1}\\
  \simeq & \D_{n+1}\boxtimes \Z(\cA_0\boxtimes \ldots \boxtimes \cA_n ).  \qedhere
\end{align*} 
\end{proof}

Following \cite{davydov2013witt} and \cite{laugwitz2023constructing}, we define:
\begin{definition}
If two non-degenerate braided finite tensor categories $\C$ and $\D$ satisfy the equivalent conditions of Theorem~\ref{thm:Witt-equivalence}, we say that $\C$ and $\D$ are \textit{Witt equivalent} and write $\C\sim\D$.
\end{definition}

As in the semisimple case, the above definition yields an equivalence relation on the class of non-degenerate braided FTCs. We denote the Witt equivalence class of a non-degenerate braided finite tensor category $\C$ as $[\C]$. We denote the set of Witt equivalence classes of braided finite tensor categories as $\W^{\ns}$. 
Then, we get the following results immediately.
\begin{lemma}\label{lem:Witt-consequences}
\begin{enumerate}
    \item[\textnormal{(a)}] The set $\W^{\ns}$ is an abelian group.
    \item[\textnormal{(b)}] Let $A$ be a haploid, commutative, exact algebra in $\C$. Then $[\C]=[\C_A^{\loc}]$ in $\W^{\ns}$.
    \item[\textnormal{(c)}] Let $\D$ be a non-degenerate topologizing braided tensor subcategory of some non-degenerate braided FTC $\C$. Then $[\D][\Z_{(2)}(\D\subset \C)] = [\C]$ in $\W^{\ns}$.
    \item[\textnormal{(d)}] For a FTC $\C$ and a nondegenerate braided FTC $\D$ with fully faithful tensor functor $\D\rightarrow\Z(\C)$, we have $[\Z_{(2)}(\D\subset\Z(\C))] = [\oD]$ in $\W^{\ns}$.
\end{enumerate}
\end{lemma}
\begin{proof}
(a) The group operation in $\W^{\ns}$ is given by Deligne product and the identity element is $[\Vect]$. Since $\C\boxtimes\oC\simeq \Z(\C)$, the inverse of $[\C]$ is $[\oC]$.

(b) By Theorem~\ref{cor:ndBFTC-Z(C-A)}(a), $\Z(\C_A)\simeq \oC\boxtimes\C_A^{\loc}$. Now, the claim follows by Theorem~\ref{thm:Witt-equivalence}(c).

(c) This follows by \cite[Theorem~4.17]{laugwitz2022relative} which says that $\C\simeq \D\boxtimes \Z_{(2)}(\D\subset \C)$. 

(d) This is a corollary of (c).
\end{proof}

Next we present an example of a relation in the Witt group.
For this, we first discuss an example from \cite{laugwitz2022relative}. 
We assume that $G$ and $q$ as in the setting of \cite[\S5.2]{laugwitz2022relative}. Namely, $G=\langle g_1,\ldots,g_n\rangle$ is a finite abelian group, and $\bq = (q_{ij})\in\mathrm{Mat}_n(\kk)$ is a family of elements of $\kk^{\times}$ such that $r(g_i,g_j)=q_{ji}$ defines a universal $r$-form on $K=\kk G$. One can define the category $\B_{\bq}=\mathrm{Corep}(K)$ which is a semisimple braided monoidal category with the braiding induced by $r$.
Now consider a vector space $V$ with basis $\{x_1,\ldots,x_n\}$. The following action and coaction of $K$ on $V$ make it a Yetter-Drinfeld module:
\[ g_i\cdot x_j = q_{ij} x_j,\qquad \delta(x_i) = g_i\otimes x_i . \]
Consequently, one can form the Nichols algebra $\fB_{\bq}:=\fB(V)$ which is a Hopf algebra in the category of $K$-Yetter Drinfeld modules. If $\fB_{\bq}$ is finite-dimensional, then one forms the finite tensor category $\C$ of $\fB_{\bq}$-modules in $\B$. In \cite{laugwitz2022relative}, an explicit quasitriangular Hopf algebra $H=\mathrm{Drin}_{K^*}(\fB_{\bq}^*,\fB_{\bq})$ is provided such that 
$\Rep(H) \simeq \Z_{(2)}(\B_{\bq}\subset \Z(\C))$.

The category $\B_{\bq}$ is non-degenerate if and only if the pairing $b$ on $K$ is nondegenerate. In this case, by Lemma~\ref{lem:Witt-consequences} (d), we get the following relation in the Witt group $\W^{\ns}$:
\[ [\Rep(\mathrm{Drin}_{K^*}(\fB_{\bq}^*,\fB_{\bq}))] = [\overline{\B_{\bq}}]. \] 

The class of Hopf algebras $\mathrm{Drin}_{K^*}(\fB_{\bq}^*,\fB_{\bq})$ includes in particular the small quantum groups $u_q(\fg)$ at odd order roots of unity $q$ (see \cite[\S5.4]{laugwitz2022relative}). Thus, by the above discussion, we get that the categories $\Rep(u_q(\fg))$ are Witt equivalent to semisimple categories $\B_{\bq}$ for some $\bq$.


\subsection{Completely anisotropic categories}\label{Witt-anisotropic}
The notion of completely anisotropic braided fusion categories was introduced in \cite{davydov2013witt}. Generalizing it we introduce the following notion.
\begin{definition}\label{def:completely-anisotropic}
We call a braided finite tensor category \textit{completely anisotropic} if the only haploid, commutative, exact algebra in it is $\unit$.
\end{definition}

\begin{remark}
An alternate definition of completely anisotropic categories, using rigid Frobenius algebras, was suggested in \cite{laugwitz2023constructing}. There the algebra is assumed to be separable, which (as explained in Section~\ref{subsec:relation-to-recent-work}) may be restrictive. Note that a completely anisotropic category as in Definition~\ref{def:completely-anisotropic} will be completely anisotropic in the sense of \cite{laugwitz2023constructing}.
\end{remark}

Next we consider some properties of completely anisotropic categories.

\begin{lemma}
A braided tensor subcategory $\C$ of a completely anisotropic category $\D$ is also completely anisotropic. 
\end{lemma}
\begin{proof}
If $\C$ is not completely anisotropic then it will contain a haploid, commutative and exact algebra $A$. By Lemma~\ref{lem:exact-preserve}, $A$ will also be a non-trivial haploid, commutative, exact algebra in $\D$, thereby contradicting our assumption.
\end{proof}

\begin{lemma}\label{lem:aniso-central-functor}
  A braided FTC $\C$ is completely anisotropic if and only if every central functor out of it is fully faithful. 
\end{lemma}
\begin{proof}
($\Rightarrow$)
Let $F:\C\rightarrow\D$ be a central tensor functor.
The restriction $F':\C\rightarrow \mathrm{Im}(F)$ is a surjective central tensor functor. 
By Corollary~\ref{cor:comm-alg-via-adjoint-2}, $A=(F')^{\ra}(\unit_{\D})$ is a haploid, commutative and exact algebra in $\C$ and $\C_A\simeq \mathrm{Im}(F)$. But $\C$ is anisotropic, so $A\simeq \unit_{\C}$ which makes $F':\C\xrightarrow{\sim} \mathrm{Im}(F)$ a tensor equivalence. Thus, $F$ is fully faithful.

($\Leftarrow$)
We will prove the contrapositive.
Suppose that $\C$ is not completely anisotropic. Then we have a haploid, commutative, exact algebra $A\neq \unit$ in $\C$. The free functor $F_A:\C\rightarrow\C_A$ is then a central tensor functor that is not fully faithful, thereby showing that not every central functor is fully faithful.
\end{proof}

Consequently, to show that $\C$ is not completely anisotropic, it suffices to exhibit a central tensor functor $\C\rightarrow\D$ where $\FPdim(\C)>\FPdim(\D)$ (cf. \cite[Proposition~6.3.3]{etingof2016tensor}).

\begin{lemma}\label{lem:criteria-anisotropic}
  Suppose that $\B$ is a braided finite tensor category admitting an injective tensor functor $\B\rightarrow\Z(\C)$ for a finite tensor category $\C$. If $\FPdim(\B)<\FPdim(\C)$, then the relative Drinfeld center $\Z_{(2)}(\B\subset\Z(\C))$ is not completely anisotropic.
\end{lemma}
\begin{proof}
Note that the forgetful functor $U:\Z_{(2)}(\B\subset\Z(\C))\rightarrow\C$ is a central functor. Since 
\[ \FPdim(\Z_{(2)}(\B\subset\Z(\C))) = \frac{\FPdim(\C)^2}{\FPdim(\B)} > \FPdim(\C),\]
$U$ cannot be fully faithful. Thus, by Lemma~\ref{lem:aniso-central-functor}, the claim follows.
\end{proof}

\begin{example}\label{ex:YD-aniso}
From \S\ref{subsec:Yetter-Drinfeld-relative} recall that for a nondegenerate braided finite tensor category $\B$ and a non-trivial Hopf algebra $H\in \B$, we have an equivalence
\[ {}^H_H\YD(\B) \simeq \Z_{(2)}(\B\subset \Z(\C)), \quad \text{where} \quad \C=\Rep(H).\]
As $\FPdim(\C)=\FPdim(\B)\FPdim(H)>\FPdim(\B)$, the category of $H$-Yetter-Drinfeld modules in $\B$ is not completely anisotropic. 
\end{example}

A finite tensor category is called \emph{incompressible} if every tensor functor out of it is fully faithful. By Lemma~\ref{lem:aniso-central-functor}, every braided incompressible category is completely anisotropic.
\begin{example}
Let $\kk$ be of characteristic $p > 0$. Examples of incompressible categories include:
\begin{enumerate}
    \item Certain incompressible symmetric tensor categories $\mathrm{Ver}_p$ constructed in \cite{benson2023new}. 
    \item Non-symmetric analogues of $\mathrm{Ver}_p$ constructed in \cite{decoppet2024verlinde}, generalizing \cite{sutton2023sl}. For odd $p$, certain classes of these categories are incompressible \cite[Proposition~3.53]{decoppet2024verlinde} (We thank Thibault D\'ecoppet for this observation.).
\end{enumerate}
\end{example}

Next we consider completely anisotropic categories that are non-degenerate. For the remainder of this section, we fix char$(\kk)=0$.

\begin{example}
Take an abelian group $G$ with a non-degenerate quadratic form $q$ such that $q(g)\neq 1$ for all $g\neq 1$ in $G$. Then the nondegenerate braided fusion category $\C(G,q)$ is completely anisotropic \cite[\S A.7.1]{drinfeld2010braided}. Other known examples of anisotropic categories come from the categories $\C(\mathfrak{g},l)$ associated to affine Lie algebras at positive integer levels \cite[\S6.4]{davydov2013witt}.
\end{example}

In the next section, we show that some of the currently known examples of non-degenerate BFTCs obtained from vertex operator algebras are not completely anisotropic.


\subsection{VOA extensions}\label{subsec:VOA-extensions}
In this section, we explain the connection between Witt equivalence and VOA extensions. The following discussion can be considered as a non-semisimple (or logarithmic) analogue of \cite[\S6]{davydov2013witt}.

Let $V$ be a vertex operator algebra \cite{frenkel1989vertex}. To each VOA, one can associate a category $\Rep(V)$ of grading-restricted generalized $V$-modules. In a series of papers \cite{HLZ}, Huang, Lepowsky, and Zhang developed a theory determining when $\Rep(V)$ admits a braided monoidal structure. However, in general, verifying the HLZ conditions to prove that $\Rep(V)$ is braided monoidal (and establishing further properties like rigidity or modularity) is difficult.

An inclusion $V\subset W$ of VOAs with coinciding conformal vectors is called a \textit{VOA extension}. By the work of many authors \cite{kirillov2002q,huang2015braided,creutzig2017tensor}, it is known that if $\C=\Rep(V)$ admits the braided monoidal structure of HLZ theory, then $A:=W$ is a haploid commutative algebra in $\C$ (assuming $W$ is of CFT type), and there is a monoidal equivalence $\C_A^{\loc} \simeq \Rep(W)$.

A crucial subtlety concerns the duality structure (see \cite[\S4]{mcrae2025rigidity}). $\Rep(W)$ possesses a natural duality given by the contragredient dual of modules. Conversely, $\C_A^{\loc}$ inherits a categorical duality from $\C$ via the internal Hom. While these structures are a priori distinct, the equivalence $\C_A^{\loc} \simeq \Rep(W)$ identifies them under suitable conditions. Therefore, verifying the rigidity of $\C_A^{\loc}$ (e.g., by checking that $A$ is an exact algebra) implies that $\Rep(W)$ is rigid with respect to standard contragredient duals.

We focus on the setting where $V$ is a \textit{strongly finite} VOA \cite{creutzig2017logarithmic} (i.e., $\mathbb{N}$-graded, simple, self-contragredient, and $C_2$-cofinite). In this case, $\Rep(V)$ is braided monoidal and finite abelian \cite{huang2009cofiniteness}. McRae \cite{mcrae2021rationality} proved that if such a $\Rep(V)$ is rigid, it is necessarily a MTC.

Next we discuss the connection between Witt equivalence and VOA extensions.
\begin{lemma}
  Let $V$ be a strongly finite VOA such that $\Rep(V)$ is rigid. If $V\subset W$ is a VOA extension with $W$ being a simple VOA, then $[\Rep(V)]=[\Rep(W)]$ in the Witt group.
\end{lemma}
\begin{proof}
  First, we note that by \cite{mcrae2021rationality}, $\Rep(V)$ is an MTC. Moreover, $W$ being a simple VOA implies that $A:=W$ is a simple (commutative) algebra in $\Rep(V)$. Thus, by \cite[Theorem~7.1]{coulembier2025simple}, simplicity of $A$ implies that it is exact. Now, thanks to Corollary~\ref{cor:ndBFTC-Z(C-A)}(b), $\Rep(W)\simeq \Rep(V)_A^{\loc}$ is a nondegenerate braided FTC. Finally, the claim follows by Lemma~\ref{lem:Witt-consequences}(b).
\end{proof}

This generalizes the same observation for strongly rational VOAs that was made in \cite[\S6.1]{davydov2013witt}. Using this lemma, one can obtain examples of relations in the Witt group. For example, the following are examples of non-semisimple braided FTCs whose Witt class contains a semisimple representative.

\begin{example}
\begin{enumerate}
  \item Set $q=e^{\pi \iota/p}$ for some $p\in\mathbb{Z}_{>0}$, where $\iota=\sqrt{-1}$, and let $u_q^{\phi}(\fg)$ denote the quasi-Hopf modification of the small quantum group introduced in \cite{creutzig2020quasi,negron2021log}. Let $\W(p)$ denote the triplet VOA \cite{adamovic2008triplet}. 
  By recent works \cite{gannon2021quantum,creutzig2023algebraic}, it is now known that $\Rep(u_q^{\phi}\fsl_2) \simeq \Rep(\W(p))$ as braided tensor categories. As the triplet VOA $\W(p)$ is a strongly finite VOA which admits a VOA extension $\W(p)\subset V$ (where $V$ is the lattice VOA) and $\Rep(V)\simeq \C(\bZ_{2p},q)$ for some quadratic form $q$ on $\bZ_{2p}$, we get that $[\Rep(u_q^{\phi}(\fsl_2))] = [\C(\bZ_{2p},q)]$.

  A similar Witt relation can also be obtained as follows. In \cite{gainutdinov2018modularization}, certain quasi-quantum groups $\tilde{u}_q^{\phi}(\fg)$ are constructed (for $q=e^{\pi \iota/p}$) and it is shown that $\Rep(\tilde{u}_q^{\phi}(\fg)) \simeq \Z_{(2)}(\C\subset \Z(\D))$ where $\C$ is the braided fusion category of $\Lambda$-graded vector spaces (for some finite abelian group $\Lambda$) with non-trivial associator. Negron \cite[\S10.3]{negron2021log} explained that $\Rep(u_q^{\phi}(\fg))\simeq \Rep(\tilde{u}_q^{\phi}(\fg))$ as braided tensor categories. Thus, using Lemma~\ref{lem:Witt-consequences}(d), the above observations yield that $[\Rep(u_q^{\phi}(\fsl_2))] = [\oC]$. 
  
  \item Let $SF^+_d$ denote the even part of the symplectic fermion VOA \cite{abe2007orbifold}. Its representation category $\Rep(SF^+_d)$ is a non-degenerate braided tensor category. 
  In fact, we have an extension of VOAs $\W(2)^{\otimes d} \subset SF^+_d$. Thus, we get that
  \[ [\Rep(SF^+_d)] = [\Rep(\W(2)^{\otimes d})] = [\Rep(\W(2))^{\boxtimes d}] = [\Rep(\W(2))]^d = [\Rep(u_\iota^{\phi}(\fsl_2))]^d \]
  Here the second equality uses \cite[Corollary~1.3]{mcrae2023deligne}.
  As explained in (a), $[\Rep(u_{\iota}^{\phi}(\fsl_2))]$ contains a semisimple representative. This implies that $[\Rep(SF^+_d)]$ contains a semisimple representative for all $d$.

  In \cite{farsad2022symplectic} a factorizable quasi-Hopf algebra $Q(d)$ is introduced for every $d\in\mathbb{Z}_{\geq 0}$. Runkel conjectured that, $\Rep(Q(d))$ is braided tensor equivalent to $\Rep(SF^+_d)$. 
  Thus, if Runkel's conjecture is true, then the Witt class of $[Q(d)]$ contains a semisimple representative for all $d$.
\end{enumerate} 
\end{example}

\begin{lemma}\label{lem:VOA-aniso}
The categories $\Rep(\W(p))$ and $\Rep(SF^+_d)$ are not completely anisotropic.   
\end{lemma}
\begin{proof}
By \cite[Corollary~10.6]{creutzig2023algebraic}, $\Rep(\W(p))$ is equivalent to a category of Yetter-Drinfeld modules. Using this and the VOA extension $\W(2)^{\otimes d} \subset SF^+_d$, one can also write $\Rep(SF^+_d)$ as a category of Yetter-Drinfeld modules (see \cite[Theorem~10.2]{creutzig2023algebraic}).
Thus, by Example~\ref{ex:YD-aniso}, $\Rep(\W(p))$ and $\Rep(SF^+_d)$ are not completely anisotropic. 
\end{proof}

A non-degenerate braided FTC $\C$ is called \textit{prime} if every topologizing non-degenerate braided tensor subcategory is equivalent to either itself or $\Vect$. 
Our next result answers \cite[Question~6.25(1)]{laugwitz2023constructing}.
\begin{lemma}
The category $\Rep(Q(1))$ is a non-semisimple prime modular tensor category that is not completely anisotropic.
\end{lemma}
\begin{proof}
The equivalence $\Rep(Q(1)) \simeq \Rep(SF_1^+)$ was proved in \cite{gainutdinov2017symplectic}. Thus, from Lemma~\ref{lem:VOA-aniso}, it follows that $\Rep(Q(1))$ is not completely anisotropic. On the other hand, by \cite[Proposition~5.3]{berger2023non} it is a prime modular tensor category.
\end{proof}


\subsection{Questions and remarks}\label{subsec:Witt-questions}
We conclude this article by outlining several questions and remarks arising from our results. Some of these are related to questions raised in \cite[\S4]{brochier2021invertible}.


\subsubsection{} \label{subsubsec:WittQ1}
In the semisimple setting, the conditions in Theorem~\ref{thm:Witt-equivalence} are equivalent to the existence of haploid commutative exact algebras $A\in\C$ and $B\in\D$ such that there is an equivalence of local module categories $\C_A^\loc\simeq\D_B^\loc$ \cite[Proposition~5.15]{davydov2013witt}. In non-semisimple setting, the existence of such algebras implies that $\C$ and $\D$ are Witt equivalent, we don't know if the converse is true.

\begin{question}\label{que:WittQ1}
    Is every Witt equivalence between finite tensor categories $\C$ and $\D$ induced by an equivalence $\C_A^\loc \simeq \D_B^\loc$ for some haploid commutative exact algebras $A$ and $B$?
\end{question}


\subsubsection{} \label{subsubsec:WittSur}
Let $\W$ denote the Witt group of nondegenerate braided fusion categories defined in \cite{davydov2013witt}. There exists a natural group homomorphism $\W\to \W^{\ns}$.
\begin{question}\label{que:Witt-sur}
   Is the homomorphism $\W\to \W^{\ns}$ surjective?
\end{question}
Phrased differently, this question asks whether there exists a Witt class in $\W^\ns$ that does not contain a semisimple representative. This question has been raised recently in several works, see e.g., \cite{brochier2021invertible,laugwitz2023constructing}.

While the surjectivity of the map $\W\to\W^\ns$ appears difficult to address, the injectivity is equivalent to the following problem:
\begin{question}\label{que:Witt-trivial}
  If $\C \in [\Vect]$ in the non-semisimple Witt group, is $\C$ necessarily the Drinfeld center of a finite tensor category? 
\end{question}


\subsubsection{} 
Throughout this section, we analyzed various non-semisimple BFTCs whose Witt classes admit semisimple representatives. We subsequently demonstrated that these categories fail to be completely anisotropic. This observation motivates the following:
\begin{question}
If the Witt class of a non-semisimple braided FTC $\C$ contains a semisimple category, can $\C$ be completely anisotropic?
\end{question}
Current evidence suggests the answer is no. Another related question is:
\begin{question}\label{que:Witt-aniso}
Does there exist a non-semisimple completely anisotropic tensor category? 
\end{question}


\subsubsection{} 
Given non-degenerate braided finite tensor categories $\C$ and $\D$, we write $\C \mathrel{\dot{\sim}} \D$ if $\C \boxtimes \overline{\D} \cong \Z(\cA)$ for some finite tensor category $\cA$. Condition (c) in Theorem~\ref{thm:Witt-equivalence} corresponds to the equivalence relation generated by $\mathrel{\dot{\sim}}$. In the semisimple setting, the relation $\mathrel{\dot{\sim}}$ is inherently transitive and thus forms an equivalence relation. 

\begin{question}\label{que:Witt-condition-C}
    In the non-semisimple setting, is $\dot{\sim}$ an equivalence relation?
\end{question}


\subsubsection{Relationship between the posed questions}
The questions raised in this section are deeply interconnected, and a resolution to some would provide a pathway to solving the others. We discuss some of these connections below.

If the answer to Question~\ref{que:Witt-condition-C} is yes, there will be many consequences: 
\begin{enumerate}
  \item A similar argument as in \cite[Theorem~5.13]{davydov2013witt} could be used to show that each Witt class contains a unique minimal representative. This could then be used to answer Question~\ref{que:WittQ1} in the affirmative. This would imply that each Witt class contains a unique minimal representative and this representative will be completely anisotropic. 
  \item As a consequence of (a), Questions~\ref{que:Witt-sur} and \ref{que:Witt-aniso} will be closely related: a negative answer to Question~\ref{que:Witt-aniso} would yield a positive answer to Question~\ref{que:Witt-sur}.
  \item This would also imply a positive answer to Question~\ref{que:Witt-trivial}.
\end{enumerate}

\begin{remark}
  Since this paper was posted, all of the questions raised in this subsection have been answered in \cite{ostrik2026non}.
\end{remark}


\section*{Acknowledgements}
The second author (H.Y.) would like to thank Chelsea Walton for providing financial support that made this collaboration happen. H.Y. would also like to thank Terry Gannon, Thomas Creutzig, Thibault D\'ecoppet, Victor Ostrik and Chelsea Walton for helpful discussions. Lastly, we are grateful to the referee for the careful reading of our manuscript. 

The first author (K.S.) is supported by JSPS KAKENHI Grant Number 24K06676.
H.Y. is supported by a start-up grant from the University of Alberta and an NSERC Discovery Grant. Part of the work was finished while H.Y. was in residence at the Mathematical Sciences Research Institute in Berkeley, California, during the Quantum symmetries reunion in 2024; this was supported by NSF grant DMS-1928930.


\appendix


\section{Small quantum sl2: odd order case}\label{app:odd-order-uqsl2}
We will work over an algebraically closed field $\kk$ of characteristic $\text{char}(\kk)=0$.
Let $q=e^{2\pi\iota/N}$ be a primitive root of unity of order $N=2h+1$, where $h \ge 1$. Then, $u_q(\fsl_2)$ is a Hopf algebra with generators $E,F,K$ and subject to the following relations:
\begin{equation*}
K^N=1,\quad E^N=F^N=0, \quad KE = q^2 EK, \quad KF = q^{-2} FK, \quad EF - FE = \frac{K-K^{-1}}{q-q^{-1}}. 
\end{equation*}
The Hopf algebra structure is given by:
\begin{align*}
\Delta(K) = K\otimes K, \; \quad  \Delta(E) = E\otimes K + 1\otimes E, \quad \Delta(F) = F\otimes 1 + K^{-1}\otimes F, \;\qquad \\
\varepsilon(K) = 1 , \quad \varepsilon(E) = \varepsilon(F) = 0, \quad  S(K) = K^{-1}, \quad S(E) = -EK^{-1}, \quad S(F) = -KF. 
\end{align*}
Recall that the $R$-matrix of $u_q(\fsl_2)$ is (see \cite[\S9.7]{kassel2012quantum}):
\begin{equation*}
R = R_{(1)}\otimes R^{(2)} = \frac{1}{N} \sum_{0\leq p,s,r \leq N-1} \frac{(q-q^{-1})^p}{[p]!} q^{\frac{p(p-1)}{2} + 2p(s-r)-2sr} E^p K^s\otimes F^p K^r.
\end{equation*}
The ribbon element is given by:
\begin{equation*}
  \theta = \beta\left(\sum_{r=0}^{N-1} q^{h r^2}\right)\left(\sum_{0 \leq m, j \leq N-1} \frac{\left(q-q^{-1}\right)^m}{[m]!}(-1)^m q^{-\frac{1}{2} m+m j+\frac{1}{2}(j+1)^2} F^m E^m K^j\right).
\end{equation*}
where $\beta = \frac{(-\iota)^h}{\sqrt{N}(1+2 q^{1/2})}$. Note that this ribbon element is a rescaling of the typically used ribbon element by a scalar factor.
This makes the category $\C=\Rep(u_q(\fsl_2))$ a non-degenerate finite ribbon category \cite{lyubashenko1995invariants}.

In this section and the next, we will frequently use the following elementary lemma.
\begin{lemma}\label{lem:sum-roots-of-unity}
  Let $n>1$ and $q$ a primitive $n^{\text{th}}$ root of unity. Then $\sum_{j=0}^{n-1} q^{rj}$ equals $n$ if $r\equiv 0 \bmod n$ and $0$ otherwise. 
\end{lemma}
\begin{proof}
  If $r\equiv 0 \bmod n$, then $q^r=1$ and hence the sum is $n$. Otherwise, the sum is given by the geometric series formula $\frac{1-q^{rn}}{1-q^r}=0$. As $q^{rn}=(q^n)^r=1$, the claim follows.
\end{proof}


\subsection*{An algebra object in \texorpdfstring{$\Rep(u_q(\fsl_2))$}{Rep(uqsl2)}}
Next we consider the algebra $A=\kk[x]/(x^N)$. It admits the following $u_q(\fsl_2)$-action making it an algebra object in $\Rep(u_q(\fsl_2))$
\begin{equation*}
  K\cdot x = q^{-2} x, \quad E\cdot x =  1, \quad F\cdot x = -q x^2 .
\end{equation*}


\begin{lemma}\label{lem:theta-id-odd-order-sl2}
  The algebra $A$ satisfies $\theta_A=\id_A$.
\end{lemma}
\begin{proof}
We show below that $\theta$ acts on $x$ as identity. Then, as $\theta_A:A\rightarrow A$ is an algebra map, $\theta\cdot x^i = \theta_A(x^i) = \theta_A(x)^i =(\theta\cdot x)^i =x^i$, for all $i$. 
\begin{align*}
  \theta\cdot x = \;
  & \beta\left(\sum_{r=0}^{N-1} q^{h r^2}\right)
  \left(\sum_{0 \leq m, j \leq N-1} \frac{\left(q-q^{-1}\right)^m}{[m]!}(-1)^m q^{-\frac{1}{2} m+m j+\frac{1}{2}(j+1)^2} F^m E^m K^j\right)\cdot x \\
  = \; & \beta\left(\sum_{r=0}^{N-1} q^{h r^2}\right) \left(\sum_{0 \leq m, j \leq N-1} \frac{\left(q-q^{-1}\right)^m}{[m]!}(-1)^m q^{-\frac{1}{2} m+m j+\frac{1}{2}(j+1)^2} F^m E^m \cdot (q^{-2j}x) \right) \\
  \stackrel{(*)}{=} \; & \beta\left(\sum_{r=0}^{N-1} q^{h r^2}\right) \left(\sum_{j=0}^{N-1} 
   q^{\frac{(j+1)^2}{2} - 2j} \right) x \;
  = \;  \beta\left(\sum_{r=0}^{N-1} q^{h r^2}\right) 
  \left(\sum_{j=0}^{N-1} q^{\frac{(j-1)^2}{2}} \right) x.
\end{align*}
The equality $(*)$ is obtained by observing that $F^m E^m\cdot x$ is nonzero if and only if $m=0$. 

The first sum $(\sum_{r=0}^{N-1} q^{h r^2})$ is given by standard Gauss sum formula and is $\sqrt{N}(\iota)^h$.
After cancellation, the second sum $(\sum_{j=0}^{N-1} q^{\frac{(j-1)^2}{2}})$ reduces to $(1+2q^{1/2})$. Thus, the action of $\theta$ on $x$ is scalar multiplication by $\beta \sqrt{N}(\iota)^h(1+2q^{1/2}) = (-\iota^2)^h = 1^h =1$. 
\end{proof}


\section{Small quantum sl2: even order case}\label{app:even-order-uqsl2}
We will work over an algebraically closed field $\kk$ of characteristic $\text{char}(\kk)=0$.
Let $q = e^{\pi \iota/p}$ and set $N=2p$ for some positive integer $p>1$. Then the small quantum group $u_q(\fsl_2)$ is a Hopf algebra that does not admit a quasitriangular structure. Thus, one has to modify the Hopf algebra structure to obtain a quasitriangular quasi-Hopf algebra $u_q^{\phi}(\fsl_2)$ as explained in \cite{creutzig2020quasi}. 
The algebra structure is given by:
\begin{equation*}
K^{2p}=1,\quad E^p=F^p=0, \quad KE = q^2 EK, \quad KF = q^{-2} FK, \quad EF - FE = \frac{K-K^{-1}}{q-q^{-1}}. 
\end{equation*}
The co-multiplication is modified as follows:\footnote{This quasi-Hopf structure corresponds to the choice $t=1$ in the setting of \cite[\S4.2]{creutzig2020quasi}.}
\[  
\Delta(E) = E\otimes K + (e_0 + q e_1) \otimes E, \quad \Delta(F) = F\otimes 1 + (e_0 + q^{-1}e_1)K^{-1}\otimes F, \quad \Delta(K) = K\otimes K,
\]
where $e_0 = (1+K^p)/2$ and $e_1 = (1-K^p)/2$. The co-associator is given by:
\[ \Phi = 1\otimes 1\otimes 1 + e_1 \otimes e_1 \otimes (K^{-1} - 1).\]
The antipode is given by:
\[ S(K) = K^{-1}, \quad S(E) = EK^{-1}(e_0+qe_1), \quad S(F) = -KF(e_0 + q^{-1}e_1)\]
with (co)evaluation elements $\alpha = 1$ and $\beta = e_0 + K^{-1}e_1$. The $R$-matrix is given by:
\[ 
R = \frac{1}{4p}\sum_{n=0}^{p-1}\sum_{s,r=0}^{2p-1} \frac{(q-q^{-1})^n}{[n]!} q^{\frac{n(n-1)}{2} -2sr} (1+q^r+q^{-(n+s)} + q^{1/2 + r - n-s}) K^s E^n \otimes K^r F^n.
\]
Lastly, it admits a twist element given by:
\[ 
\theta = \frac{1-\iota}{2\sqrt{p}} \sum_{n=0}^{p-1} \sum_{j\in\mathbb{Z}_{2p}} \frac{(q-q^{-1})^n}{[n]!} q^{n(j-1/2)+1/2(j+p+1)^2} F^nE^nK^j.
\]
With the above data, the category $\C=\Rep(u_q^{\phi}(\fsl_2))$ is a modular tensor category \cite[Theorem~4.1]{creutzig2020quasi}.


The unit object is $\unit=\kk$ with the actions of $E,F,K$ given by $0,0,1$ respectively. In addition, $\C$ admits another simple invertible object $\psi$ which is $\kk$ as a vector and actions of $E,F,K$ are given by $0,0,-1$ respectively. Next, we note some properties of $\psi$. We will use that $q^p = e^{\pi\iota} = -1$ and $q^{p/2} = e^{\pi\iota/2} = \iota$.


\subsection{\texorpdfstring{$\psi\otimes\psi\cong \unit$}{Psi-tensor-Psi}}
As $\Delta(K)=K\otimes K$, it acts on $\psi\otimes \psi$ by $(-1)(-1) = 1$. As $\Delta(E)$ and $\Delta(F)$ both have either a $E$ or $F$ term in each summand, they both act on $\psi\otimes\psi$ by $0$.


\subsection{\texorpdfstring{$c_{\psi,\psi}=(\iota)^{p}\id_{\unit}$}{braiding-psi}}
Using the formula of $R$-matrix, clearly $c_{\psi,\psi} = \zeta \id_{\unit}$ where
\begin{align*}
  \zeta & 
   =  \frac{1}{4p}\sum_{n=0}^{p-1}\sum_{s,r=0}^{2p-1} \frac{(q-q^{-1})^n}{[n]!} q^{\frac{n(n-1)}{2} -2sr} (1+q^r+q^{-(n+s)} + q^{1/2 + r - n-s}) [(K^r F^n)\cdot\psi \, (K^s E^n)\cdot\psi]\\
  & = \frac{1}{4p}\sum_{s,r=0}^{2p-1} q^{-2sr} (1+q^r+q^{-s} + q^{1/2 + r -s}) (-1)^{r+s}\\
  = & \, \frac{1}{4p} \left[ \sum_{s,r=0}^{2p-1} q^{-2sr+pr+ps} + \sum_{s,r=0}^{2p-1} q^{-2sr+r+pr+ps} + \sum_{s,r=0}^{2p-1} q^{-2sr-s+pr+ps} + \sum_{s,r=0}^{2p-1} q^{1/2 + r -s -2sr+pr+ps}  \right].
\end{align*}
The last equality is obtained by replacing $(-1)^{r+s}$ by $q^{pr+ps}$.
Each of the $4$ sums can be computed by separating the $s$ and $r$ parts and using Lemma~\ref{lem:sum-roots-of-unity}. For instance, the second term is $\sum_{s=0}^{2p-1} q^{ps} \left(\sum_{r=0}^{2p-1} q^{r(-2s+1+p)}\right)$. When $p$ is even, $-2s+1+p$ is never congruent to $0$ mod $2p$ making the sum zero. When $p$ is odd, the inner sum equals $2p$ when $s=(p+1)/2$ or $s=p+ (p+1)/2$ and $0$ otherwise. Thus, the total is 
\[ q^{p(p+1)/2} +  q^{p(p+1)/2+p^2} = q^{p(p+1)/2}(1+ (q^p)^p) = q^{p(p+1)/2} (1 + (-1)^p) = 0 \] 
where the last equality holds because $p$ is odd. By a similar argument, the third sum is zero.

The first term equals $\sum_{s=0}^{2p-1} q^{ps} \left(\sum_{r=0}^{2p-1} q^{r(p-2s)}\right)$. Using Lemma~\ref{lem:sum-roots-of-unity}, the inner sum is zero if $p$ is odd. If $p$ is even and $s$ is $p/2$ or $3p/2$, the inner sum equals $2p$. Thus, the first term is  
\[ 2p((q^p)^{p/2}+ (q^p)^{3p/2}) =  2p \, q^{p^2/2}(1+ (-1)^{p}) = 4p \, q^{p^2/2} = 4p (q^{p/2})^p = 4p \, \iota^p.\]

The fourth term is $q^{1/2} \sum_{s=0}^{2p-1} q^{s(p-1)} \sum_{r=0}^{2p-1} q^{r(1-2s+p)}$ which is zero if $p$ is even. When $p$ is odd, the inner sum is non-zero only if $s=(p+1)/2$ or $s=p+(p+1)/2$. Thus, it equals 
\[(2p) q^{1/2} q^{(p-1)(p+1)/2} (1 + q^{(p-1)p}) = (2p) q^{p^2/2} (1 + q^{p(p-1)}) = (2p)(q^{p/2})^p (1+1) = 4p \, \iota^p.\]

Thus, depending on whether $p$ is even or odd, only one of the first and fourth terms is non-zero. Also, the second and third terms are always zero. Thus, to conclude, we get that $\zeta = (\iota)^p$.


\subsection{\texorpdfstring{$\theta_{\psi} = -(\iota)^p \, \id_{\psi}$}{twist-psi}} 
Our proof is along the lines of \cite{MSE}. 
Using the formula of $\theta$, the actions of $K,E,F$ on $\psi$ and $q=e^{\frac{2\pi\iota}{2p}}$, we get $\theta_{\psi} = \xi\id_{\psi}$ where
\[\xi = \frac{1-\iota}{2\sqrt{p}} \sum_{j=0}^{2p-1} e^{\frac{2\pi\iota}{4p}(j+p+1)^2} (-1)^j .\]

As $(j+p+1)^2 \equiv (j+p+1+2p)^2 \bmod 4p$, we have that 
\[e^{\frac{2\pi\iota}{4p}(j+p+1)^2} = e^{\frac{2\pi\iota}{4p}(j+p+1+2p)^2}.\]
Thus, 
\[\sum_{j=0}^{2p-1} e^{\frac{2\pi\iota}{4p}(j+p+1)^2} (-1)^j = \sum_{j=0}^{2p-1} e^{\frac{2\pi\iota}{4p}(j+p+1+2p)^2} (-1)^{j+2p} = \sum_{j=2p}^{4p-1} e^{\frac{2\pi\iota}{4p}(j+p+1)^2} (-1)^j. \]
The last equality is obtained by changing variable from $j+2p$ to $j$. Consequently, 
\begin{align*}
  \xi 
  = \frac{1-\iota}{2\sqrt{p}} \left( \frac{1}{2} \sum_{j=0}^{4p-1} e^{\frac{2\pi\iota}{4p}(j+p+1)^2} (-1)^j \right) 
  = \frac{1-\iota}{4\sqrt{p}} \left( \sum_{k=0}^{4p-1} e^{\frac{2\pi\iota}{4p}k^2} (-1)^{k-p-1} \right)
\end{align*}
where the second equality is a change of variables $k=j+p+1$. As $k$ runs through $\mathbb{Z}/4p \mathbb{Z}$, so does $k+p+1$, so the indices remain the same. Next, we replace $(-1)^k$ by $e^{\frac{2\pi\iota}{4p}2pk}$ to get
\begin{align*}
  \xi 
  = \frac{(1-\iota)(-1)^{p+1}}{4\sqrt{p}} \sum_{k=0}^{4p-1} e^{\frac{2\pi \iota}{4p}k^2} e^{\frac{2\pi \iota}{4p} 2pk}
  = \frac{(1-\iota)(-1)^{p+1} e^{\frac{2\pi \iota}{4p}(-p^2)}}{4\sqrt{p}} \sum_{k=0}^{4p-1} e^{\frac{2\pi \iota}{4p}(k+p)^2} .
\end{align*}
As $k$ runs through $\mathbb{Z}/4p \mathbb{Z}$, so does $k+p$. So $\sum_{k=0}^{4p-1} e^{\frac{2\pi \iota}{4p}(k+p)^2} = \sum_{k=0}^{4p-1} e^{\frac{2\pi \iota}{4p}k^2}$. Using the quadratic Gauss sum formula, this equals $(1+\iota)\sqrt{4p}$ (see e.g. \cite[Exercise~5, page 43]{berndt1998gauss}). 
To conclude,
\[ \xi =  \frac{(1-\iota)(-1)^{p+1} (e^{\frac{-\pi \iota}{2}})^p}{4\sqrt{p}} (1+\iota)\sqrt{4p} = -(-1)^p(-\iota)^p = -(\iota)^p. \]


\section{Pointed braided fusion categories}\label{app:pointed}
In this section, we recall some known facts about pointed braided fusion categories. While these are well-known when char$(\kk)=0$ \cite{joyal1993braided,drinfeld2010braided}, we discuss the case char$(\kk)=p>0$ for completeness. 

Let $G$ be a finite abelian group. A \textit{quadratic form} on $G$ is a function $q:G\rightarrow \kk^{\times}$ such that $q(g)=q(g^{-1})$ and 
\[ \frac{q(g_1g_2h)}{q(g_1g_2)} = \frac{q(g_1h)q(g_2h)}{q(g_1)q(g_2)q(h)} . \] 
This implies $q(g^n)=q(g)^{n^2}$.
One can use this data to construct a braided fusion category $\C(G,q)$ whose objects are given by $\{\delta_g\}_{g\in G}$ and $\Hom(\delta_g,\delta_h)=\delta_{g,h}\id_g$. The category $\C(G,q)$ is braided monoidal with $\delta_g\otimes \delta_h=\delta_{gh}$ and unit object $\unit=\delta_e$. 
By \cite[Theorem~8.4.9]{etingof2016tensor}, $q$ determines an abelian $3$-cocycle, namely a pair $(\omega:G^{\times 3}\rightarrow \kk^{\times}, c: G^{\times 2}\rightarrow\kk^{\times})$ satisfying certain conditions, on $G$. In fact, the associator and braiding of $\C(G,q)$ are given by $\omega$ and $c$ respectively. 

Let $H\subset G$ be a subgroup satisfying $q|_H=1$. Then there exists a commutative exact Frobenius algebra $A_H$ in $\C(G,q)$.
\begin{itemize}
  \item Since $q|_H=1$, the restriction of the abelian $3$-cocycle $(\omega,c)$ to $H$ is trivial in abelian cohomology. Thus, one can choose a $2$-cochain $\psi: H^{\times 2}\rightarrow \kk^{\times}$ with $d\psi=\omega|_H^{-1}$ such that the algebra
  \[
    A_H=\oplus_{h\in H} \delta_h
  \]
  with multiplication given by $h\otimes h'\mapsto \psi(h,h')hh'$ is haploid and commutative. The corresponding $\C(G,q)$-module category $\M=\C(G,q)_{A_H}$ is semisimple. Since $\C(G,q)$ is fusion, $\M$ is exact. Therefore, $A_H$ is exact. Different choices of $\psi$ give isomorphic algebras.
  \item Moreover, we can construct the following maps
  \begin{align*}  
    \Delta: A_H \rightarrow A_H\otimes A_H, & \quad \delta_h\mapsto \oplus_{k\in H} \frac{\delta_{hk}\otimes \delta_{k^{-1}}}{\psi(hk,k^{-1})}\\
    \varepsilon:A_H\rightarrow \unit \qquad, & \quad \delta_h\mapsto \delta_{h,e} \delta_e
  \end{align*}
  which make $A_H$ a Frobenius algebra in $\C(G,q)$. 
\end{itemize}
In fact, every commutative exact algebra in $\C(G,q)$ is of the form $A_H$ (up to isomorphism) for some subgroup $H\subset G$ satisfying $q|_H=1$.

Lastly, observe that $m\circ\Delta(\delta_h) = |H| \, \id_{\delta_h}$. Thus, $A_H$ is separable Frobenius if and only if $|H|\neq 0$ in $\kk$. Equivalently, $H$ has no $p$-torsion (see also \cite[Lemma~2.1.8]{decoppet2023relative}). In fact, under this condition, we get a special Frobenius algebra.


\bibliographystyle{alpha}
\bibliography{references}

\end{document}